\documentclass[11pt,twoside]{journal}
\newcommand{\myheader}[4]{ \vspace*{-25mm}
	\noindent 
	\parbox[t]{\textwidth}{
	\parbox[t]{6cm}{
	\vspace*{0cm}
	\fontsize{8}{8}\selectfont
	(\the\year)
	}
	\hfill
	\parbox[t]{5.5cm}{\centering
	\fontsize{11}{12}\selectfont \bf      	     
	     \vskip2.5pt%
\fontsize{8}{8}\selectfont \sl 
      	     \vskip2.5pt%
	}
	}
\vspace{1.0cm}
\\
\addtocounter{page}{#3}
\addtocounter{page}{-1}
}

\newcommand{\mytitle}[1]{\noindent\parbox[t]{\textwidth}{\fontsize{17}{17}\selectfont \bf #1}\vspace{0.5cm}}

 \newcommand{\myauthor}[2]{%
 \noindent\parbox[t]{\textwidth}{\fontsize{11}{11}\selectfont #1}\newline \parbox[t]{13cm}{\fontsize{8}{9.5}\selectfont #2}\vspace{0.4cm}}

\newcommand{\summary}[1]{\vspace{0.4cm}\noindent\parbox[t]{\textwidth}{\fontsize{12}{12}\selectfont Summary} \\[1.5ex] \noindent\parbox[t]{\textwidth}{\fontsize{9}{10}\selectfont #1} \vspace{0.5cm}}

\input epsf
\usepackage{psfig}
\usepackage{latexsym}
\usepackage[dvips]{graphicx}
\usepackage{amsbsy,pstricks}
\usepackage{amsmath,amssymb,amsfonts}
\usepackage{Bold_ams,stmaryrd,psfrag,subfigure}
\usepackage[ruled,section]{algorithm}
\usepackage{algorithmic,url,rotating}
\usepackage{pi,multirow,slashbox}
\newtheorem{remark}{Remark}[section]
\newcommand{\comment}[1]{ }

\begin{document}

\year2007
\received={March 2005}
\myheader{13}{4}{1}{58}
\markboth{P. Gosselet and C. Rey}{Non-overlapping domain decomposition methods \\in structural mechanics}
\thispagestyle{myfirstpage}

\mytitle{Non-overlapping domain decomposition methods in structural mechanics}

\myauthor{Pierre Gosselet}{LMT Cachan\\
ENS Cachan / Universit\'e Pierre et Marie Curie / CNRS UMR 8535 \\
61 Av. Pr\'esident Wilson\\
94235 Cachan FRANCE\\
Email: gosselet@lmt.ens-cachan.fr}

\myauthor{Christian Rey}{LMT Cachan\\
ENS Cachan / Universit\'e Pierre et Marie Curie / CNRS UMR 8535 \\
61 Av. Pr\'esident Wilson\\
94235 Cachan FRANCE\\
Email: rey@lmt.ens-cachan.fr}

\summary{ The modern design of industrial structures leads to very
complex simulations characterized by nonlinearities, high
heterogeneities, tortuous geometries... Whatever the modelization
may be, such an analysis leads to the solution to a family of
large ill-conditioned linear systems. In this paper we study
strategies to efficiently solve to linear system based on
non-overlapping domain decomposition methods. We present a review
of most employed approaches and their strong connections. We
outline their mechanical interpretations as well as the practical
issues when willing to implement and use them. Numerical
properties are illustrated by various assessments from academic to
industrial problems.
An hybrid approach, mainly designed for multifield problems, is
also introduced as it provides a general framework of such
approaches.}

\section{INTRODUCTION}
\begin{figure}[ht]\centering
\begin{minipage}{.4\textwidth}\centering \includegraphics[width=0.8\textwidth]{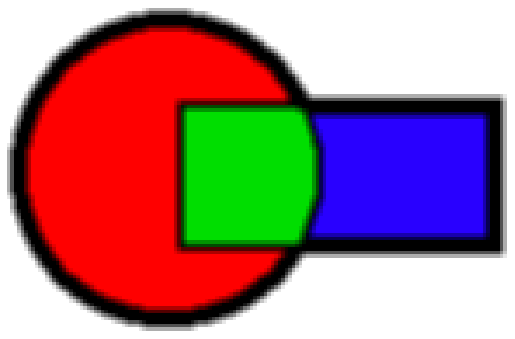}\caption{Schwarz' original
problem}\label{fig:ddm}
\end{minipage}
\begin{minipage}{.55\textwidth}\centering
\includegraphics[width=1.\textwidth]{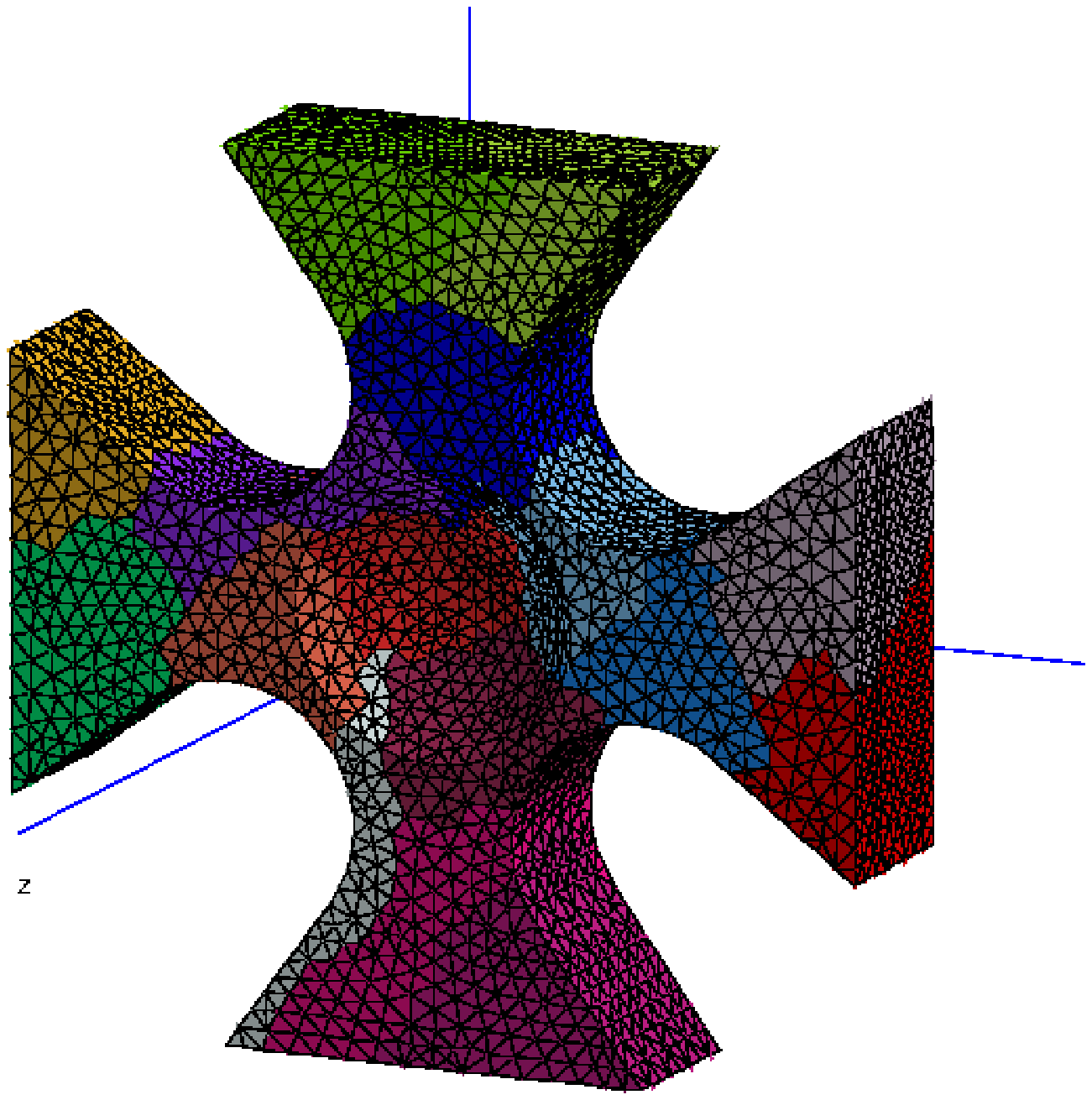}
\par\caption{16 subdomains bitraction test specimen}
(courtesy of ONERA -- Pascale Kanout\'e)\label{fig:croix}
\end{minipage}
\end{figure}

Hermann Schwarz (1843-1921) is often referred to as the father of
domain decomposition methods. In a 1869-paper he proposed an
alternating method to solve a PDE equation set on a complex domain
composed the overlapping union of a disk and a square (fig.
\ref{fig:ddm}), giving the mathematical basis of what is nowadays
one of the most natural ways to benefit the modern hardware
architecture of scientific computers.

In fact the growing importance of domain decomposition methods in
scientific computation is deeply linked to the growth of parallel
processing capabilities (in terms of number of processors, data
exchange bandwidth between processors, parallel library
efficiency, and of course performance of each processor). Because
of the exponential increase of computational resource requirement
for numerical simulation of more and more complex physical
phenomena (non-linearities, couplings between physical mechanisms
or between physical scales, random variables...) and more and more
complex problems (optimization, inverse problems...), parallel
processing appears to be an essential tool to handle the resulting
numerical models.

Parallel processing is supposed to take care of two key-points of
modern computations, the amount of operations and the required
memory. Let us consider the simulation of a physical phenomenon,
classically modeled by a PDE $\mathcal{L}(x)=f,\ x\in H(\Omega)$.
To take advantage of the parallel architecture of a calculator, a
reflexion has to be carried out on how the original problem could
be decomposed into collaborating subprocesses. The criteria for
this decomposition will be: first, the ability to solve
independent problems (on independent processors); second, how
often processes have to be synchronized; and last what quantity of
data has to be exchanged when synchronizing. When tracing back the
idle time of resolution processes and analyzing hardware and
software solutions, it is often observed that inter-processor
communications are the most penalizing steps.

If we now consider the three great classes of mathematical
decomposition of our reference problem which are operator
splitting (for instance $\mathcal{L}=\sum_i\mathcal{L}_i$
\cite{FORTIN:1982:MLA}), function-space decomposition (for
instance $H(\Omega)=\operatorname{span}(v_i)$, an example of which
is modal decomposition) and domain decomposition
($\Omega=\bigcup\Omega_i$), though the first two can lead to very
elegant formulations, only domain decompositions ensure (once one
subdomain or more have been associated to one processor) that
independent computations are limited to small quantities and that
the data to exchange is limited to the interface (or small
overlap) between subdomains which is always one-to-one exchange of
small amount of data.

So domain decomposition methods perfectly fit the criteria for
building efficient algorithms running on parallel computers. Their
use is very natural in engineering (and more precisely design)
context : domain decomposition methods offer a framework where
different design services can provide the virtual models of their
own parts of a structure, each assessed independently, domain
decomposition can then evaluate the behavior of the complete
structure just setting specific behavior on the interface (perfect
joint, unilateral contact, friction). Of course domain
decomposition methods also work with one-piece structure (for
instance fig. \ref{fig:croix}), then decomposition can be
automated according to criteria which will be discussed later.

From an implementation point of view, programming domain
decomposition methods is not an overwhelming task. Most often it
can be added to existing solvers as a upper level of current code
using the existing code as a black-box. The only requirement to
implement domain decomposition is to be able to detect the
interface between subdomains and use a protocol to share data on
this common part. In this paper we will mostly focus on domain
decomposition methods applied to finite element method
\cite{CIARLET:1979:FEM,ZIENKIEWICZ:1989:FEM}, anyhow they can be
applied to any discretization method (among others meshless
methods \cite{BELYTSCHKO:1996:MMO,BREITKOPT:2002:MPA,LIU:2005:IMM}
and discrete element methods
\cite{DADDETTA:2004:PCH,BOLANDER:2005:ILM,DELAPLACE:2005:FDF}).

Though domain decomposition methods were more than one century
old, they had not been extensively studied. Recent interest arose
as they were understood to be well-suited to modern engineering
and modern computational hardware. An important date in reinterest
in domain decomposition methods is 1987 as first international
congress dedicated to these methods occurred and DDM association
was created (see \url{http://www.ddm.org}).

Yet the studies were first mathematical analysis oriented and
emphasized on Schwarz overlapping family of algorithms. As
interest in engineering problems grew, non-overlapping Schwarz and
Schur methods, and coupling with discretization methods (mainly
finite element) were more and more studied. Indeed, these methods
are very natural to interpret mechanically, and moreover
mechanical considerations often resulted in improvement to the
methods. Basically the notion of interface between neighboring
subdomains is a strong physical concept, to which is linked a set
of conservation principles and phenomenological laws: for instance
the conservation of fluxes (action-reaction principle) imposes the
pointwise mechanical equilibrium of the interface and the equality
of incoming mass (heat...) from one subdomain to the outgoing mass
(heat...) of its neighbors; the "perfect interface" law consists
in supposing that displacement field (pressure, temperature) is
continuous at the interface, contact laws enable disjunction of
subdomains but prohibit interpenetration.

In this context two methods arose in the beginning of the 90's :
so-called Finite Element Tearing and Interconnecting (FETI) method
\cite{FARHAT:1991:FETI} and Balanced Domain Decomposition (BDD)
\cite{LETALLEC:1991:DDM}. From a mechanical point of view BDD
consists in choosing the interface displacement field as main
unknown while FETI consists in privileging the interface effort
field. BDD is usually referred to as a primal approach while FETI
is a dual approach. One of the interests of these methods, beyond
their simple mechanical interpretation, is that they can easily be
explained from a purely algebraic point of view (\textit{ie}
directly from the matrix form of the problem). In order to fit
parallelization criteria, it clearly appeared that the interface
problem should be solved using an iterative solver, each iteration
requiring local (\textit{ie} independent on each subdomain)
resolution of finite element problem, which could be done with a
direct solver. Then these methods combined direct and iterative
solver trying to mix robustness of the first and cheapness of the
second. Moreover the use of an iterative solver was made more
efficient by the existence of relevant preconditioners (based on
the resolution of a local dual problem for the primal approach and
a primal local problem for the dual approach).

When first released, FETI could not handle floating substructures
(\textit{ie} substructures without enough Dirichlet conditions),
thus limiting the choice of decomposition, while the primal
approach could handle such substructures but with loss of
scalability (convergence decayed as the number of floating
substructures increased). A key point then was the introduction of
rigid body motions as constraints and the use of generalized
inverses. Because of its strong connections with multigrid methods
\cite{WESSELING:2004:IMM}, the rigid body motions constraint took
the name of "coarse problem", it made the primal and dual methods
able to handle most decompositions without loss of scalability
\cite{FARHAT:1994:ADV,MANDEL:1993:BAL,LETALLEC:1994:DDM}. From a
mechanical point of view, the coarse problem enables
non-neighboring subdomains to interact without requiring the
transmission of data through intermediate subdomains, it then
enables to spread global information on the whole structure scale.

Once equipped with their best preconditioners and coarse problems,
mathematical results
\cite{FARHAT:1994:OPT,KLAWONN:2001:FNN,BRENNER:2005:LBD} provide
theoretical scalability of the methods. For instance for 3D
elasticity problems, if $h$ is the diameter of finite elements and
$H$ the diameter of subdomains, condition number $\kappa$ of the
interface problem reads ($C$ is a real constant):
\begin{equation}
\kappa \simeq C \left(1+\log\left(\frac{H}{h}\right)\right)^2
\end{equation}
which proves that the condition number only depends
logarithmatically on the number of elements per subdomain. Many
numerical assessment campaigns confirmed the good properties of
the methods, their robustness compared to iterative solvers
applied to the complete structure and their low cost (in terms of
memory and CPU requirements) compared to direct solvers. Thus
because they are well-suited to modern hardware (like PC clusters)
they enable to achieve computations which could not be realized on
classical computers because of too high memory requirement or too
long computational time: these methods can handle problems with
several millions of degrees of freedom.

Primal and dual methods were extended to heterogeneous problems by
a cheap intervention on the preconditioners
\cite{RIXEN:1998:SUPERL} and on the initialization
\cite{GOSSELET:2003:IEI}, and to forth order elasticity (plates
and shells) problems by the adjunction of so-called "second level
problem" in order to regularize the displacement field around the
corners of subdomains
\cite{FARHAT:1996:COR1,LETALLEC:1997:SHELL,FARHAT:2000:FETINL}. As
it became clear that the regularization of the displacement field
was sufficient to suppress rigid body motions, specific algorithms
which regularized \textit{a priori} the subdomain problems were
proposed: FETIDP \cite{FARHAT:2001:FETI_DP} and its primal
counterpart BDDC \cite{CROS:2002:PSC}, first in the plates and
shells context, then in the second order elasticity context
\cite{KLAWONN:2005:SCR}. Now FETIDP and BDDC are considered as
efficient as original FETI and BDD.

Methods were employed in many other contexts: transient dynamics
\cite{FRAGAKIS:2004:MHP2}, multifield problems (multiphysic
problems such as porous media \cite{GOSSELET:2003:AHD} and
constrained problems such as incompressible flows
\cite{LI:2002:DPF,GOLDFELD:2002:BNN}), Helmotz equations
\cite{DELABOURDONNAYE:1998:FETIH,FARHAT:1999:FETIH,FARHAT:2000:TLD}
and contact \cite{DUREISSEIX:2001:CONTACT,DOSTAL:2005:OSF}. The
use of domain decomposition methods in structural dynamic analysis
is a rather old idea though it can now be confronted to well
established methods in static analysis; the Craig-Bampton
algorithm \cite{CRAIG:1968:CSD} is somehow the application of the
primal strategy to such problems, the dual version of which was
proposed in \cite{RIXEN:2004:DCB}, moreover ideas like the
adjunction of coarse problems enabled to improve these methods.

Because of the strong connection between primal and dual
approaches, some methods try to propose frameworks which
generalize the two methods. The hybrid approach
\cite{GOSSELET:2003:THESE,GOSSELET:2004:DDM} enables to select
specific treatment (primal or dual) for each interface degree of
freedom; if all degrees of freedom have the same treatment, the
hybrid approach is exactly a classical approach. For certain
multifield problems the hybrid approach enables to define
physic-friendly solvers. The hybrid approach can also be obtained 
from specific optimality considerations \cite{dolean-2005}.
Mixed approaches \cite{LAD99b,NOUY:2003:THES,SERIES:2003:MDD2} consist in searching
a linear combination of interface displacement and effort field,
depending on the artificial stiffness introduced on the interface
one can recover the classical approaches (null stiffness for the
dual approach, infinite stiffness for the primal approach).
Moreover, the mixed approaches enable to provide the interface
mechanical behavior and provide a more natural framework to handle
complex interfaces (contact, friction...) than classical
approaches.

In this paper we aim at reviewing most of non-overlapping domain
decomposition methods. To adopt a rather general point of view we
introduce a set of notations strongly linked to mechanical
considerations as it gives the interface the main role of the
methods. We try to include all methods inside the same pattern so
that we can easily highlight the connections and differences
between them. We adopt a practical point of view as we describe
the mechanical concepts, the algebraic formulations, the
algorithms and the practical implementation of the methods. At
each step we try to emphasize keypoints and not to avoid
theoretical and practical difficulties.

This paper is organized as follows. In section \ref{sec:gene} we
introduce the mechanical framework of our study, the common
notations and the notion of interface assembly operators and
mechanical operators which will play a central role in the
methods. Section \ref{sec:methods} provides a rather extensive
review of the nonoverlapping domain decomposition methods in the
framework of discretized problems: basic primal and dual
approaches (with their variations), three-field method for
conforming grids, mixed and hybrid approaches. A keypoint of the
previous methods is the adjunction of optional constraints to form
a "coarse problem" which transmits global data through the whole
structure, the strategies to introduce these optional constraints
are studied in section \ref{sec:constraints}, which leads to the
definition of "recondensed" strategies FETIDP and BDDC. Section
\ref{sec:issu} deals with practical issues which are very often
common to most of the methods. Assessments are given in section
\ref{sec:asse} to illustrate the methods and outline their main
properties. Section \ref{sec:conc} concludes the paper. As Krylov
iterative solvers are often coupled to domain decomposition
methods, main concepts and algorithms to use them are given in
appendix \ref{sec:solvers}.

\section{Formulation of an interface problem}\label{sec:gene}
To present as smoothly as possible non-overlapping domain
decomposition methods we first consider a reference continuous
mechanics problem, decompose the domain in two subdomains in order
to introduce interface  fields, then in order to describe
correctly the interface we study a $N$-subdomain decomposition.
Since our aim is not to prove theoretical performance results but
make the reader feel some key-points of substructuring, we do not
go too far in continuous formulation and quickly introduce
discretized systems.

\subsection{Reference problem}

\begin{figure}[ht]\centering
\psfrag{O}{$\Omega$} \psfrag{f}{$\u{g}$} \psfrag{g}{$\u{f}$}
\psfrag{u}{$\u{u}_0$} \psfrag{a}{$\partial_g\Omega$} \psfrag{b}{$\partial_u\Omega$}
\includegraphics[height=3.cm]{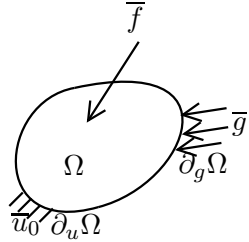}
\caption{Reference problem}\label{fig:omega:1}
\end{figure}

Let us consider domain $\Omega$ in $\mathbb{R}^n$ (n=1, 2 or 3)
submitted to a classical linear elasticity problem (see figure
\ref{fig:omega:1}) : displacement $\u{u}_0$ is imposed on part
$\partial_u\Omega$ of the boundary of the domain, effort $\u{g}$
is imposed on complementary part $\partial_g\Omega$, volumic
effort $\u{f}$ is imposed on $\Omega$, elasticity tensor is
$\mathfrak{a}$ \cite{GERMAIN:1986:M,SALENCON:1988:MMC}. The system
is governed by the following equations :
\begin{equation}\label{eq:elas:1}
\left\{
  \begin{array}{lcl}
    \u{\operatorname{div}}(\uu{\sigma})+\u{f}=\u{0} & & \text{in\ }\Omega\\
    \uu{\sigma}=\mathfrak{a}:\uu{\varepsilon}(\u{u}) & & \text{in\ }\Omega\\
    \uu{\varepsilon}(\u{u})=\frac{1}{2}\left( \uu{\operatorname{grad}}(\u{u})+\uu{\operatorname{grad}}(\u{u})^T
    \right)& & \text{in\ }\Omega\\
    \uu{\sigma}.\u{n}=\u{g}& & \text{on\ }\partial_g\Omega\\
    \u{u}=\u{u}_0  &&\text{on\ }\partial_u\Omega
  \end{array}
\right.
\end{equation}

In order to have the problem well posed, we suppose
$\operatorname{mes}(\partial_u \Omega)>0$. We also suppose that
tensor $\mathfrak{a}$ defines a symmetric definite positive
bilinear form on 2nd-order symmetric tensors. Under these
assumptions, problem \eqref{eq:elas:1} has a unique solution
\cite{DUVAUT:1990:MMC}.

\subsection{Two-subdomain decomposition}
\begin{figure}[H]\centering
\psfrag{a}{$\Omega^{(1)}$} \psfrag{b}{$\Omega^{(2)}$} \psfrag{c}{$\Upsilon$}
\includegraphics[height=3.cm]{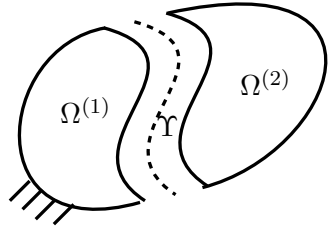}
\caption{Two-subdomain decomposition}\label{fig:omega:2}
\end{figure}
Let us consider a partition of domain $\Omega$ into $2$ substructures
$\Omega^{(1)}$ and $\Omega^{(2)}$. We define interface
$\Upsilon$ between substructures (figure \ref{fig:omega:2}) :
\begin{equation}\label{eq:interf:1}
    \Upsilon=\partial \Omega^{(1)} \bigcap \partial \Omega^{(2)} \\
\end{equation}
System \eqref{eq:elas:1} is posed on domain $\Omega$, we write its
restrictions to $\Omega^{(1)}$ and $\Omega^{(2)}$ :
\begin{equation}\label{eq:elas_dd:1}
s = 1 \text{\ or\ }2,\ \left\{
  \begin{array}{lcl}
    \u{\operatorname{div}}(\uu{\sigma}\s)+\u{f}\s=\u{0} & & \text{in\ }\Omega\s\\
    \uu{\sigma}\s=\mathfrak{a}\s:\uu{\varepsilon}(\u{u}\s) & & \text{in\ }\Omega\s\\
    \uu{\varepsilon}(\u{u}\s)=\frac{1}{2}\left( \uu{\operatorname{grad}}(\u{u}\s)+\uu{\operatorname{grad}}(\u{u}\s)^T
    \right)& & \text{in\ }\Omega\s\\
    \uu{\sigma}\s.\u{n}\s=\u{g}\s& & \text{on\ }\partial_g \Omega \bigcap \partial \Omega\s\\
    \u{u}\s=\u{u}_0\s  &&\text{on\ } \partial_u \Omega \bigcap \partial \Omega\s
  \end{array}
\right.
\end{equation}
and the interface connection conditions, continuity of displacement
\begin{equation}\label{eq:elas_dd:2}
    \u{u}^{(1)}=\u{u}^{(2)}   \qquad \text{on\ }\Upsilon
\end{equation}
and equilibrium of efforts (action-reaction principle)
\begin{equation}\label{eq:elas_dd:3}
\uu{\sigma}^{(1)} \u{n}^{(1)} + \uu{\sigma}^{(2)}
\u{n}^{(2)}=\u{0}    \qquad \text{on\ }\Upsilon
\end{equation}
Of course system (\ref{eq:elas_dd:1}, \ref{eq:elas_dd:2}, \ref{eq:elas_dd:3}) is strictly equivalent to global problem \eqref{eq:elas:1}.

\subsection{N-subdomain decomposition}

\begin{figure}[ht]\centering
\psfrag{a}{$\Omega^{(1)}$}\psfrag{b}{$\Omega^{(2)}$}\psfrag{c}{$\Omega^{(3)}$}
\psfrag{d}{$\Upsilon^{(1,2)}$}\psfrag{e}{$\Upsilon^{(1,3)}$}\psfrag{f}{$\Upsilon^{(2,3)}$}
\psfrag{h}{$\Upsilon$} \subfigure[Primal (geometric)
interface]{\includegraphics[width=0.3\textwidth]{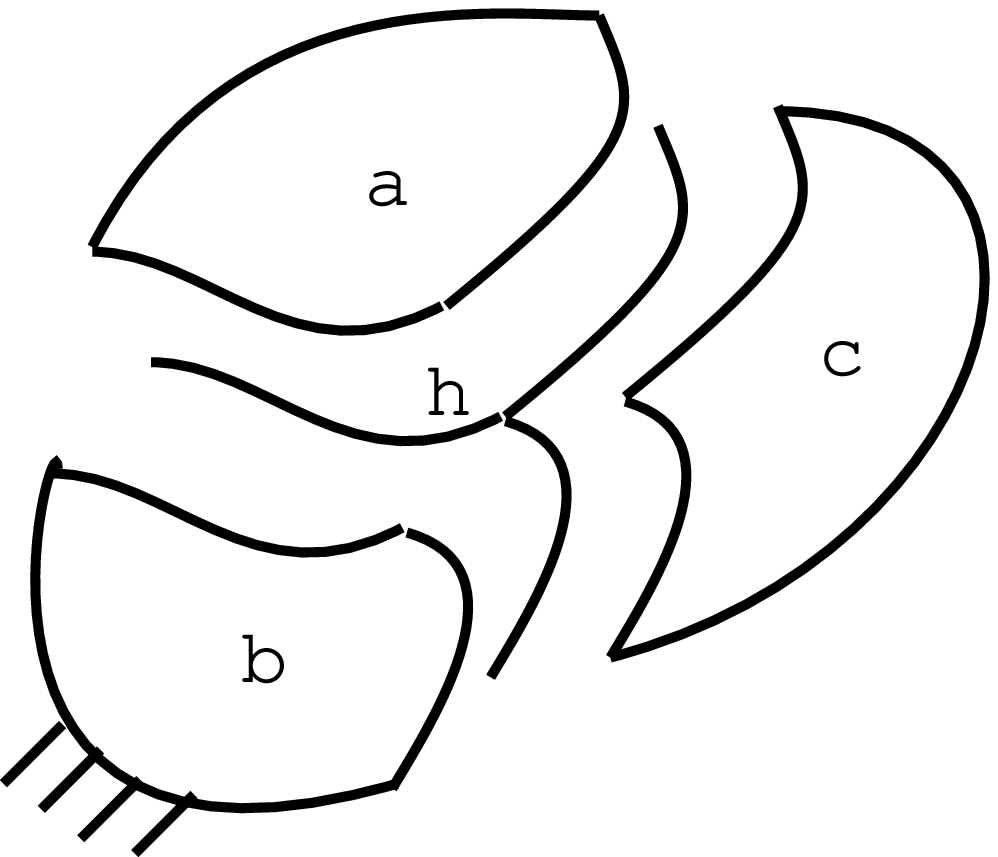}}
\qquad\subfigure[Dual interface
(connectivity)]{\includegraphics[width=0.3\textwidth]{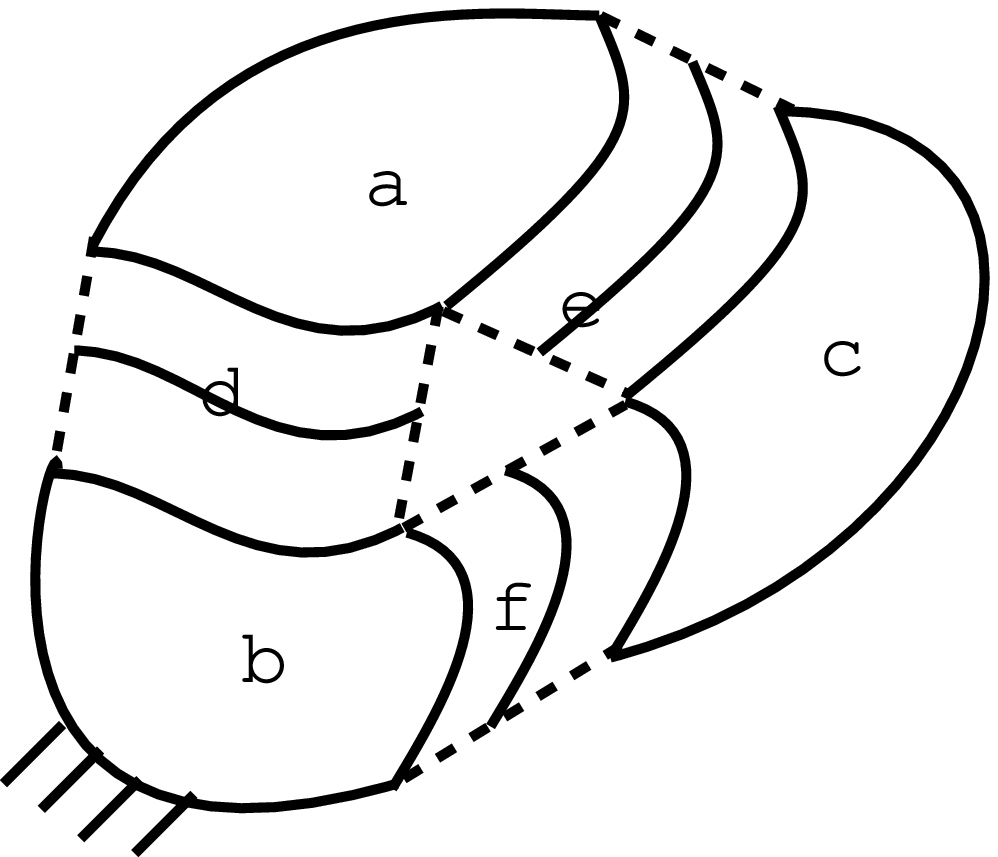}}\caption{Definition
of the interface for a $N$-subdomain
decomposition}\label{fig:omega:3}
\end{figure}

Let us consider a partition of domain $\Omega$ into $N$ subdomains
denoted $\Omega\s$. We can define the interface between two
subdomains, the complete interface of one subdomain and the
geometric interface at the complete structure scale:
\begin{equation}\label{eq:interN:1}
\left\{
  \begin{array}{lcl}
    \Upsilon^{(i,j)}=\Upsilon^{(j,i)}=\partial \Omega^{(i)} \bigcap \partial \Omega^{(j)} \\
    \Upsilon\s = \bigcup_j \Upsilon^{(s,j)}\\
    \Upsilon = \bigcup_s \Upsilon\s
  \end{array}
\right.
\end{equation}

When implementing the method, one (possibly virtual) processor is
commonly assigned to each subdomain, hence because we can tell
"local" computations (realized independently on each processor)
from "global" computations (realized by exchanging data between
processors), we often refer to values as being global or local.
Then $\Upsilon\s$ is the local interface and $\Upsilon$ the global
interface. Because exchanges are most often one-to-one,
$\Upsilon^{(i,j)}$ is the $(i-j)$-communication interface.

Using more than two subdomains (except when using
"band"-decomposition) leads to the appearance of "multiple-points"
also called "crosspoints" (which are nodes shared by more than two
subdomains). These crosspoints lead to the existence of two
descriptions of the interface (figure \ref{fig:omega:3}):
so-called geometric interface $\Upsilon$ and so-called
connectivity interface made out of the set of one-to-one
interfaces $(\Upsilon^{(i,j)})_{1\leqslant i< j \leqslant N}$.
Each of the two most classical methods exclusively uses one of
these descriptions so the geometric interface $\Upsilon$ is often
referred to as the primal interface while the connectivity
interface $\du{\Upsilon}$ is referred to as the dual interface.

\comment{
On a mathematical point of view, interfaces are described by different trace spaces:
\begin{eqnarray}\label{eq:espacetrace:1}
    \mathcal{Y}_\Upsilon &=& H^{1/2}\left(\bigcup\limits_{1\leqslant i<j\leqslant N} \Upsilon^{(i,j)}\right) \\
    \du{\mathcal{Y}}_\Upsilon&=&\bigotimes_{1\leqslant i<j\leqslant
    N}H^{1/2}(\Upsilon^{(i,j)})
\end{eqnarray}

Working on trace spaces \cite{DAUTRAY:1987:AMC} is not the purpose
of the paper, anyhow we would like to emphasize the essential
duality between geometric and connectivity interfaces: as reader
may observe on figure \ref{fig:omega:3} there is a Delaunay /
Vorono�duality between these two descriptions, which translate to
dual mathematical spaces (not explained in previous equations) and
then to dual formulations of interface problem.} How crosspoints
are handled is a fundamental key in the differentiation of domain
decomposition methods. In the remainder of the paper we will
always refer to data attached to the dual interface using
underlined notation.

\begin{remark}
Reader may have observed that the above presented connectivity
description is redundant for crosspoints:  let $x$ be a crosspoint
if $x$ belongs to  $\Upsilon^{(1,2)}$ and $\Upsilon^{(2,3)}$, it
of course belongs to  $\Upsilon^{(1,3)}$. In the general case of a
$m$-multiple crosspoint, there are $\binom{m}{2}$ connectivity
relationships while only  $(m-1)$ would be sufficient and
necessary. We will present strategies to remove these redundancies
in the algebraic analysis of the methods.
\end{remark}

\begin{remark}
Cross-points may also introduce, at the continuous level,
punctual-interfaces in 2d or edge-interfaces in 3d, which are
interface with zero measure. Most often from a physical point of
view these are not considered as interfaces. Anyhow after
discretization all relationships are written node-to-node and the
problem no longer exists.
\end{remark}

\subsection{Discretization}
We suppose that the reference problem has been discretized, leading to the resolution of $n\times n$ linear system:
\begin{equation}\label{eq:sysdis:1}
Ku=f
\end{equation}
Because of its key role in structural mechanics we will often
refer to finite element discretization though any other technique
would suit. The key points are the link between matrix $K$ and the
domain geometry and the sparse filling  of matrix $K$ (related to
the fact that only narrow nodes have non-zero interaction).

We restrict to the case of element-oriented decompositions (each
element belongs to one and only one substructure) which are
conforming to the mesh which implies three conditions
\cite{RIXEN:1997:THESE}:
\begin{itemize} \item there is a one-to-one correspondence between degrees of freedom by the interface ; \item approximation spaces are the same by the interface ; \item models (beam, shell, 3d...) are the same by the interface.
\end{itemize}
Under these assumptions connection conditions simply write as node
equalities. For non-conforming meshes, a classical solution is to
use mortar-elements for which continuity and equilibrium are
verified in a weak sense
\cite{ACHDOU:1995:SPF,ACHDOU:1996:MIS,STEFANICA:1999:FME}.

\subsubsection{Boolean operators}
In order to write communication relation between subdomains we have to introduce several operators.

The first one is the "local trace" operator $\traceop\s$ which is
the discrete projection from $\Omega\s$ to $\Upsilon\s$. It enables
to cast data from a complete subdomain to its interface, and once
transposed to extend data set on the interface to the whole
subdomain (setting internal degrees of freedom to zero). In the
remainder of the paper we will use subscript $b$ for interface
data and subscript $i$ for internal data.

\begin{figure}[H]\centering
\psfrag{a}{$1^{(1)}$}\psfrag{b}{$2^{(1)}$}\psfrag{c}{$3^{(1)}$}\psfrag{d}{$4^{(1)}$}\psfrag{e}{$5^{(1)}$}
\psfrag{f}{$1^{(2)}$}\psfrag{g}{$2^{(2)}$}\psfrag{h}{$3^{(2)}$}\psfrag{i}{$4^{(2)}$}\psfrag{j}{$5^{(2)}$}
\psfrag{k}{$1^{(3)}$}\psfrag{l}{$2^{(3)}$}\psfrag{m}{$3^{(3)}$}\psfrag{n}{$4^{(3)}$}
\psfrag{ct}{$1_b^{(1)}$}\psfrag{dt}{$2_b^{(1)}$}\psfrag{et}{$3_b^{(1)}$}
\psfrag{ht}{$1_b^{(2)}$}\psfrag{it}{$2_b^{(2)}$}\psfrag{jt}{$3_b^{(2)}$}
\psfrag{kt}{$1_b^{(3)}$}\psfrag{lt}{$2_b^{(3)}$}\psfrag{mt}{$3_b^{(3)}$}
\psfrag{o}{$1_\Upsilon$}\psfrag{p}{$2_\Upsilon$}\psfrag{q}{$3_\Upsilon$}\psfrag{r}{$4_\Upsilon$}
\psfrag{s}{$\du{1}_\Upsilon$}\psfrag{t}{$\du{2}_\Upsilon$}\psfrag{u}{$\du{3}_\Upsilon$}\psfrag{v}{$\du{4}_\Upsilon$}\psfrag{w}{$\du{5}_\Upsilon$}\psfrag{x}{$\du{6}_\Upsilon$}
\subfigure[Subdomains]{\includegraphics[width=5.cm]{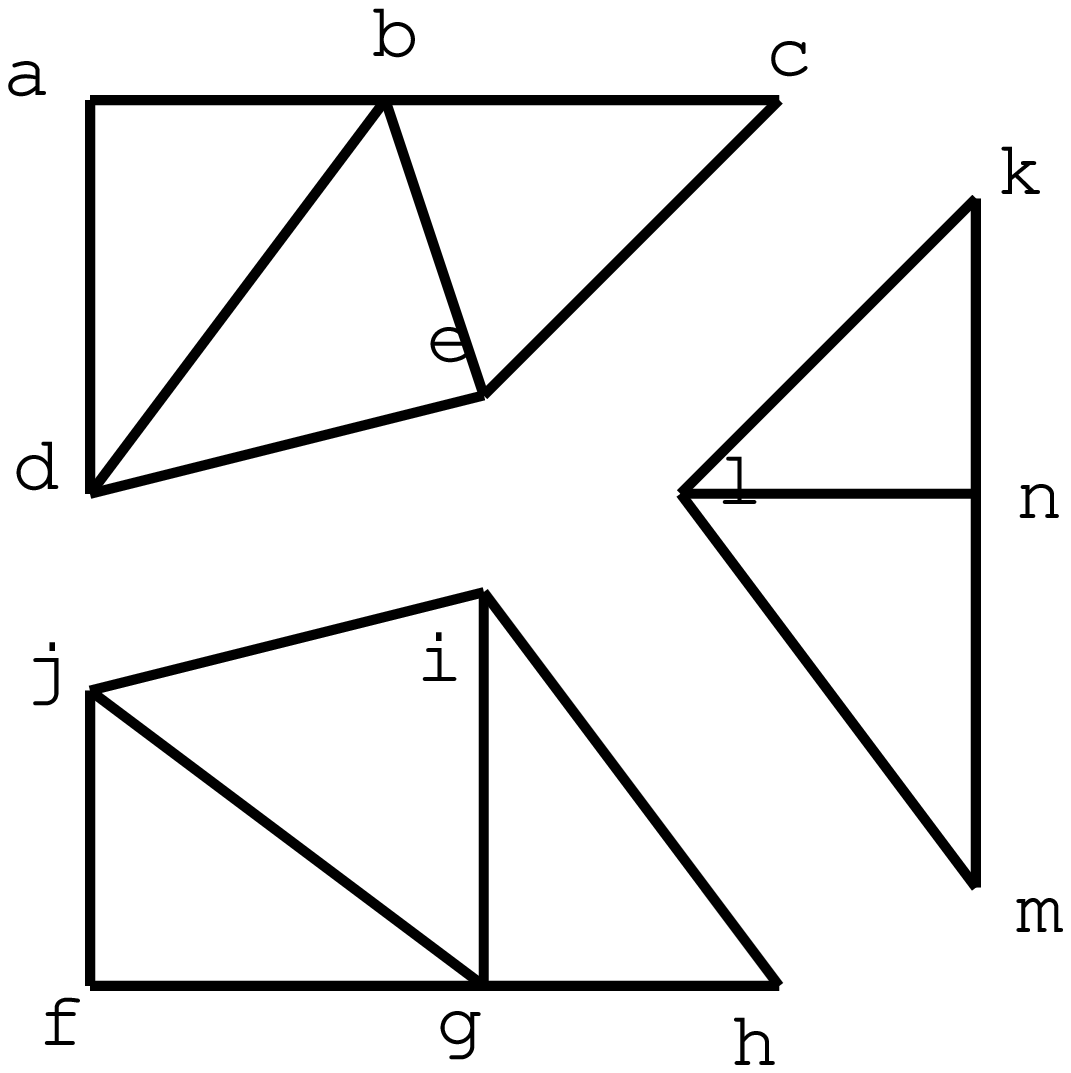}
}\qquad
\subfigure[Local interface]{\includegraphics[width=5.cm]{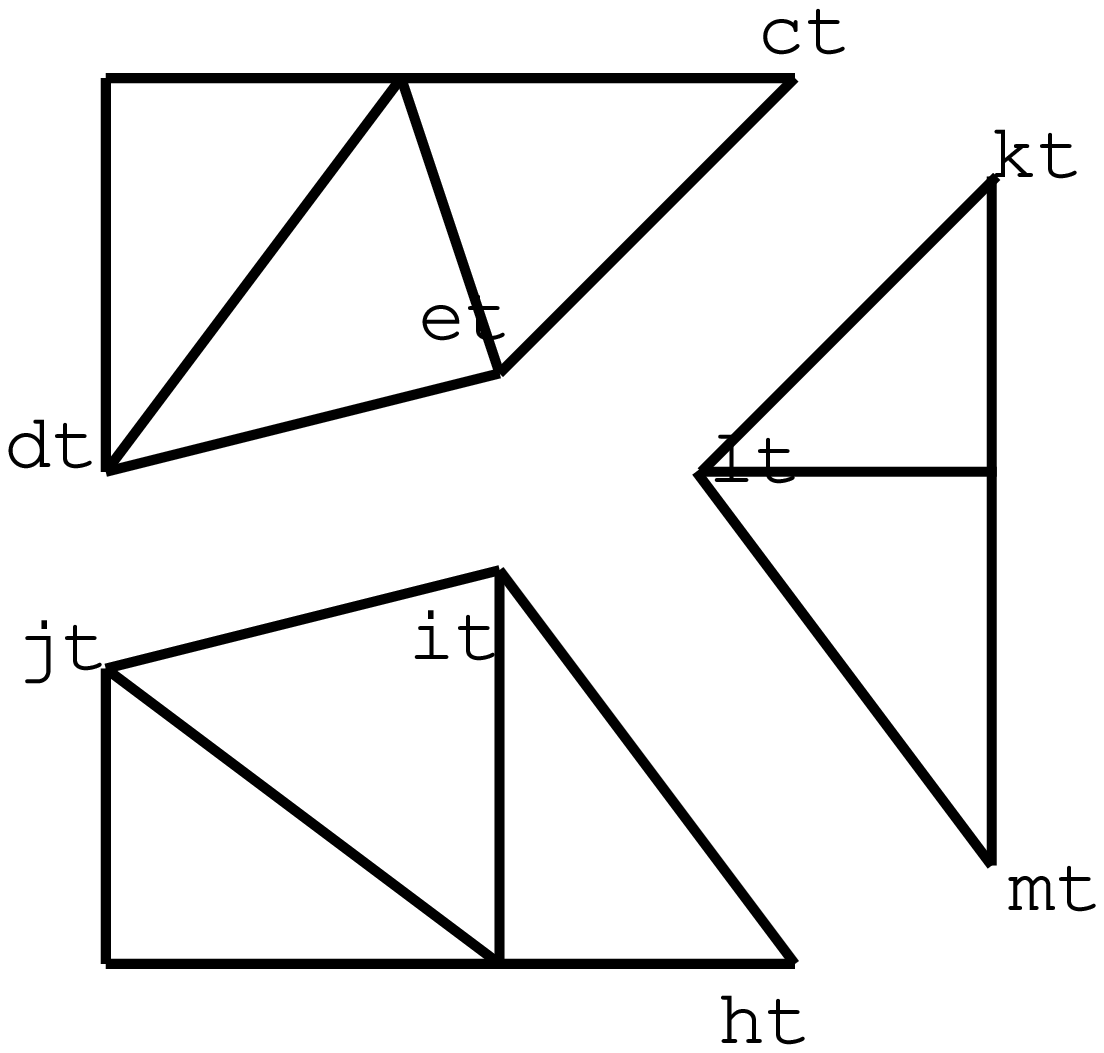} }\\
\subfigure[Primal interface]{\includegraphics[width=5.cm]{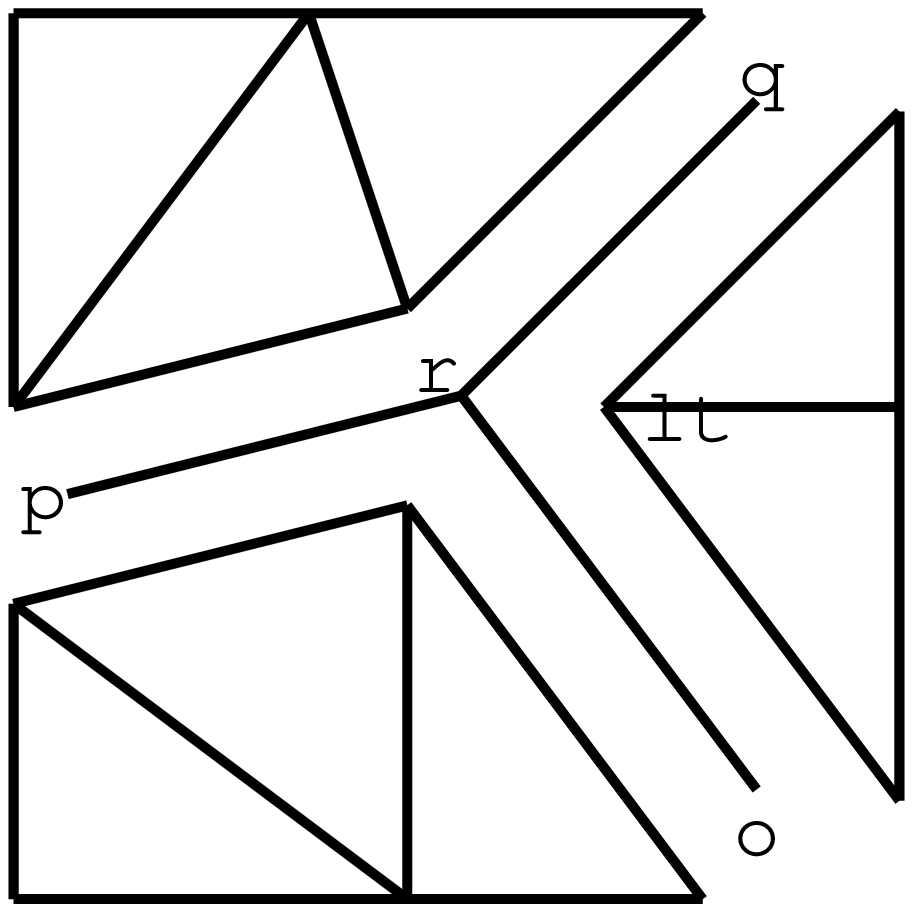}
}\qquad \subfigure[Dual interface]{\includegraphics[width=5.cm]{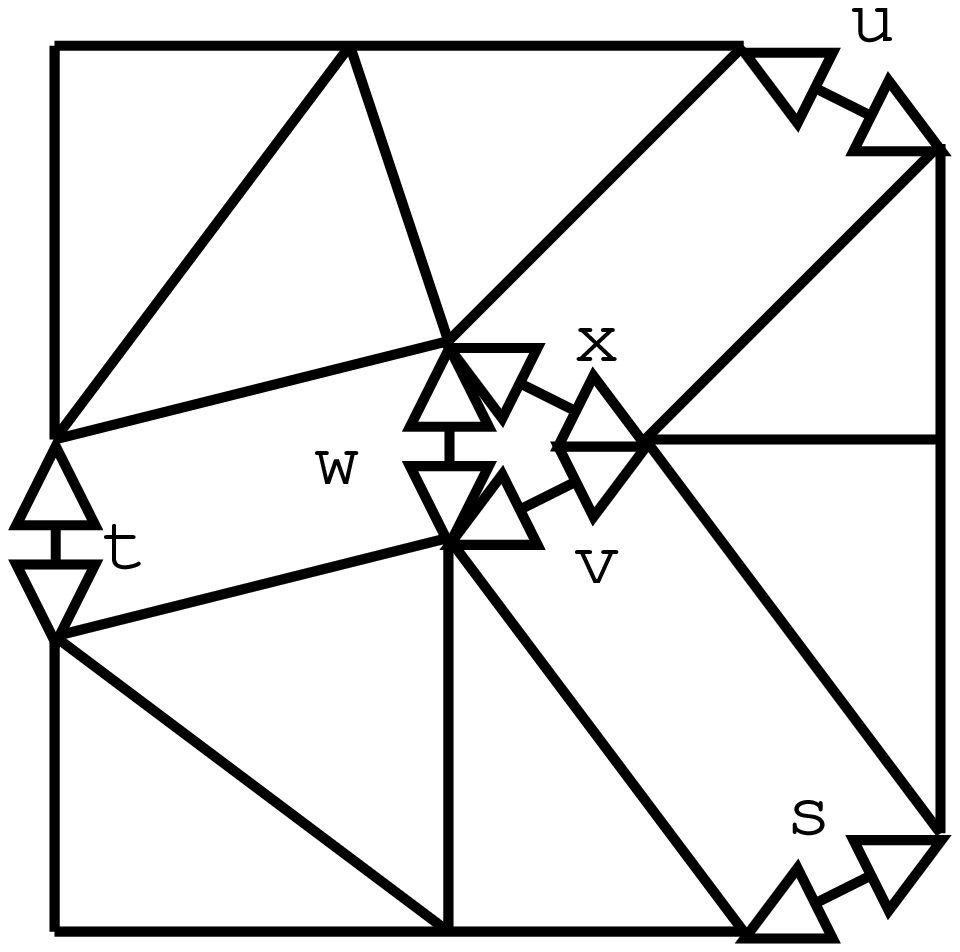}
} \scriptsize\begin{equation*}
\begin{array}{rclcrclcrcl}
\tras{1}&=&\begin{pmatrix}0&0&1&0&0\\0&0&0&1&0\\0&0&0&0&1\end{pmatrix}&&\tras{2}&=&\begin{pmatrix}0&0&1&0&0\\0&0&0&1&0\\0&0&0&0&1\end{pmatrix}&&\tras{3}&=&\begin{pmatrix}1&0&0&0\\0&1&0&0\\0&0&1&0\end{pmatrix}\\
\pass{1}&=&\begin{pmatrix}0&0&0\\0&1&0\\1&0&0\\0&0&1\end{pmatrix}&&\pass{2}&=&\begin{pmatrix}1&0&0\\0&0&1\\0&0&0\\0&1&0\end{pmatrix}&&\pass{3}&=&\begin{pmatrix}0&0&1\\0&0&0\\1&0&0\\0&1&0\end{pmatrix}\\
\dass{1}&=&\begin{pmatrix}0&0&0\\0&1&0\\1&0&0\\0&0&0\\0&0&1\\0&0&1\end{pmatrix}&&\dass{2}&=&\begin{pmatrix}1&0&0\\0&0&-1\\0&0&0\\0&1&0\\0&-1&0\\0&0&0\end{pmatrix}&&\dass{3}&=&\begin{pmatrix}0&0&-1\\0&0&0\\-1&0&0\\0&-1&0\\0&0&0\\0&-1&0\end{pmatrix}
\end{array}\end{equation*}\normalsize
\caption{Local numberings, interface numberings, trace and assembly operators}\label{fig:omegef:2}
\end{figure}

Then data lying on one subdomain interface has to be exchanged
with its neighboring subdomains. It can be either realized on the
primal interface or the dual interface, leading to two (global)
"assembly" operators: the primal one $\pass{s}$, and the dual one
$\dass{s}$. The primal assembly operator is a strictly boolean
operator while the dual assembly operator is a signed boolean
operator (see figure \ref{fig:omegef:2}): if a degree of freedom
is set to 1 on one side of the interface, its corresponding degree
of freedom on the other side of the interface is set to $-1$.
Non-boolean assembly operators can be used in order to average
connection conditions when using non-conforming domain
decomposition methods \cite{BERNARDI:1990:MORT}.

\begin{remark}
Our $(\trass,\pass{s},\dass{s})$ set of operators is not the most
commonly used in papers related to domain decomposition. The
interest of this choice is to be sufficient to explain most of the
available strategies with only three operators. Other notations
use "composed operators" like ($B\s=\dass{s}\trass$ or
${L\s}^T=\pass{s}\trass$) which are not sufficient to describe all
methods and which, in a way, omit the fundamental role played by
the interface.
\end{remark}

Boolean operators have important classical properties. Please note
the first one which expresses the orthogonality of the two
assembly operators.
\begin{subequations}
\begin{equation}\label{eq:ortho_int:1}\sum_s{\dass{s}{\pass{s}}^T}=0\end{equation}
\begin{equation}\label{eq:ortho_int:2}{\pass{s}}^T\pass{s}=I_{\Upsilon\s}\end{equation}
\begin{equation}\label{eq:ortho_int:3}{\dass{s}}^T\dass{s}=\operatorname{diag}(\text{multiplicity}-1)_{\Upsilon\s} \end{equation}
\begin{equation}\label{eq:ortho_int:4}\pass{s}{\pass{s}}^T=\left| \begin{array}{l} I \text{\ on\ }\Upsilon\s \\ 0 \text{\ elsewhere} \end{array} \right.\end{equation}
\end{subequations}
\begin{remark}
An interesting choice of description would have been to use
redundant local interface (defining some kind of
$\du{\traceop}\s$). This choice would stick to most classical
implementations where the local interface of one subdomain is
defined neighborwise. Dual assembly operator would write easily as
a simple signing operator, but handling multiple points would be
slightly more difficult for the primal assembly operator (see
section \ref{sec:implementation}).
\end{remark}

\begin{figure}[H]\centering
\psfrag{ct}{$1_b^{(1)}$}\psfrag{dt}{$2_b^{(1)}$}\psfrag{et}{$3_b^{(1)}$}
\psfrag{ht}{$1_b^{(2)}$}\psfrag{it}{$2_b^{(2)}$}\psfrag{jt}{$3_b^{(2)}$}
\psfrag{kt}{$1_b^{(3)}$}\psfrag{lt}{$2_b^{(3)}$}\psfrag{mt}{$3_b^{(3)}$}
\psfrag{s}{$\du{1}_\Upsilon$}\psfrag{t}{$\du{2}_\Upsilon$}\psfrag{u}{$\du{3}_\Upsilon$}\psfrag{v}{$\du{4}_\Upsilon$}\psfrag{w}{$\du{5}_\Upsilon$}\psfrag{x}{$\du{6}_\Upsilon$}
\subfigure[Local interface]{\includegraphics[width=5.6cm]{fig/omegef2b.eps}
}\qquad
\subfigure[Redundant connectivity]{\includegraphics[width=5.cm]{fig/omegef2d.eps} }\\
\subfigure[Non-redundant connectivity]{\includegraphics[width=5.cm]{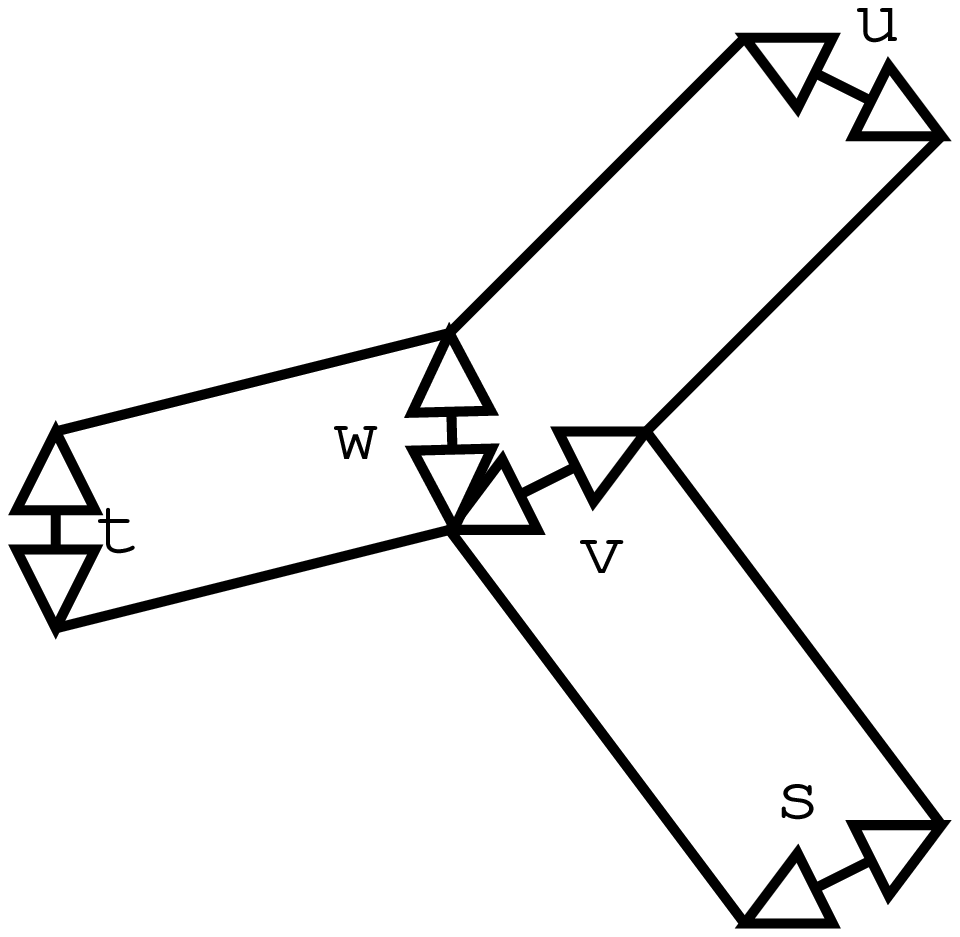} }\qquad
\subfigure[Orthonormal connectivity]{\includegraphics[width=5.cm]{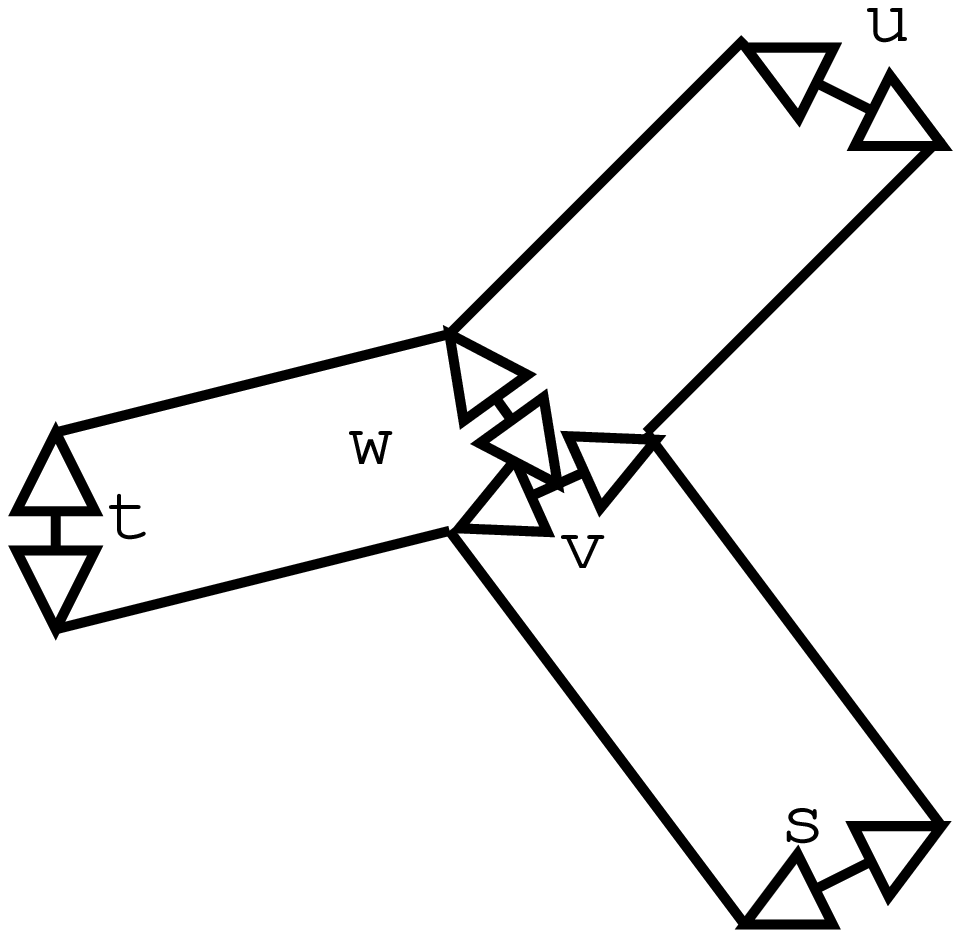} }
\scriptsize\begin{equation*}
\begin{array}{rclcrclcrcl}
\dass{1}&=&\begin{pmatrix}0&0&0\\0&1&0\\1&0&0\\0&0&0\\0&0&1\\0&0&1\end{pmatrix}&&\dass{2}&=&\begin{pmatrix}1&0&0\\0&0&-1\\0&0&0\\0&1&0\\0&-1&0\\0&0&0\end{pmatrix}&&\dass{3}&=&\begin{pmatrix}0&0&-1\\0&0&0\\-1&0&0\\0&-1&0\\0&0&0\\0&-1&0\end{pmatrix}\\
\dass{1}_{N}&=&\begin{pmatrix}0&0&0\\0&1&0\\1&0&0\\0&0&0\\0&0&1\end{pmatrix}&&\dass{2}_{N}&=&\begin{pmatrix}1&0&0\\0&0&-1\\0&0&0\\0&1&0\\0&-1&0\end{pmatrix}&&\dass{3}_{N}&=&\begin{pmatrix}0&0&-1\\0&0&0\\-1&0&0\\0&-1&0\\0&0&0\end{pmatrix}\\
\dass{1}_{O}&=&\begin{pmatrix}0&0&0\\0&\frac{1}{\sqrt{2}}&0\\\frac{1}{\sqrt{2}}&0&0\\0&0&0\\0&0&\frac{2}{\sqrt{6}}\end{pmatrix}&&\dass{2}_{O}&=&\begin{pmatrix}\frac{1}{\sqrt{2}}&0&0\\0&0&-\frac{1}{\sqrt{2}}\\0&0&0\\0&\frac{1}{\sqrt{2}}&0\\0&-\frac{1}{\sqrt{6}}&0\end{pmatrix}&&\dass{3}_{O}&=&\begin{pmatrix}0&0&-\frac{1}{\sqrt{2}}\\0&0&0\\-\frac{1}{\sqrt{2}}&0&0\\0&-\frac{1}{\sqrt{2}}&0\\0&-\frac{1}{\sqrt{6}}&0\end{pmatrix}\\
\end{array}\end{equation*}\normalsize
\caption{Suppressing redundancies of dual interface}\label{fig:omegef:3}
\end{figure}

\begin{remark}
Redundancies can easily be removed from the dual description of
the interface. One has just to modify the connectivity table, so
that one multiple point is connected only once to each subdomain.
This can be carried out introducing two different assembly
operators the "non-redundant" one and the "orthonormal" one (see
figure \ref{fig:omegef:3}) \cite{FRAGAKIS:2002:UFF}. Only
relationship \ref{eq:ortho_int:3} is modified (then
${\dass{s}}^T\dass{s}=I_{\Upsilon\s}$). The interest of the use of
these assembly operators will be discussed in section
\ref{sub:AFETI}.
\end{remark}

\subsubsection{Basic equations}
In order to rewrite equation \eqref{eq:sysdis:1} in a
domain-decomposed context, we have to introduce the reaction
unknown which is the discretization of $\uu{\sigma}^{(1)}
\u{n}^{(1)} = - \uu{\sigma}^{(2)} \u{n}^{(2)}$ in equation
\eqref{eq:elas_dd:3}. $\lambda\s$ is the reaction imposed by
neighboring subdomains on subdomain $(s)$. Commonly $\lambda\s$ is
defined on the whole subdomain $(s)$ while it is non-zero only on
its interface, so $\lambda\s={\trass}^T \lambda_b\s$.
\begin{subequations}\label{eq:ef_dd:1}
\begin{equation}\label{eq:ef_dd:1a} \forall s,\ K\s u\s = f\s + \lambda\s \end{equation}
\begin{equation}\label{eq:ef_dd:1b} \sum_s \dass{s} \trass u\s = 0 \end{equation}
\begin{equation}\label{eq:ef_dd:1c}  \sum_s \pass{s} \trass \lambda\s =0 \end{equation}
\end{subequations}

Equation \eqref{eq:ef_dd:1a} corresponds to the (local)
equilibrium of each subdomain submitted to external conditions
$f\s$ and reactions from neighbors $\lambda\s$. Equation
\eqref{eq:ef_dd:1b} corresponds to the (global) continuity of the
displacement field through the interface. Equation
\eqref{eq:ef_dd:1c} corresponds to the (global) equilibrium of the
interface (action-reaction principle).

This three-equation system \eqref{eq:ef_dd:1} is the starting
point from a rich zoology of methods, most of which possess strong
connections we will try to emphasis on. Before going further in
the exploration of these methods, we propose to introduce local
condensed operators that represent a subdomain on its interface,
then a set of synthetic notations.

\subsubsection{Local condensed operators}\label{sec:localop}
Philosophically, local condensed operators are operators that
represent how neighboring subdomains "see" one subdomain: a
subdomain can be viewed as a  black-box, the only information
necessary for neighbors is how it behaves on its interface.
Associated to this idea is the classical assumption that local
operations are "exactly" performed. From an implementation point
of view, when solving problems involving local matrices, a direct
solver is employed. As we will see, the use of exact local solvers
will be coupled with the use of iterative global solvers leading
to a powerful combination of speed and precision of computations.

In this section we will always refer to the local equilibrium of
subdomain $(s)$ under interface loading:
\begin{equation}\label{eq:localeq:1}
K\s u\s = \lambda\s = {\traceop\s}^T\lambda_b\s
\end{equation}

\begin{description}
\item[Primal Schur complement $S_p\s$: ]
If we renumber the local degrees of freedom of subdomain $(s)$ in
order to separate internal and boundary degrees of freedom, system
$\eqref{eq:localeq:1}$  writes
\begin{equation}\label{eq:localeq:2}
\begin{pmatrix}
K\s_{ii} & K\s_{ib} \\ K\s_{bi} & K\s_{bb}
\end{pmatrix}
\begin{pmatrix}
u\s_{i} \\ u\s_{b}
\end{pmatrix}=\begin{pmatrix}
0 \\ \lambda\s_{b}
\end{pmatrix}
\end{equation}
From the first line we draw
\begin{equation}\label{eq:condensation:1}
u\s_{i} = -{K\s_{ii}}^{-1}K\s_{ib}u\s_{b}
\end{equation}
then the Gauss elimination of $u\s_{i}$ leads to
\begin{equation}\label{eq:schurp:1}
\left(K\s_{bb}-K\s_{bi}{K\s_{ii}}^{-1}K\s_{ib}\right)u\s_{b}=S\s_p u\s_b = \lambda\s_b
\end{equation}
which is the condensed form of the local equilibrium of subdomains
expressed in terms of interface fields. Operator $S\s_p$ is called
local primal Schur complement. Its computation is realized by the
inversion of matrix $K\s_{ii}$ which corresponds to Dirichlet
conditions imposed on the interface of subdomain $(s)$, so the
primal Schur complement is always well defined, and commonly
called the "local Dirichlet operator". Note that the symmetry,
positivity, and definition properties are inherited by matrix
$S\s_p$ from matrix $K\s$.

An important result is that the kernel of matrices $K\s$ and
$S\s_p$ can be deduced one from the other ($I\s_b$ is the identity
matrix on the interface):
\begin{eqnarray}\label{eq:kernels}
K\s R\s = 0 & \Longrightarrow & S_p\s \traceop\s R\s = S_p\s R_b\s =0 \\
S_p\s R_b\s = 0  & \Longrightarrow & K\s \begin{pmatrix} -{K\s_{ii}}^{-1}K\s_{ib}\\ I\s_b \end{pmatrix}R_b\s = K\s R\s = 0
\end{eqnarray}

Primal Schur complement can also be interpreted as the
discretization of the Stecklov-Poincar\'e operator. From a
mechanical point of view, it is the linear operator that provides
the reaction associated to given interface displacement field.

If we consider that the subdomain is also loaded on internal
degrees of freedom, then the condensation of the equilibrium on
the interface reads:
\begin{eqnarray}
K\s u\s=f\s \Rightarrow S_p\s u_b\s &=& b_p\s  \\ b_p\s &=& f_b\s-K\s_{bi}{K\s_{ii}}^{-1} f_i\s
\end{eqnarray}
$b_p\s$ is the condensed effort imposed on the substructure.

\item[Dual Schur complement $S_d\s$: ] the dual Schur complement is a linear operator that computes interface displacement field from given interface effort.
From equation \eqref{eq:localeq:1} and \eqref{eq:schurp:1} we have:
\begin{equation}\label{eq:schurd:1}
\left({\traceop\s}{K\s}^+{\traceop\s}^T\right) \lambda_b\s ={S_p\s}^+ \lambda_b\s = S_d\s \lambda_b\s = u_b\s
\end{equation}
where ${K\s}^+$ is the generalized inverse of matrix $K\s$, and where it is  assumed that no rigid body motion is excited. If we denote by $R\s$ the kernel of matrix $K\s$ this last condition reads:
\begin{equation}\label{eq:admissibility:1}
{R\s}^T\lambda\s=0 \quad \text{\ or\ equivalently\ } \quad{R\s_b}^T\lambda_b\s=0
\end{equation}
\end{description}
\begin{remark}
A generalized inverse (or pseudo-inverse) of matrix $M$ is a
matrix, denoted $M^+$, which verifies the following property:
$\forall y\in \range(M),\ MM^+y=y$. Note that this definition
leads to non-unique generalized inverse, however all results
presented are independent of the choice of generalized inverse.
\end{remark}
Of course in order to take into account, inside
\eqref{eq:schurd:1}, the possibility of the substructure to have
zero energy modes, an arbitrary rigid displacement can be added
leading to the next expression where vector $\alpha\s$ denotes the
magnitude of rigid body motions:
\begin{equation}\label{eq:schurd:2}
u_b\s=S_d\s \lambda_b\s + R_b\s \alpha\s
\end{equation}
\begin{description}
\item[Hybrid Schur complement: $S_{pd}\s$]
this operator corresponds to an interface where degrees of freedom
are partitioned into two subsets. Suppose that the first subset is
submitted to given Dirichlet conditions and the second to Neumann
condition, $S_{pd}\s$ is the linear operator that associates
resulting reaction on the first subset and resulting displacement
on the second subset to those given conditions. We denote by
subscript $p$ data defined on the first subset and by subscript
$d$ data defined on the second subset (schematically $b=p\cup d$
and $p\cap d=\emptyset$).
\begin{equation}
S_{pd}\s \begin{pmatrix} u_p\s \\ \lambda_d\s \end{pmatrix} = \begin{pmatrix}  \lambda_p\s \\ u_d\s \end{pmatrix}
\end{equation}
The computation of operator $S_{pd}\s$ though no more complex than
the computation of $S_p\s$ or $S_d\s$, requires more notations. A
synthetic option is to denote by subscript $\bp$ the sets of
internal (subscript $i$) data and second-interface-subset
(subscript $d$) data, schematically $\bp=i\cup d$. We introduce a
modified trace operator:
\begin{equation}
\traceop\s_d v_{\bp}= \traceop\s_d \begin{pmatrix} v_i \\ v_d \end{pmatrix} = v_d
\end{equation}
Then internal equilibrium \eqref{eq:localeq:1} reads:
\begin{equation}
\begin{pmatrix}
K_{\bp\bp}\s & K_{\bp p}\s \\  K_{p \bp}\s & K_{p p}\s
\end{pmatrix}
\begin{pmatrix}
u_{\bp}\s \\  u_{p}\s
\end{pmatrix} =
\begin{pmatrix}
{\traceop\s_d}^T \lambda_{d}\s \\  \lambda_{p}\s
\end{pmatrix}
\end{equation}
Then hybrid Schur complement is:
\begin{equation}\label{eq:loc_hyb_sc:1}
S_{pd}\s = \begin{pmatrix}
K_{pp}\s -  K_{p\bp}\s{K_{\bp\bp}\s}^+K_{\bp p}\s & K_{p \bp}\s{K_{\bp\bp}\s}^+ {\traceop\s_d}^T \\
- {\traceop\s_d} {K_{\bp\bp}\s}^+ K_{\bp p} & {\traceop\s_d} {K_{\bp\bp}\s}^+ {\traceop\s_d}^T
\end{pmatrix}
\end{equation}
As can be noticed the diagonal blocks of $S_{pd}\s$ look like
fully primal and fully dual Schur complements, while extradiagonal
blocks are antisymmetric (assuming $K\s$ is symmetric). Of course
if all interface degrees of freedom belong to the same subset, the
hybrid Schur complement equals "classical" fully primal or fully
dual Schur complement. Moreover it stands out clearly that:
\begin{equation}
{S_{pd}\s}^+ = S_{dp}\s = \begin{pmatrix}
{\traceop\s_p} {K_{\bd\bd}\s}^+ {\traceop\s_p}^T & - {\traceop\s_p} {K_{\bd\bd}\s}^+ K_{\bd d}  \\K_{d \bd}\s{K_{\bd\bd}\s}^+ {\traceop\s_p}^T
& K_{dd}\s -  K_{d\bd}\s{K_{\bd\bd}\s}^+K_{\bd d}\s
\end{pmatrix}
\end{equation}
$S_{dp}\s$ is the operator which associates displacement on the
first subset and reaction on the second subset to given effort on
the first subset and given displacement on the second subset.

As both matrices $K_{\bp\bp}\s$ and $K_{\bd\bd}\s$ may not be
invertible, only their pseudo-inverse has been introduced. The
invertibility is strongly dependent on the choice of interface
subsets.

\comment{If $R_{\bp}\s$ denotes the kernel of matrix $K_{\bp\bp}\s$ then the kernel of $S_{pd}\s$ can be computed:
\begin{equation}
K_{\bp\bp}\s R_{\bp}\s =0 \Longrightarrow S_{pd}\s \begin{pmatrix} K_{p \bp}\s \\ \traceop_d\s
 \end{pmatrix} R_{\bp}\s= 0
\end{equation}
Notice that since $K_{p \bp}\s  R_{\bp}\s$ represents the reaction
of the structure to rigid body motion $R_{\bp}\s$, this term is
often equal to zero (for instance in linear elasticity) though it
may be non-zero in some cases (for instance buckling).}
\end{description}

\subsubsection{Block notations}
While condensed operators simplify the writing of local
operations, the block notations make it easer to understand the
global operations of domain decomposition. We propose to denote by
superscript $\block{\square}$ the row-block repetition of local
vectors and the diagonal-block repetition of matrices, block
assembly operators are written in one row (column-block) and
denoted by special font, for instance:
\begin{gather*}
  \block{u} = \begin{pmatrix} u^{(1)} \\ \vdots \\ u^{(N)}
  \end{pmatrix}\qquad
  \block{f} = \begin{pmatrix} f^{(1)} \\ \vdots \\ f^{(N)}
  \end{pmatrix} \qquad
  \block{\lambda} = \begin{pmatrix} \lambda^{(1)} \\ \vdots \\ \lambda^{(N)} \end{pmatrix}\\
  \block{K} = \begin{pmatrix} K^{(1)} & 0 &\ldots &0 \\0 & \ddots & \ddots &\vdots \\ \vdots &\ddots & \ddots  &0 \\ 0&\ldots & 0& K^{(N)}
  \end{pmatrix} \qquad
  \block{\traceop} =\begin{pmatrix} \tras{1} & 0 &\ldots &0 \\0 & \ddots & \ddots &\vdots \\ \vdots &\ddots & \ddots  &0 \\ 0&\ldots & 0& \tras{N}
  \end{pmatrix}\\
  \plock{\assemop} = \begin{pmatrix} \pass{1} &\ldots & \pass{N} \end{pmatrix}
  \qquad
  \plock{\du{\assemop}} = \begin{pmatrix} \dass{1} &\ldots & \dass{N} \end{pmatrix}
\end{gather*}
\begin{remark}
The specific notation for assembly operators aims at emphasizing
at their specific role in terms of parallelism for the methods.
Moreover, the only operation that requires exchange of data
between subdomains is the use of non-transposed assembly
operators.
\end{remark}
Fundamental system \eqref{eq:ef_dd:1} then reads:
\begin{subequations}\label{eq:block_sys:2}
\begin{equation}\label{eq:block_sys:2a}
  \block{K} \block{u} = \block{f} + \block{\lambda}
  \end{equation}
  \begin{equation}\label{eq:block_sys:2b}
  \plock{\assemop}\block{\traceop}\block{\lambda} = 0
  \end{equation}
  \begin{equation}\label{eq:block_sys:2c}
  \plock{\du{\assemop}}\block{\traceop}\block{u} = \du{0}
  \end{equation}
\end{subequations}
or in condensed form:
\begin{subequations}\label{eq:block_sys:3}
\begin{equation}\label{eq:block_sys:3a}
  \block{S}_p \block{u}_b = \block{b}_p + \block{\lambda}_b
  \end{equation}
  \begin{equation}\label{eq:block_sys:3b}
  \plock{\assemop}\block{\lambda}_b = 0
  \end{equation}
  \begin{equation}\label{eq:block_sys:3c}
  \plock{\du{\assemop}}\block{u}_b = \du{0}
  \end{equation}
\end{subequations}
The orthogonal property of assembly operators \eqref{eq:ortho_int:1} simply reads:
\begin{equation}\label{eq:block_ortho_assem:1}
    \plock{\du{\assemop}}{\plock{\assemop}}^T=0
\end{equation}
Relation \eqref{eq:ortho_int:4} reads:
\begin{equation}\label{eq:block_assem_prim:1}
    \plock{\assemop}\plock{\assemop}^T=\operatorname{diag}(\text{multiplicity})
\end{equation}
\begin{remark}
For improved readability, we will denote by bold font objects
defined in a unique way on the interface (\textit{ie} "assembled"
quantities). Schematically, assembly operators enable to go from
block notations to bold notations and transposed assembly
operators realize the reciprocal operations.
\end{remark}
\subsubsection{Brief review of classical strategies}
We can define general strategies to solve system \eqref{eq:block_sys:2} or  \eqref{eq:block_sys:3}:
\begin{description}
\item[Primal approaches] \cite{LETALLEC:1991:DDM,LETALLEC:1993:DDM,MANDEL:1993:BAL,LETALLEC:1994:DDM,MANDEL:1996:COEF,LETALLEC:1996:BERG,LETALLEC:1997:SHELL}
a unique interface displacement unknown $\ubold_b$ satisfying
equation \eqref{eq:block_sys:3c} is introduced, then an iterative
process enables to satisfy \eqref{eq:block_sys:3b} while always
verifying \eqref{eq:block_sys:3a}.
\item[Dual approaches]
\cite{FARHAT:1991:FETI,FARHAT:1992:SADDLE,FARHAT:1994:RFETI,FARHAT:1994:ADV,FARHAT:1994:OPT,MANDEL:1996:OPT,BHARDWAJ:1999:PQ_LRV}
a unique interface effort unknown $\lambdabold_b$ satisfying
equation \eqref{eq:block_sys:3b} is introduced, then an iterative
process enables to satisfy \eqref{eq:block_sys:3c} while always
verifying \eqref{eq:block_sys:3a}.
\item[Three fields approaches] \cite{BREZZI:1993:3FIELD,PARK:1997:ALG, RIXEN:1999:AFETI}
a unique interface displacement $\ubold_b$ is introduced, then
relation \eqref{eq:block_sys:3c} is dualized so that interface
efforts $\block{\lambda}_b$ are introduced as Lagrange multipliers
which yet have to verify relation \eqref{eq:block_sys:3b}. Then
the iterative process looks simultaneously for
$\left(\block{\lambda}_b,\ubold_b,\block{u}_b\right)$ verifying
exactly equation \eqref{eq:block_sys:3a}. As this method is mostly
designed for non-matching discretizations it will not be exposed
in the remaining of this paper, anyhow a variant of the dual
method which is equivalent to the three-field method with
conforming grids will be described.
\item[Mixed approaches] \cite{GLO90,LAD99b,SERIES:2003:MDD1,SERIES:2003:MDD2}
new interface unknown is introduced which is a linear combination
of interface displacement and effort,
$\block{\mu}_b=\block{\lambda}_b+\block{T}_b\block{u}_b$, then the
interface system is rewritten in terms of unknown $\block{\mu}_b$,
this new system is solved iteratively and then $\block{\lambda}_b$
and $\block{u}_b$ are postprocessed. Of course matrix
$\block{T}_b$ is an important parameter of these methods.
\item[Hybrid approaches] \cite{KLAWONN:1999:DP,FARHAT:2000:FETI_DP,MANDEL:2000:DP,FARHAT:2001:FETI_DP,GOSSELET:2004:DDM}
interface is split into parts where primal, dual or mixed
approaches are applied, specific recondensation methods may then
be applied.
\end{description}
Many different methods can be deduced from these large strategies,
the most common will be presented and discussed in section
\ref{sec:methods}. Anyhow since iterative solvers are used to
solve interface problems, we recommend the reader to refer to
appendix \ref{sec:solvers} where most used solvers are presented,
including important details about constrained resolutions.

\section{Classical solution strategies to the interface problem}\label{sec:methods}
The aim of this section is to give extended review of classical
domain decomposition methods, the principle of which has just been
exposed. The association with Krylov iterative solvers is an
important point of these methods, appendix \ref{sec:solvers}
provides a summary of important results and algorithms that are
used in this section.

\subsection{Primal domain decomposition method}
The principle of primal domain decomposition method is to write
the interface problem in terms of one unique unknown interface
displacement field $\ubold_b$. The trace of local displacement
fields then writes $\block{u}_b=\plock{\assemop}^T\ubold_b$.
Because of the orthogonality between assembly operators
\eqref{eq:block_ortho_assem:1}, equation \eqref{eq:block_sys:3c}
is automatically verified. Using equation \eqref{eq:block_sys:3b}
to eliminate unknown reaction $\block{\lambda}_b$ inside
\eqref{eq:block_sys:3a}, we get the primal formulation of the
interface problem:
\begin{equation}\label{eq:block_primal:1}
  \Sbold_p\ubold_b=\left(\plock{\assemop}\block{S}_p\plock{\assemop}^T\right)\ubold_b =  \plock{\assemop}\block{b}_p = \bbold_p
\end{equation}

Operator $\Sbold_p$ is the global primal Schur complement of the
decomposed structure, it results as the sum of local contributions
(with non-block notations $\Sbold_p=\sum_s\pass{s} S_p\s
{\pass{s}}^T$). Using a direct solver to solve system
\eqref{eq:block_primal:1} implies the exact computation of local
contribution, the sum of these contributions (in a parallel
computing context, this step correspond to large data exchange
between processors) and the inversion of the global primal Schur
complement which size is the global geometric interface (the size
of which is far from being neglectable) and which sparsity is very
poor (each interface degree of freedom is connected to degrees of
freedom belonging  to the same subdomains). Using an iterative
solver is much less expensive since the only required operations
are matrix-vector products which can be realized in parallel
because of the assembled structure of global primal Schur
complement; moreover excellent and rather cheap preconditioner
exists. Note that if global matrix $K$ is symmetric positive
definite then so is operator $\Sbold_p$ and then popular conjugate
gradient algorithm can be used to solve the primal interface
problem, in other cases solvers like GMRes or orthodir have to be
employed.

\subsubsection{Preconditioner to the primal interface problem}
An efficient preconditioner $\tilde{\Sbold}_p^{-1}$ is an
interface operator giving a good approximation of the inverse of
$\Sbold_p$. Various strategies are possible. For instance, a
direct preconditioning method is based on the construction of an
approximate Schur complement from a simplified structure defined
by degrees of freedom "near" the interface. Anyhow such a method
does not respect the repartition of the data through processors. A
good parallel preconditioner has to minimize data exchange.

Since operator $\Sbold_p$ is the sum of local contributions, the
most classical strategy is then to define $\tilde{\Sbold}_p^{-1}$
as a scaled sum of the inverse of local contributions:
\begin{equation}
\tilde{\Sbold}_p^{-1}= \plock{\tilde{\assemop}} {\block{S}_p}^+ \plock{\tilde{\assemop}}^T =\plock{\tilde{\assemop}} {\block{S}_d} \plock{\tilde{\assemop}}^T
\end{equation}
Since ${\block{S}_p}^+=\block{S}_d$ requires the computation of
local problems with given effort on the interface, this
preconditioner is called the Neumann preconditioner. Scaled
assembly operator $\plock{\tilde{\assemop}}$ can be defined the
following way \cite{KLAWONN:2001:FNN}:
\begin{equation}
\plock{\tilde{\assemop}} = \left(\plock{\assemop} \block{M}\plock{\assemop}^T\right)^{-1}\plock{\assemop}\block{M}
\end{equation}
where $\block{M}$ is a parameter which enables to take into
account the heterogeneity of the subdomains connected by the
interface. It should make matrix $\left(\plock{\assemop}
\block{M}\plock{\assemop}^T\right)$ easily invertible and give a
representation of the stiffness of the interface, most commonly:
\begin{itemize}
\item $\block{M}=\block{I}$ for homogeneous structures,
\item $\block{M}=\diag(\block{K}_{bb})$ for compressible heterogeneous structures,
\item $\block{M}=\block{\mu}$ for incompressible heterogeneous structures ($\block{\mu}$
is the diagonal matrix the coefficients of which are the shearing
modulus of interface degrees of freedom).
\end{itemize}

The $(s)$ notation makes it easier to understand implementation of scaled assembly operators:
\begin{equation}
\tilde{\Sbold}_p^{-1}= \sum_s M\s \assemop\s {S\s_d} {\assemop\s}^T M\s
\end{equation}
\begin{itemize}
\item $M\s=\diag(\frac{1}{\text{multiplicity}})$ for homogeneous structures,
\item $M\s=\diag(\frac{\diag(K_{bb}\s)_i}{\sum_j\diag(K_{bb}^{(j)})_i})$
for compressible heterogeneous structures (assuming $i$ represents
the same degree of freedom shared by the $j$ subdomains)
\item $M\s=\diag(\frac{\mu\s_i}{\sum_j{\mu^{(j)}_i}})$
for incompressible heterogeneous structures (assuming $i$
represents the same degree of freedom shared by  the $j$
subdomains)
\end{itemize}

The following partition of unity result clearly holds:
\begin{eqnarray}
\plock{\tilde{\assemop}}\plock{\assemop}^T&=&I_\Upsilon \\
\sum_s M\s & = & I_\Upsilon
\end{eqnarray}

\subsubsection{Coarse problem}
\begin{figure}[ht]\centering
\includegraphics[width=0.8\textwidth]{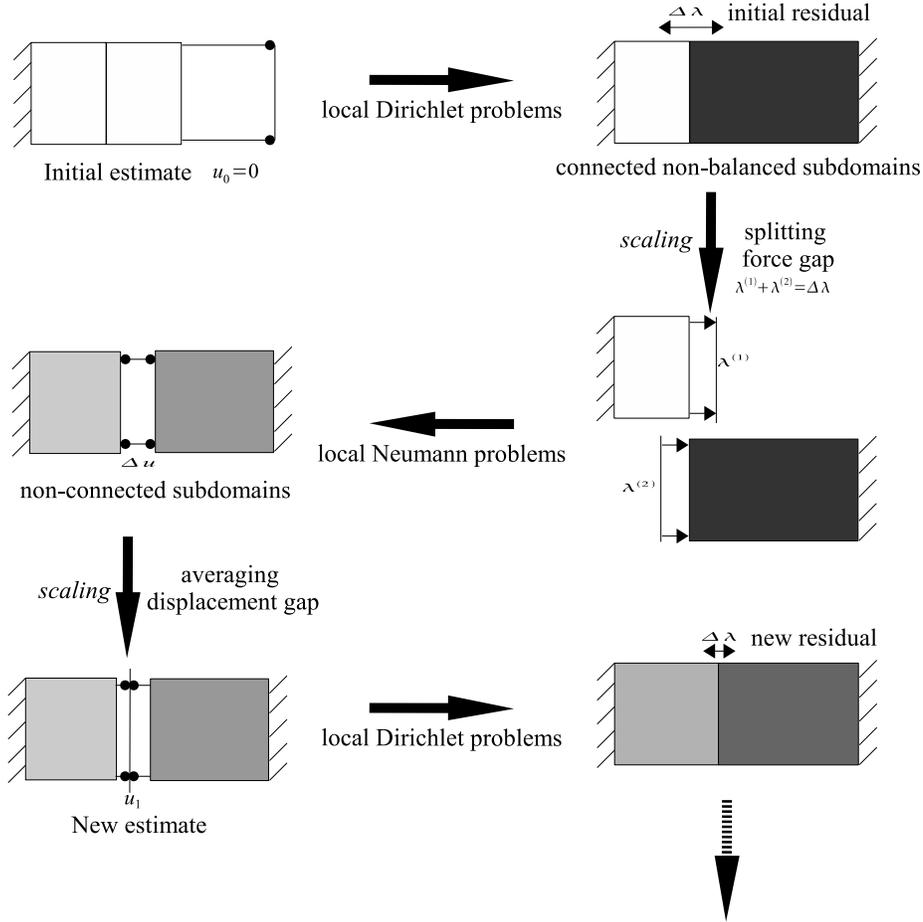}\caption{Representation of first iteration of preconditioned primal approach}\label{fig:sch_primal}
\end{figure}

The use of dual Schur complement is associated to an optimality
condition, as said earlier vector being multiplied by the pseudo
inverse should lie inside the image of $\block{S}_p$. Since
preconditioning is applied to residual $\rbold$, the optimality
condition reads:
\begin{equation}
{\block{R}_b}^T\plock{\tilde{\assemop}}^T \rbold=0
\end{equation}
and introducing classical notation
$\Gbold=\plock{\tilde{\assemop}}\block{R}_b$, $\Gbold^T \rbold=0$.
Such a condition can then be interpreted as an augmented-Krylov
algorithm (see section \ref{sub:augmkry}). Once equipped with that
augmentation problem, the primal Schur complement method is
referred to as the "balanced domain decomposition" (BDD
\cite{MANDEL:1993:BAL,LETALLEC:1994:DDM}). Algorithm
\ref{alg:bdd:1} summarizes the classical BDD approach, and figure
\ref{fig:sch_primal} provides a schematic representation of the
first iteration of the preconditioned primal approach.
\begin{algorithm}\caption{Primal Schur complement with conjugate gradient}\label{alg:bdd:1}
\begin{algorithmic}[1]
\STATE Set $P=I-G(G^T S_p G)^{-1} G^T S_p$
\STATE Compute $u_0 = G (G^T S_p G)^{-1} G^T b_p$
\STATE Compute $r_0=b_p-S_p u_0=P^T b_p$
\STATE $z_0 = \tilde{S}_p^{-1}r_0$ set $w_0=z_0$%
\FOR{$j=0,\ldots,m$} %
  \STATE $p_j = S_p P w_j $ (notice $S_p P = P^T S_p = P^T S_p P$ )
  \STATE $\alpha_j=(z_j,r_j)/(p_j,w_j)$%
  \STATE $u_{j+1}=u_j+\alpha_j w_j$%
  \STATE $r_{j+1}=r_j-\alpha_j p_j$%
  \STATE $z_{j+1}=\tilde{S}_p^{-1}r_{j+1}$
  \STATE For $0\leqslant i\leqslant j$, $\beta_j^i=-(z_{j+1},p_i)/(w_i,p_i)$
  \STATE $w_{j+1}=z_{j+1}+\sum_{i=1}^j {\beta_j^i w_i}$%
\ENDFOR
\end{algorithmic}
\end{algorithm}

\subsubsection{Error estimate}

The reference error estimate is the one linked to the convergence
over the complete structure: $\frac{\|Ku-f\|}{\|f\|}$. Assuming
local inversions are exact, we reach the following result:
\begin{equation}\label{eq:res_pri:1}
\frac{\|Ku-f\|}{\|f\|} = \frac{\|\Sbold_p \ubold_b - \bbold_p\|}{\|f\|}
\end{equation}
During the iterative process $\|\Sbold_p \ubold_b - \bbold_p\|$ is
the norm of residual $\rbold$ as computed line 9 of algorithm
\ref{alg:bdd:1}, so the global convergence can be controlled by
the convergence of the interface iterative process.

\subsubsection{P-FETI method}
The P-FETI method is a variation of BDD proposed by
\cite{FRAGAKIS:2004:MHP1,FRAGAKIS:2004:MHP2} inspired by the dual
approach (the reader should refer to the dual method before going
further inside P-FETI). Its principle is to provide another
assembly operator which incorporate rigid body elimination by a
dual-like projector.
\begin{eqnarray}
\tilde{\Sbold}_p ^{-1} &=& \plock{H} \block{S}_d \plock{H}^T \\
 \plock{H}^T & = & \plock{\tilde{\assemop}}^T - \plock{\du{\assemop}}^T \du{\Qbold} \du{\Gbold} \left( \du{\Gbold}^T\du{\Qbold} \du{\Gbold} \right)^{-1}\Gbold^T
\end{eqnarray}
The choice of matrix $\du{\Qbold}$ is guided by the same
considerations as in the dual method. It is worth noting that when
$\du{\Qbold}$ is chosen equal to the Dirichlet preconditioner of
the dual method ($\du{\Qbold}=\plock{\du{\tilde{\assemop}}}^T
\block{S}_p\plock{\du{\tilde{\assemop}}}$) then the P-FETI method
is equivalent to the classical balanced domain decomposition.

\subsection{Dual domain decomposition method}
The principle of dual domain decomposition method is to write the
interface problem in terms of one unique unknown interface effort
field $\du{\lambdabold}_b$. The trace of local effort fields then
writes
$\block{\lambda}_b=\plock{\du{\assemop}}^T\du{\lambdabold}_b$.
Because of the orthogonality between assembly operators
\eqref{eq:block_ortho_assem:1}, equation \eqref{eq:block_sys:3b}
is automatically verified. In order to eliminate unknown interface
displacement field using \eqref{eq:block_sys:3c}, we first obtain
it from equation \eqref{eq:block_sys:3a} (or equivalently
\eqref{eq:block_sys:2a}): as seen in \eqref{eq:schurd:2} the
inversion of local systems may require the use of generalized
inverse and the introduction of rigid body motions the magnitude
of which is denoted by vector $\alpha\s$, the use of generalized
inverse is then submitted to compatibility condition.
\begin{eqnarray}
\block{u}_b & = &  \block{S}_d (\block{b}_p + \plock{\du{\assemop}}^T\du{\lambdabold}_b) + \block{R}_b \block{\alpha} \\
{\block{R}_b}^T (\block{b}_p + \plock{\du{\assemop}}^T\du{\lambdabold}_b)& = & 0
\end{eqnarray}
The first line is then premultiplied by $\plock{\du{\assemop}}$ (same expressions could be obtain from non condensed notations).
\begin{gather*}
\du{\Sbold}_d  = \plock{\du{\assemop}} \block{S}_d \plock{\du{\assemop}}^T \\
\block{b}_d  =   \block{S}_d \block{b}_p = \block{\traceop}{\block{K}}^+ \block{f} \\
\du{\Gbold} =\plock{\du{\assemop}} \block{R}_b \\
\block{e} = {\block{R}_b}^T\block{b}_p = {\block{R}}^T \block{f}
\end{gather*}
we get the dual formulation of the interface problem:
\begin{equation}\label{eq:block_dual:1}
  \begin{pmatrix} \du{\Sbold}_d & \du{\Gbold} \\ \du{\Gbold}^T & 0 \end{pmatrix}
  \begin{pmatrix} \du{\lambdabold}_b \\ \block{\alpha} \end{pmatrix} =
  \begin{pmatrix} - \bbold_d \\ -\block{e}
  \end{pmatrix}
\end{equation}
This is the basic dual Schur complement method, also called Finite
Element Tearing and Interconnecting method (FETI
\cite{FARHAT:1991:FETI,FARHAT:1994:ADV}). For similar reasons to
the primal Schur complement method, this system is most often
solved using an iterative solver, then we will soon discuss the
preconditioning issue and how the
$\du{\Gbold}^T\du{\lambdabold}_b+\block{e}=0$ constraint is
handled. Let us first remark that global dual Schur complement
$\du{\Sbold}_d$ is non-definite as soon as redundancies appear in
the connectivity description of the interface, anyhow it is easy
to prove \cite{FARHAT:1991:FETI} that local contributions
$\block{\lambda}_b=\plock{\du{\assemop}}^T \du{\lambdabold}_b$ are
unique (non-definition only affect the "artificial" splitting of
forces on multiple points), and that because the right hand side
lies in $\range(\plock{\du{\assemop}})$ iterative process
converges; other considerations on the splitting of physical
efforts between subdomains will lead to improved initialization
(see section \ref{subsub:fetiini} and \cite{GOSSELET:2003:IEI}).
\subsubsection{Preconditioner to the dual interface problem}
Like it is done in the primal approach, the most interesting
preconditioners are researched as assembly of local contributions,
and the global dual Schur complement being a sum of local contributions,
optimal preconditioner is a scaled sum of local inverses.
\begin{equation}\label{eq:precdual:1}
\du{\tilde{\Sbold}}_d^{-1} = \plock{\tilde{\du{\assemop}}} {\block{S}_d}^+ \plock{\tilde{\du{\assemop}}}^T =\plock{\tilde{\du{\assemop}}} {\block{S}_p} \plock{\tilde{\du{\assemop}}}^T
\end{equation}
Because this preconditioner uses local primal Schur complement,
which corresponds to the local resolution of imposed displacement
problems, it is commonly called the Dirichlet preconditioner. One
interesting point is the possibility to give approximation of the
local Schur complement operator leading to the following
preconditioners:
\begin{eqnarray}
\block{S}_p \approx \block{K}_{bb} & \qquad & \text{lumped\ preconditioner} \\
\block{S}_p \approx \diag(\block{K}_{bb}) & \qquad & \text{superlumped\ preconditioner}
\end{eqnarray}
These preconditioners have very low computational cost (they do
not require the computation and storage of the inverse of local
internal matrices ${\block{K}_{ii}}^{-1}$), even if their
numerical efficiency is not as strong as the Dirichlet
preconditioner, they can lead to very reduced computational time.

Scaled assembly operator $\plock{\tilde{\du{\assemop}}}$ can be defined the following way \cite{KLAWONN:2001:FNN}:
\begin{equation}
\plock{\tilde{\du{\assemop}}} = \left(\plock{\du{\assemop}} {\block{M}}^{-1}\plock{\du{\assemop}}^T\right)^{+}\plock{\du{\assemop}}{\block{M}}^{-1}
\end{equation}
where $\block{M}$ is the same parameter as for the primal
approach. Such a definition is not that easy to implement, an
almost equivalent strategy is then used, easily described using
the $(s)$ notation:
\begin{equation}
\tilde{\Sbold}_d^{-1}= \sum_s \du{M}\s \du{\assemop}\s {S\s_p} {\du{\assemop}\s}^T \du{M}\s
\end{equation}
\begin{itemize}
\item $\du{M}\s=\diag(\frac{1}{\text{multiplicity}})$ for homogeneous structures,
\item $\du{M}\s=\diag(\frac{\diag(K_{bb}^{(r)})_i}{\sum_j\diag(K_{bb}^{(j)})_i})$ for compressible heterogeneous structures (assuming $i$ represents the same degree of freedom shared by the $j$ subdomains and $(r)$ is the subdomain connected to $(s)$ on degree of freedom $i$),
\item $\du{M}\s=\diag(\frac{\mu^{(r)}_i}{\sum_j{\mu^{(j)}_i}})$ for incompressible heterogeneous structures (assuming $i$ represents the same degree of freedom shared by  the $j$ subdomains and $(r)$ is the subdomain connected to $(s)$ on degree of freedom $i$).
\end{itemize}
We have the following partition of unity result:
\begin{equation}
\plock{\tilde{\du{\assemop}}}\plock{\du{\assemop}}^T=I_{\du{\Upsilon}}
\end{equation}
and the following complementarity between primal and dual scalings \cite{GOSSELET:2003:IEI}:
\begin{eqnarray}\label{eq:scalings:1}
\plock{\tilde{\assemop}}^T\plock{\assemop}+\plock{\du{\assemop}}^T\plock{\du{\tilde{\assemop}}}&=&\block{I}
\\{\pass{s}}^TM\s\pass{s}+{\dass{s}}^T\du{M}\s\dass{s}&=&I_{\Upsilon\s}
\end{eqnarray}

\subsubsection{Coarse problem}
Admissibility condition
$\du{\Gbold}^T\du{\lambdabold}_b+\block{e}=0$, can be handled
with an initialization / projection algorithm (see section
\ref{sub:conskry}):
$\du{\lambdabold}_b=\du{\lambdabold}_0+\du{\Pbold}\du{\lambdabold}^*$
with $\du{\Gbold}^T\du{\lambdabold}_0=-\block{e}$ and
$\du{\Gbold}^T \du{\Pbold}=0$.
\begin{eqnarray}
\du{\lambdabold}_0 & = & -\du{\Qbold} \du{\Gbold} \left(\du{\Gbold}^T \du{\Qbold} \du{\Gbold}\right)^{-1} \block{e} \\
\du{\Pbold} & = & \du{\Ibold} - \du{\Qbold} \du{\Gbold} \left(\du{\Gbold}^T \du{\Qbold} \du{\Gbold}\right)^{-1} \du{\Gbold}^T
\end{eqnarray}
The easiest choice for operator $\du{\Qbold}$ is the identity
matrix, projector $\du{\Pbold}$ is then orthogonal, this choice is
well suited to homogeneous structures. For heterogeneous
structures, matrix $\du{\Qbold}$ has to provide information on the
stiffness of subdomains, then $\du{\Qbold}$ is chosen to be a
version of the preconditioner leading to "superlumped projector"
($\du{\Qbold}=\plock{\tilde{\du{\assemop}}}\operatorname{diag}(\block{K}_{bb})\plock{\tilde{\du{\assemop}}}^T$),
"lumped projector"
($\du{\Qbold}=\plock{\tilde{\du{\assemop}}}\block{K}_{bb}\plock{\tilde{\du{\assemop}}}^T$)
and "Dirichlet projector"
($\du{\Qbold}=\plock{\tilde{\du{\assemop}}}\block{S}_{p}\plock{\tilde{\du{\assemop}}}^T$).
Superlumped projector is often a good compromise between numerical
efficiency and computational cost.

\begin{figure}[ht]\centering
\includegraphics[width=0.8\textwidth]{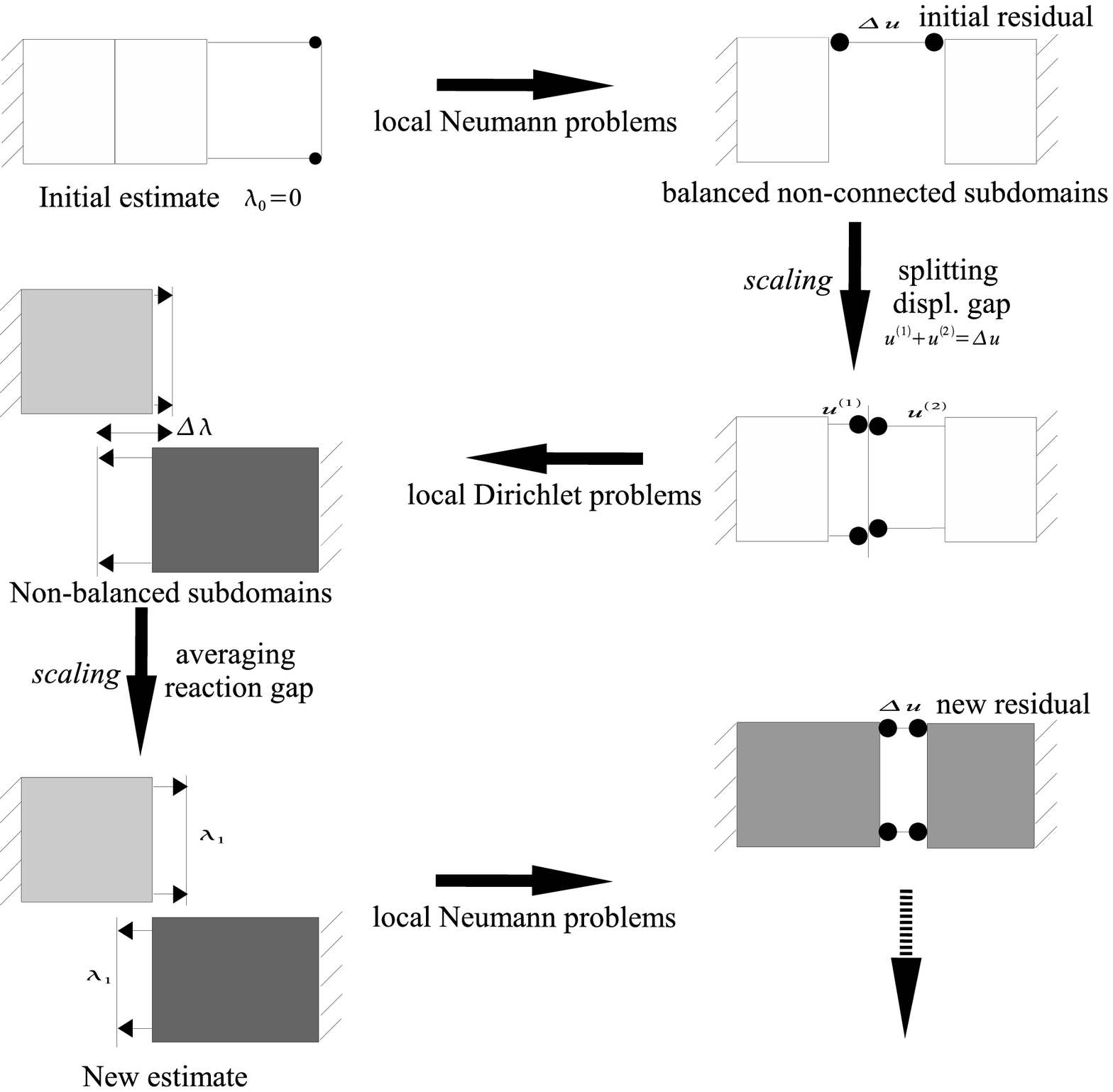}\caption{Representation of first iteration of preconditioned dual approach}\label{fig:sch_dual}
\end{figure}

Algorithm \ref{alg:feti:1} presents a classical implementation of
FETI method, and figure \ref{fig:sch_dual} provides a schematic
representation of the first iteration of the preconditioned dual
approach.

\begin{algorithm}\caption{Dual Schur complement with conjugate gradient}\label{alg:feti:1}
\begin{algorithmic}[1]
\STATE Set $P={I}- {Q} {G}({G}^T {Q} {G})^{-1} {G}^T $
\STATE Compute ${\lambda}_0 = -{Q} {G} ({G}^T {Q} {G})^{-1} e$
\STATE Compute ${r}_0=P^T{b}_d-S_d \lambda_0)$
\STATE $z_0 = P\tilde{S}_d^{-1}r_0$ set $w_0=z_0$%
\FOR{$j=0,\ldots,m$} %
  \STATE $p_j = P^T S_d w_j $
  \STATE $\alpha_j=(z_j,r_j)/(p_j,w_j)$%
  \STATE $\lambda_{j+1}=\lambda_j+\alpha_j w_j$%
  \STATE $r_{j+1}=r_j-\alpha_j p_j$%
  \STATE $z_{j+1}=P \tilde{S}_d^{-1}r_{j+1}$
  \STATE For $0\leqslant i\leqslant j$, $\beta_j^i=-(z_{j+1},p_i)/(w_i,p_i)$
  \STATE $w_{j+1}=z_{j+1}+\sum_{i=1}^j {\beta_j^i w_i}$%
\ENDFOR
\STATE $\block{\alpha}= (G^T Q G)^{-1} G^t r_m$
\STATE $\block{u}= {\block{K}}^+\lambda_m + \block{R}\block{\alpha}$
\end{algorithmic}
\end{algorithm}

\subsubsection{Error estimate}\label{subsub:error_dual}

The convergence of the dual domain decomposition method is
strongly linked to physical considerations. After projection, the
residual can be interpreted as the jump of displacement between
substructures:
\begin{eqnarray}\label{eq:dualresid}
\du{\rbold} & = & \du{\Pbold}^T (-\du{\bbold}_d - \du{\Sbold}_d \du{\lambdabold}) = \plock{\du{\assemop}} \block{u} = \du{\Deltabold}(u) \\
\du{\Deltabold}(u)_{|\Upsilon^{(i,j)}}&=&u^{(i)}_{|\Upsilon^{(i,j)}}-u^{(j)}_{|\Upsilon^{(i,j)}}
\end{eqnarray}
Anyhow, such an interpretation cannot be connected to the global
convergence of the system. In order to evaluate the global
convergence, a unique interface displacement field has to be
defined (most often using a scaled average of local displacement
fields) and used to evaluate the global residual. When using the
Dirichlet preconditioner, it is possible to cheaply evaluate that
convergence criterion. Average interface displacement $\ubold_b$
can be defined as follow:
\begin{equation}
\ubold_b= \plock{\assemop}\plock{\du{\tilde{\assemop}}}^T\du{\Deltabold}
\end{equation}
then, from equation \eqref{eq:res_pri:1}, convergence criterion
reads: $\|K u -f \|=\| \plock{\assemop} \block{S}_p
\plock{\du{\tilde{\assemop}}}^T\du{\rbold} \|$. So when using the
Dirichlet preconditioner, the evaluation of the global residual
only requires the use of a geometric assembly after the local
Dirichlet resolution.

\subsubsection{Interpretation and improvement on the initialization}\label{subsub:fetiini}
Let us come back to the original dual system \eqref{eq:block_sys:2}.
\begin{equation}\label{eq:FETINI_FETI:1}
\left\{ \begin{array}{l}  \block{K} \block{u} = \block{f} +{\block{\traceop}}^T
{\plock{\du{\assemop}}}^T\du{\lambdabold}_b \\
  \plock{\du{\assemop}}\block{\traceop}\block{u} = \du{0} \end{array}\right.
\end{equation}
And suppose this system is being initialized with non zero effort ${\du{\lambdabold}_b}_0$:
\begin{eqnarray}\label{eq:FETINI_FETI:2}
\block{K} \block{u}& =&   \block{f} + {\block{\traceop}}^T {\plock{\du{\assemop}}}^T\du{\lambdabold}_b\\
\nonumber \text{let}\ \du{\lambdabold}_b &=& \du{\widetilde{\lambdabold}}_b + {\du{\lambdabold}_b}_0\\
\nonumber \block{K} \block{u} & =& \block{f} +
{\block{\traceop}}^T {\plock{\du{\assemop}}}^T \du{\widetilde{\lambdabold}}_b +{\block{\traceop}}^T {\plock{\du{\assemop}}}^T {\du{\lambdabold}_b}_0 \\
& =& \block{\widetilde{f}} + {\block{\traceop}}^T
{\plock{\du{\assemop}}}^T\du{\widetilde{\lambdabold}}_b \\
\text{with\ } \block{\widetilde{f}} &=& \block{f} +
{\block{\traceop}}^T{\plock{\du{\assemop}}}^T
{\du{\lambdabold}_b}_0
\end{eqnarray}
So initialization  ${\du{\lambdabold}_b}_0$ can be interpreted as modification
${\block{\traceop}}^T{\plock{\du{\assemop}}}^T
{\du{\lambdabold}_b}_0$ of the intereffort between substructures:
local problems are defined except for an equilibrated interface effort field;
 the only field that makes mechanical sense (and that is uniquely defined) is the assembly of interface efforts.
\begin{gather}\label{eq:FETIINI_FETI:2a}
\plock{\assemop}\block{\traceop}\block{\widetilde{f}} =
\plock{\assemop}\block{\traceop}\block{f}=\fbold_b\text{\ \ \ global\
interface\ effort}\\
\text{because}\
\plock{\assemop}\block{\traceop}{\block{\traceop}}^T{\plock{\du{\assemop}}}^T
{\du{\lambdabold}_b}_0 = \plock{\assemop}
{\plock{\du{\assemop}}}^T {\du{\lambdabold}_b}_0 = 0
\end{gather}

Non-zero initialization then can be interpreted as a repartition
of global interface effort $\fbold_b$. Two strategies can be
defined in order to realize that splitting.
\begin{description}
\item[Classical effort splitting]
Though splitting is hardly ever interpreted as a specific
initialization, it is commonly realized, based on the
difference of stiffness between neighboring substructures (that
idea is strongly connected to the definition of scaled assembly
operators): the aim is to guide the stress flow inside the stiffer
substructure, sticking to what mechanically occurs.

Global interface effort $\fbold_b$ is then split according to the
stiffness scaling ($\block{M}=\diag(\block{K}_{bb})$), which leads
to modified local effort $\block{\tilde{f}}_b$.
\begin{equation}\label{eq:repart:1}
\block{\widetilde{f}}_b =\block{M}\plock{\assemop}
\left(\plock{\assemop}\block{M}{\plock{\assemop}}^T
\right)^{-1} \fbold_b
\end{equation}
Complete effort $\block{\widetilde{f}}$ is constituted by
$\block{f}$ inside the substructure
($(I-{\block{\traceop}}^T\block{\traceop})\block{f}$) and split effort on its interface
(${\block{\traceop}}^T\block{\widetilde{f}}_b$).
\begin{equation}\label{eq:repart:1b}
\block{\widetilde{f}}=
(I-{\block{\traceop}}^T\block{\traceop})\block{f} +
{\block{\traceop}}^T\block{\widetilde{f}}_b
\end{equation}
Because of the complementarity between scaled assembly operators \eqref{eq:scalings:1}, final effort reads
\begin{equation}\label{eq:repart:2}
\block{\widetilde{f}}=\block{f}-{\block{\traceop}}^T{\plock{\du{\assemop}}}^T
\left(\plock{\du{\assemop}}{\block{M}}^{-1}{\plock{\du{\assemop}}}^T
\right)^{+}\plock{\du{\assemop}}{\block{M}}^{-1}{\block{\traceop}}\block{f}
\end{equation}
\item[Interface effort splitting] \cite{GOSSELET:2003:IEI}
If we start from condensed dual system \eqref{eq:block_sys:3}
\begin{equation}\label{eq:feti_condens:2}\left\{\begin{array}{l}
    \block{S}_p \block{u}_b=\block{b}_p+{\plock{\du{\assemop}}}^T\du{\lambdabold}_b\\
    \plock{\du{\assemop}} \block{u}_b = 0 \end{array}\right.
\end{equation}
condensed efforts can be split along the interface as long as
global condensed effort remains unique. Assembled condensed
interface effort reads $\bbold_p={\plock{\assemop}\block{b}_p}$,
if it is split according to the stiffness of the substructures:
\begin{equation}
\block{\widetilde{b}}_p = \block{M}\plock{\assemop}
\left(\plock{\assemop}\block{M}{\plock{\assemop}}^T
\right)^{-1} \bbold_p
\end{equation}
We have, using the complementarity between scalings:
\begin{equation}\label{eq:repart:3}
\block{\widetilde{b}}_p=\block{b}_p-{\plock{\du{\assemop}}}^T\left(\plock{\du{\assemop}}
{\block{M}}^{-1}{\plock{\du{\assemop}}}^T
\right)^{+}\plock{\du{\assemop}}{\block{M}}^{-1}\block{b}_p
\end{equation}
Or in a non-condensed form:
\begin{equation}\label{eq:repart:4}
    \block{\widetilde{f}}=\block{f}-{\block{\traceop}}^T{\plock{\du{\assemop}}}^T
    \left(\plock{\du{\assemop}}{\block{M}}^{-1}{\plock{\du{\assemop}}}^T
\right)^{+}\plock{\du{\assemop}}{\block{M}}^{-1}\block{b}_p
\end{equation}
\end{description}
As will be shown in assessments, the classical splitting leads to
almost no improvement of the method while the condensed splitting
can be very efficient for heterogeneous structures. In fact the
initialization associated to this splitting can be proved to be
optimal in a mechanical sense; it can also be obtained from the
assumptions used for the primal approach.

Primal approach initialization is realized supposing that
interface displacement field is zero on the condensed problem;
from \eqref{eq:feti_condens:2} we get:
\begin{equation}{\plock{\du{\assemop}}}^T{\du{\lambdabold}_b}_0+\block{b}_p\simeq 0\end{equation}
which could only be the solution if null interface displacement
was the solution. Then local interface efforts are split into an
equilibrated part and its remaining $\block{\varrho}$:
\begin{equation}\left\{\begin{array}{l}
    \block{b}_p={\plock{\du{\assemop}}}^T \du{\gammabold} + \block{\varrho}\\ \du{\gammabold}= \left(\plock{\du{\assemop}}{\block{D}}{\plock{\du{\assemop}}}^T
\right)^{+}\plock{\du{\assemop}}{\block{D}} \block{b}_p
\end{array}\right.
\end{equation}
$\block{D}$ is a symmetric definite matrix, remaining
$\block{\varrho}$ is orthogonal to
$\range(\block{D}\plock{\du{A}})$. If the system is initialized
by:
\begin{equation}
    {\du{\lambdabold}_b}_{00}=-\du{\gammabold}
\end{equation}
then initial residual
${\plock{\du{\assemop}}}^T{\du{\lambdabold}_b}_{00}+\block{b}_p=-\block{\varrho}$
is minimal in the sense of the norm associated to $\block{D}$. If
$\block{D}=\operatorname{diag}(\block{\K{bb}})^{-1}$ then
initialization is equivalent to the splitting of condensed efforts
according to the stiffness of substructures ;
$\operatorname{diag}(\block{\K{bb}})$ being an approximation of
$\block{S_p}$ that norm can be interpreted as an energy.

The initialization by the splitting of condensed efforts has to be made compatible with solid body motions by the computation of:
\begin{equation}
    {\du{\lambdabold}_b}_0=\du{\Pbold}{\du{\lambdabold}_b}_{00}-\du{\Qbold}\du{\Gbold}\left(\du{\Gbold}^T\du{\Qbold}\du{\Gbold}\right)^{-1}\block{e}
\end{equation}

\begin{remark}
If $\block{D}=\block{S}_d$ was not computationally too expensive
then improved initialization with Dirichlet preconditioner would
lead to immediate convergence.
\end{remark}
\begin{remark}
The recommended choice
$\block{D}=\operatorname{diag}(\block{\K{bb}})^{-1}$ is
computationally very cheap, the heaviest operation is the
computation of condensed efforts (one application of Dirichlet
operator). Then if the Dirichlet preconditioner is used, new
initialization is just as expensive as one preconditioning step
but it can lead to significant reduction of iterations, so it
should be employed. Of course if light preconditioner is preferred
classical splitting should be used.
\end{remark}

\subsection{Three fields method /  A-FETI method}\label{sub:AFETI}
The A-FETI method \cite{PARK:1997:PERF,PARK:1997:ALG} can be
explained as the application of the three-field strategy
\cite{BREZZI:1993:3FIELD} to conforming grids, this method is
widely studied in \cite{RIXEN:1999:AFETI}. Back to
\eqref{eq:block_sys:3}, we have
\begin{eqnarray}
  \block{S}_p \block{u}_b &=& \block{b}_p + \block{\lambda}_b\\
  \plock{\assemop}\block{\lambda}_b &= &0\\
  \plock{\du{\assemop}}\block{u}_b &=& \du{0}
\end{eqnarray}
A-FETI method is based, like in the primal approach, on the
introduction of unknown interface displacement field $\ubold_b$,
the continuity of displacement then reads:
\begin{equation}
\block{u}_b = \plock{\assemop}^T \ubold_b
\end{equation}
but local displacements are not eliminated like in the primal
approach, complete system reads:
\begin{equation}
\begin{pmatrix}
\block{S}_p & -I & 0 \\
-I & 0 & \plock{\assemop}^T \\
0  & \plock{\assemop} & 0
\end{pmatrix}
\begin{pmatrix}
\block{u}_b \\ \block{\lambda}_b \\ \ubold_b
\end{pmatrix} = \begin{pmatrix} \block{b}_p \\ 0 \\ 0 \end{pmatrix}
\end{equation}
In order to eliminate interface displacement $\ubold_b$ a specific symmetric projector is introduced:
\begin{equation}
B = I - \plock{\assemop}^T \left(\plock{\assemop}\plock{\assemop}^T\right)^{-1} \plock{\assemop}
\end{equation}
$B$ realizes the orthogonal projection on $\ker(\plock{\assemop})$
($\plock{\assemop}B=0$). Since
$\plock{\assemop}\block{\lambda}_b=0$ then $\block{\lambda}_b$ can
be written as
\begin{equation}
\block{\lambda}_b= B \block{\mu}_b
\end{equation}
$\block{\mu}_b$ is a new interface effort, corresponding (recall
$\plock{\assemop}\plock{\assemop}^T=\operatorname{diag}(\text{multiplicity})$)
to an average of original $\block{\lambda}_b$. Introducing last
result and using $B^T\plock{\assemop}^T=0$ to eliminate interface
displacement we have
\begin{equation}
\begin{pmatrix}
\block{S}_p & -B  \\
-B^T & 0
\end{pmatrix}
\begin{pmatrix}
\block{u}_b \\ \block{\mu}_b
\end{pmatrix} = \begin{pmatrix} \block{b}_p \\ 0 \end{pmatrix}
\end{equation}
Then using classical elimination of local displacement by the
inversion of the first line of the previous system, we get
\begin{eqnarray}
\block{u}_b &=& {\block{S}_p}^+ \left(\block{b}_p + B \block{\mu}_b \right) + \block{R}_b\block{\alpha} \\
{\block{R}_b}^T\left(\block{b}_p + B \block{\mu}_b \right) &=&0
\end{eqnarray}
which leads to
\begin{equation}
\begin{pmatrix}
B^T\block{S}_dB &  B^T\block{R}_b  \\
 {\block{R}_b}^T B & 0
\end{pmatrix}
\begin{pmatrix}
\block{\mu}_b \\ \block{\alpha}
\end{pmatrix} = \begin{pmatrix} -B^T{\block{S}}_p^+\block{b}_p \\ -{\block{R}_b}^T\block{b}_p \end{pmatrix}
\end{equation}
This system is very similar to the classical dual approach system,
and in consequence is solved in the same way (using projected
algorithm). Anyhow the main difference is that Lagrange multiplier
$\block{\mu}_b$ is defined locally on each subdomain and not
globally on the interface.

It was proved in \cite{RIXEN:1999:AFETI} that A-FETI is
mathematically equivalent to classical FETI with special choice of
the $\du{Q}$ matrix parameter of the rigid body motion projector.
In fact if $\du{Q}=\diag(\frac{1}{\text{multiplicity}})$ then FETI
leads to the same iterates as A-FETI. Moreover operator $B$ is an
orthonormal projector which realizes the interface equilibrium of
local reactions $\block{\mu}_b$, it can be analyzed as an
orthonormal assembly operator as described in figure
\ref{fig:omegef:3}.

To sum up, A-FETI can be viewed as the conforming grid version of
the three-field approach, a specific case of classical FETI, and a
dual approach with non-redundant description of the connectivity
interface with orthonormal assembly operator.

\subsection{Mixed domain decomposition method}
Mixed approaches offer a rich framework for domain decomposition
methods. It enables to give a strong mechanical sense to the
method, mostly by providing a behavior to the interface. The mixed
approach is one of the bases of the LaTIn method
\cite{LAD99b,LAD01,NOUY:2003:THES,LADEVEZE.2006.1}, a very interesting strategy designed for
nonlinear analysis; as we have restrained our paper to linearized
problems, we do not go further inside this method which would deserve
extended survey. Several studies were realized on mixed
approaches, these methods possess strong similarities, we here
mostly refer to works on so-called "FETI-2-fields" method
\cite{SERIES:2003:MDD1,SERIES:2003:MDD2}.

The principle of the method is to rewrite the interface conditions:
\begin{equation}\label{eq:interface:1}
\left\{\begin{array}{l} \plock{\assemop}\block{\lambda}_b=0 \\ \plock{\du{\assemop}}\block{u}_b=\du{0} \end{array}\right.
\end{equation}
in terms of a new local interface unknown, which is a linear
combination of interface effort and displacement.
\begin{equation}\label{eq:interface:2}
\block{\mu}_b = \block{\lambda}_b + \block{T}_b\block{u}_b
\end{equation}
$\block{\mu}_b$ is homogeneous to an effort and $\block{T}_b$ can
be interpreted as an interface stiffness. Mixed methods thus
enable to give a mechanical behavior to the interface, in our case
(perfect interfaces) it can be mechanically interpreted as the
insertion of springs to connect substructures. New interface
condition reads:
\begin{eqnarray}
\plock{\assemop}^T\plock{\assemop}\block{\lambda}_b+\block{T}_b\plock{\du{\assemop}}^T\plock{\du{\assemop}}\block{u}_b&=&\block{0}\\
\text{or}\qquad \plock{\assemop}^T\plock{\assemop}\block{\mu}_b-\left(\plock{\assemop}^T\plock{\assemop}\block{T}_b-\block{T}_b\plock{\du{\assemop}}^T\plock{\du{\assemop}}\right)\block{u}_b&=&\block{0}\label{eq:mixt_int:1}
\end{eqnarray}
It is of course necessary to study the condition for this system
being equivalent to system \eqref{eq:interface:1}. It is important
to note that two conditions lying on the global interfaces
(geometric and connectivity) were traced back to the local
interfaces, so up to a zero-measure set (multiple points) the
conditions have the same dimension. The new condition is
equivalent to the former if facing local interfaces do not hold
the same information which is the case if matrix
$\left(\plock{\assemop}^T\plock{\assemop}\block{T}_b -
\block{T}_b\plock{\du{\assemop}}^T\plock{\du{\assemop}}\right)$ is
invertible. An easy method to construct such matrices will be soon
discussed

If unknown $\block{\mu}_b$ is introduced inside local equilibrium equation, the local system reads:
\begin{equation}
\left(\block{S}_p+\block{T}_b\right)\block{u}_b = \block{b}_p+\block{\mu}_b
\end{equation}
If we assume that $\block{T}_b$ is chosen so that $\left(\block{S}_p+\block{T}_b\right)$ is invertible then we have:
\begin{equation}
\block{u}_b=\left(\block{S}_p+\block{T}_b\right)^{-1}\left(\block{\mu}_b+\block{b}_p\right)
\end{equation}
Then substituting this expression inside interface condition \eqref{eq:mixt_int:1}, interface system reads:
\begin{multline}\label{eq:mixt_sys:1}
\left( \plock{\assemop}^T\plock{\assemop}-\left(\plock{\assemop}^T\plock{\assemop}\block{T}_b-\block{T}_b\plock{\du{\assemop}}^T\plock{\du{\assemop}}\right)\left(\block{S}_p+\block{T}_b\right)^{-1}\right)\block{\mu}_b \\ =\left(\plock{\assemop}^T\plock{\assemop}\block{T}_b-\block{T}_b\plock{\du{\assemop}}^T\plock{\du{\assemop}}\right)\left(\block{S}_p+\block{T}_b\right)^{-1}\block{b}_p
\end{multline}
so mixed approaches have the originality to rewrite global
interface conditions on the local interfaces and to look for
purely local unknown (which means that the size of the unknown is
about twice the size of the unknown in classical primal or dual
methods).

This general scheme for mixed methods has, as far as we know,
never been employed. A first reason is that it leads to certain
programming complexity, second the manipulation of zero-measure
interfaces is not easy for methods aiming at introducing strong
mechanical sense and hard to justify from a mathematical point of
view. So most often a simplified method is preferred, which takes
only into account non-zero-measure interfaces. Such an approach
simplifies the connectivity description of the interface, every
relationship on the interface only deals with couples of
subdomains. In order to have the clearer expression possible, we
present the algorithm in the two subdomains case. Interface
equilibrium reads:
\begin{equation}
\left\{\begin{array}{l}
u^{(1)}_b-u^{(2)}_b=0 \\
\lambda^{(1)}_b+\lambda^{(2)}_b=0
\end{array}\right.
\end{equation}
which is equivalent to
\begin{equation}
\left\{\begin{array}{l}
\lambda^{(1)}_b+\lambda^{(2)}_b+T^{(1)}\left(u^{(1)}_b-u^{(2)}_b\right)=0\\
\lambda^{(1)}_b+\lambda^{(2)}_b+T^{(2)}\left(u^{(2)}_b-u^{(1)}_b\right)=0\\
\end{array}\right.
\end{equation}
under the condition of invertibility of $\left(T^{(1)}+T^{(2)}\right)$. Introducing unknown $\mu_b\s=\lambda_b\s+T\s u\s$ interface system reads:
\begin{equation}
\left\{\begin{array}{l}
\mu^{(1)}_b+\mu^{(2)}_b-\left(T^{(1)}+T^{(2)}\right)u^{(2)}=0\\
\mu^{(1)}_b+\mu^{(2)}_b-\left(T^{(1)}+T^{(2)}\right)u^{(1)}=0\\
\end{array}\right.
\end{equation}
Local equilibrium reads:
\begin{equation}
\left\{\begin{array}{l}
\left(S_p^{(1)}+T^{(1)}\right)u^{(1)}_b = \mu^{(1)}_b + b_p^{(1)} \\
\left(S_p^{(2)}+T^{(2)}\right)u^{(2)}_b = \mu^{(2)}_b + b_p^{(2)}
\end{array}\right.
\end{equation}
Assuming $T\s$ is chosen so that matrix $\left(S_p\s+T\s\right)$ is invertible, we can express displacement $u\s_b$ from local equilibrium equation, and suppress it from global interface conditions, which leads to:
\begin{multline}\label{eq:mixt_sys:2}
\begin{pmatrix}
I & I - \left(T^{(1)}+T^{(2)}\right)\left(S_p^{(2)}+T^{(2)}\right)^{-1}\\ \left(T^{(1)}+T^{(2)}\right)\left(S_p^{(1)}+T^{(1)}\right)^{-1} & I
\end{pmatrix}
\begin{pmatrix}
\mu^{(1)} \\ \mu^{(2)}
\end{pmatrix} =\\ \begin{pmatrix} \left(T^{(1)}+T^{(2)}\right)\left(S_p^{(2)}+T^{(2)}\right)^{-1} b_p^{(2)} \\ \left(T^{(1)}+T^{(2)}\right)\left(S_p^{(1)}+T^{(1)}\right)^{-1}b_p^{(1)} \end{pmatrix}
\end{multline}
This expression enables to give better interpretation of the
stiffness parameters $T\s$. Suppose $T^{(1)}=S_p^{(2)}$ and
$T^{(2)}=S_p^{(1)}$ then matrix \eqref{eq:mixt_sys:2} is equal to
identity and solution is directly achieved. So the aim of matrix
$T\s$ is to provide one substructure with the interface stiffness
information of the other substructures.

If we generalize to $N$-subdomain system \eqref{eq:mixt_sys:1}, we
can deduce that the optimal choice for $T\s$ is the Schur
complement of the remaining substructures on the interface of
domain $(s)$ (some kind of $S_p^{(\bar{s})}$ where $\bar{s}$
denotes all the substructures but $s$). Of course such a choice is
not computationally feasible (mostly because it does not respect
the localization of data), and approximations have to be
considered. In decreasing numerical efficiency and computational
cost order, we have:
\begin{itemize}
\item Approximate the Schur complement of the remaining of the substructure by the Schur complement of the neighbor;
\item approximate the Schur complement of the neighbor by the Schur complement of the nearer nodes of the neighbor ("strip"-preconditioners which idea is developed in \cite{PAZ:2005:ISP} in another context);
\item approximate the Schur complement of the neighbor by the stiffness matrix of the interface of the neighbor (strategy of dual approach lumped preconditioner).
\end{itemize}
The second strategy is quite a good compromise: it respects data
localization, it is not computationally too expensive and yet it
enables the propagation of the information beyond the interface.
Of course an important parameter is the definition of elements
"near the interface", which can be realized giving an integer $n$
representing the number of layers of elements over the interface.
\subsubsection{Coarse problem}
Because the interface stiffness parameter $\block{T}$ regularizes
local operators $\block{S}_p$, local operator
$\left(\block{S}_p+\block{T}\right)$ is always invertible. Such a
property can be viewed as an advantage because it simplifies the
implementation of the method introducing no kernel and generalized
inverse; but it also can be considered as a disadvantage because
no more coarse problem enables global transmission of data among
the structure. Then the communications inducted by this method are
always neighbor-to-neighbor which means that the transmission of a
localized perturbation to far substructures is always a slow
process. It is then necessary to add an optional coarse problem
(see section \ref{sub:augmkry}). Most often the optional coarse
problem is constituted of would-be rigid body motions (if
subdomains had not been regularized). Another possibility, which
is proposed inside the LaTIn method is to use rigid body motions
and extension modes of each interface as coarse problems, this
leads to much larger coarse space. The coarse matrix corresponds
to the virtual works of first order of deformation of
substructures; so mechanically it realizes and propagates a
numerical first order homogenization of the substructures.
\subsection{Hybrid approach}
The hybrid approach (see \cite{GOSSELET:2004:DDM} for a specific
application) is a proposition to provide a unifying scheme for
primal and dual approaches though it could easily be extended to
other strategies. It relies on the choice for each interface
degree of freedom of its own treatment (for now primal or dual).
So let us define two subsets of interface degrees of freedom: the
first is submitted to primal conditions (subscript $p$) and the
second to Neumann conditions (subscript $d$). Local equilibrium
then reads ($\bp=i\cup d$, $b=d\cup p$, $p \cap d = \emptyset$):
\begin{equation}\label{eq:localhybrid:1}
\begin{pmatrix}
\block{K}_{\bp\bp} & \block{K}_{\bp p} \\  \block{K}_{p \bp} & \block{K}_{p p}
\end{pmatrix}
\begin{pmatrix}
\block{u}_{\bp} \\  \block{u}_{p}
\end{pmatrix} =
\begin{pmatrix}
 \block{f}_{\bp} \\  \block{f}_p
\end{pmatrix}+
\begin{pmatrix}
{\block{\traceop}_d}^T \block{\lambda}_{d} \\  \block{\lambda}_{p}
\end{pmatrix}
\end{equation}
Preferred interface unknowns are unique displacement on the first
subset $\ubold_p$ and unique effort on the second subset
$\du{\lambdabold}_d$. Local contributions then reads:
\begin{eqnarray}\label{eq:localhybrid:2}
\block{u}_p &=& \plock{\assemop}_p^T \ubold_p \\
\block{\lambda}_d &=& \plock{\du{\assemop}}_d^T \du{\lambdabold}_d
\end{eqnarray}
which ensure the continuity of displacement on the $p$ degrees of
freedom and the action-reaction principle on the $d$ degrees of
freedom, of course operators $\plock{\assemop}_p$ and
$\plock{\du{\assemop}}_d$ have been restricted respectfully to the
$p$ and $d$ subsets. Remaining interface conditions read:
\begin{eqnarray}\label{eq:localhybrid:3}
\plock{\du{\assemop}}_d \block{u}_d = \plock{\du{\assemop}}_d \block{\traceop}_d \block{u}_{\bp} &=& 0 \\
\plock{\assemop}_p \block{\lambda}_p &=& 0
\end{eqnarray}
To obtain the global interface system, first local unknown $\block{u}_{\bp}$ has to be
eliminated:
\begin{equation}
\block{u}_{\bp} = {\block{K}_{\bp\bp}}^+ \left( \block{f}_{\bp}
+{\block{\traceop}_d}^T \block{\lambda}_{d} - \block{K}_{\bp p}
\block{u}_{p} \right) + \block{R}_{\bp}\block{\alpha}
\end{equation}
Applying continuity condition to previous result and equilibrium
condition to the second row of \eqref{eq:localhybrid:1}, interface
system reads:
\begin{equation}
\begin{pmatrix}
\Sbold_{p\du{d}} & \begin{pmatrix} \Gbold_p \\ \du{\Gbold}_d \end{pmatrix} \\
\begin{pmatrix} -\Gbold_p^T & \du{\Gbold}_d^T \end{pmatrix} &0 \end{pmatrix}
\begin{pmatrix} \ubold_p \\ \du{\lambdabold}_d \\ \block{\alpha} \end{pmatrix} =
\begin{pmatrix} \bbold_p \\ -\du{\bbold}_d \\ -\block{e} \end{pmatrix}
\end{equation}
with the following notations:
\begin{gather*}
\Sbold_{p\du{d}}=\begin{pmatrix} \plock{\assemop}_p & 0 \\ 0 & \plock{\du{\assemop}}_d
\end{pmatrix}  \block{S}_{pd} \begin{pmatrix} \plock{\assemop}_p & 0 \\ 0 & \plock{\du{\assemop}}_d
\end{pmatrix}^T\\
\Gbold_p= \plock{\assemop}_p \block{K_{p\bp}} \block{R}_{\bp},\
\du{\Gbold}_d=\plock{\du{\assemop}}_d \block{\traceop}_d
\block{R}_{\bp}\\
\block{e}= {\block{R}_{\bp}}^T\block{f}_{\bp} \\
\bbold_p= \plock{\assemop}_p \left(\block{f}_p- \block{K}_{p\bp} {\block{K}_{\bp\bp}}^+
\block{f}_{\bp}\right),\ \du{\bbold}_d=\plock{\du{\assemop}}_d \block{\traceop}_d {\block{K}_{\bp\bp}}^+ \block{f}_{\bp}
\end{gather*}
This interface problem corresponds to the constrained resolution
of one linear system. The constraint is linked to the possible
non-invertibility of matrix $\block{K}_{\bp\bp}$ and thus to the
choice of primal subset. Notice that $\Gbold_p$ represents the
reaction of primal degrees of freedom to zero energy modes of
$\block{K}_{\bp\bp}$ and then should be zero in most cases (it may
be non zero in buckling cases). Moreover this system may represent
classical primal approach (if all interface degrees of freedom are
in subset $p$) or classical dual approach (if all interface
degrees of freedom are in subset $d$). Operator $\Sbold_{p\du{d}}$
is a primal/dual Schur complement, it is the sum of local
contributions $\block{S}_{pd}$ \eqref{eq:loc_hyb_sc:1}.

The above system is nonsymmetric semi-definite (because of
redundancies on the dual subset) positive, it has to be solved by
GMRes-like algorithm.
\subsubsection{Hybrid preconditioner}
Inspired by primal and dual preconditioners, we propose to
approximate the inverse of the sum of local contributions by a
scaled sum of local inverses.
\begin{equation}
\tilde{\Sbold}_{p\du{d}}^{-1}= \begin{pmatrix}\plock{\tilde{\assemop}}_p&0\\ 0 & \plock{\tilde{\du{\assemop}}}_d \end{pmatrix} {\block{S}_{dp}} \begin{pmatrix}\plock{\tilde{\assemop}}_p&0\\ 0 & \plock{\tilde{\du{\assemop}}}_d \end{pmatrix}^T
\end{equation}
Scaled assembly operator are defined in the same way as in primal
and dual approaches.
\subsubsection{Coarse problems}
As said earlier, depending on the choice of subset $p$, local
operator $\block{K}_{\bp\bp}$ (involved in the computation of
$\block{S}_{pd}$) may not be invertible and, like in dual
approach, a first coarse correction has been incorporated inside
the hybrid formulation. Anyhow local operator $\block{K}_{\bd\bd}$
involved in preconditioning step may also not be invertible and,
like in primal approach, a second coarse problem has to be added
to make the preconditioner optimal. Then, the optimal version of
the hybrid system incorporates two coarse problems handled by
specific initialization/projection algorithm. The admissibility
condition reads:
\begin{equation}
\Gbold^T \xbold = e
\end{equation}
with $\Gbold=\begin{pmatrix}\Gbold_p\\ \du{\Gbold}_d\end{pmatrix}$, $\xbold=\begin{pmatrix}\ubold_p \\ \du{\lambdabold}_d\end{pmatrix}$ and $\bbold=\begin{pmatrix}\bbold_p \\ -\du{\bbold}_d\end{pmatrix}$.
If $\rbold=\bbold-\Sbold_{p\du{d}}\xbold$ stands for the residual before preconditioning, the optimality condition reads:
\begin{equation}
\Hbold^T \rbold = \begin{pmatrix}\Hbold_p \\ \du{\Hbold}_d \end{pmatrix}^T \begin{pmatrix} \rbold_p \\ \rbold_d\end{pmatrix}= 0
\end{equation}
with $\Hbold_p=\plock{\tilde{\assemop}}_p{\block{\traceop}_p}{\block{R}_{\bd}}$ and $\du{\Hbold}_d=\plock{\tilde{\du{\assemop}}}_d\block{K}_{d \bd}{\block{R}_{\bd}}$ (as said before most often $\du{\Hbold}_d=0$).
To sum up constraints:
\begin{equation}
\left\{\begin{array}{l} \Gbold^T \xbold= -\block{e} \\ \Hbold^T \Sbold_{p\du{d}} \xbold = \Hbold^T \bbold \end{array} \right.
\end{equation}
Handling such constraints is described in section \ref{sub:consaugmkry}.

\begin{figure}[ht]\centering
\includegraphics[width=0.8\textwidth]{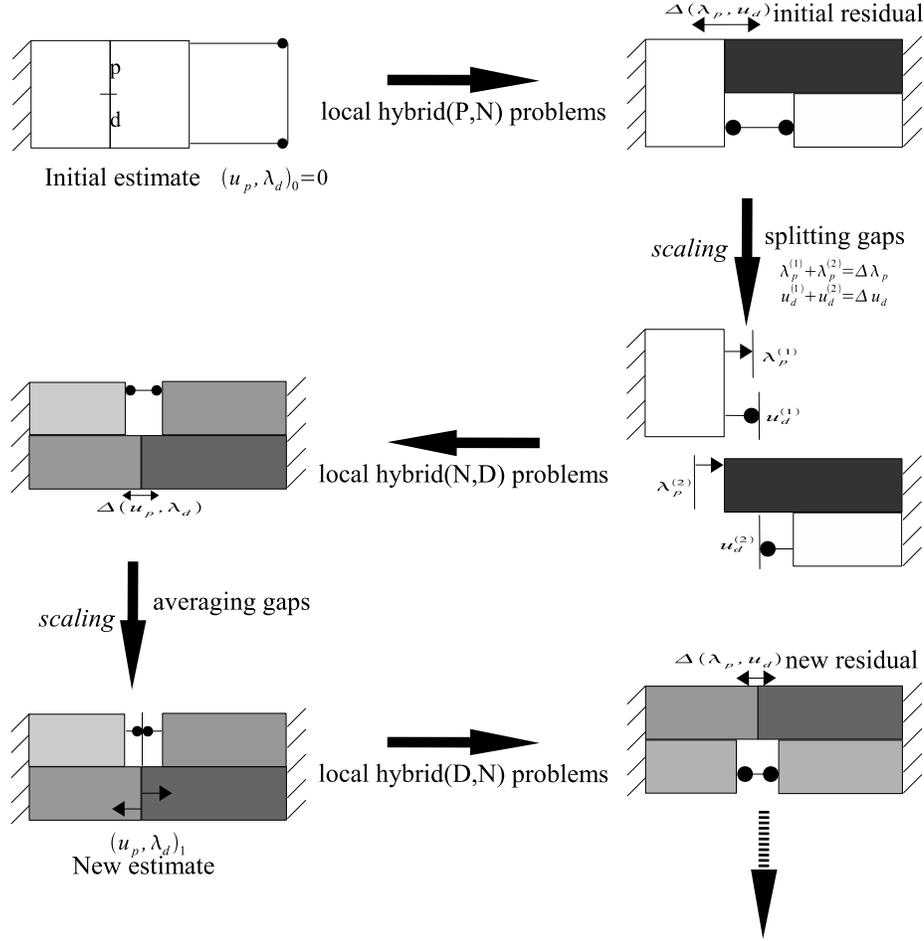}\caption{Representation of first iteration of preconditioned hybrid approach}\label{fig:sch_hybrid}

\end{figure}

Figure \ref{fig:sch_hybrid} provides a schematic representation of
the first iteration of the preconditioned hybrid approach, in the
specific case of a nodal partition of the interface. Assessments
will deal with partition at the degree of freedom level.

\subsubsection{Error estimate}
Because GMRes-like solver is used, the euclidian norm of the
residual is directly available, such a norm is the sum of
displacement gap on the $d$ part of the interface and an effort
gap on the $p$ part of the interface. For now, no other estimator
with better physical sense is available.

\section{Adding optional constraints}\label{sec:constraints}
The aim of optional constraints is to ensure that the research
space for the iterative solver possesses a certain regularity. The
choice of constraints will be discussed in a later section. Going
back to the interface system:
\begin{equation}
\begin{pmatrix}
\block{S}_p & -\block{I} \\ \plock{\du{\assemop}} & 0 \\ 0 & \plock{\assemop}  \end{pmatrix} \begin{pmatrix} \block{u}_b \\ \block{\lambda}_b\end{pmatrix}=\begin{pmatrix} \block{b}_p  \\ \du{0} \\ 0\end{pmatrix}
\end{equation}
any constraint of the form
$\du{\Cbold}^T\plock{\du{\assemop}}\block{u}_b=0$ or
$\Cbold^T\plock{\assemop}\block{\lambda}_b=0$ is trivially
verified by the solution fields, it is just a restriction of
continuity/equilibrium conditions. From an iterative process point
of view, these conditions will be reached once converged ; the
principle of optional constraints is to have every iteration
verify that condition.

There are two classical solutions to ensure these optional
constraints: either to realize a splitting of research space and
ensure, using a projector, that the resolution is limited to
convenient subspace, or to realize a condensation of constraints
and make iterations on a smaller space. In other words suppose
there are $n_c$ independent constraints in a $n$-dimension space
the first strategy researches $n$-sized solution in a
$(n-n_c)$-ranked space, while the second solution researches
$(n-n_c)$-sized solution in a $(n-n_c)$-dimension space then
deduce the $n$-sized solution. From a numerical point of view both
solutions are equivalent, they are just two ways of handling the
same constraints, anyhow from implementation and computational
points of view they have strong differences.

We will essentially focus on application to primal and dual domain decomposition methods.

\subsection{Augmentation strategy}
For this strategy the constraint is reinterpreted in terms of
constraint on the residual. Typically, the primal approach can be
augmented by constraints on the effort field while the dual
approach can have constraints on the displacement field:
\begin{eqnarray}
\Cbold^T\plock{\assemop}\block{\lambda}_b &=&- \Cbold^T\left(\plock{\assemop} \block{b}_p - \plock{\assemop}  \block{S}_p \plock{\assemop} ^T \ubold_b \right) \\
\du{\Cbold}^T\plock{\du{\assemop}}\block{u}_b&=&\du{\Cbold}^T\left(\plock{\du{\assemop}}\block{S}_d{\plock{\du{\assemop}}}^T\du{\lambdabold}_b+\block{\du{b}}_d \right)
\end{eqnarray}
The constraint is then handled as an augmentation inside the
iterative solver, its is classically realized using a projector
(see sections \ref{sub:augmkry} and \ref{sub:consaugmkry}).

\subsection{Recondensation strategy}
This strategy was recently introduced in the framework of the dual
approach, leading to the FETIDP algorithm
\cite{FARHAT:2001:FETI_DP,FARHAT:2000:FETI_DP,Lesoinne:1999:FES,KLAWONN:2005:SCR}.
Because for now only constraints on the $\block{u}_b$ field have
been considered we will restrain to this kind of constraints, the
application of constraints on $\block{\lambda}_b$ is
straightforward.

\subsubsection{Basic method}
We first consider constraints which impose continuity of specific
degrees of freedom, in other terms we suppose matrix $\du{C}$ is
identity on certain degrees of freedom and zero elsewhere; we will
show how any constraint can be rewritten in such a form. Because
these degrees of freedom will be submitted to a primal treatment
we will denote them with subscript $p$ while the remaining of the
interface will be denoted with subscript $d$. Constraint reads:
\begin{equation}
0 = \du{\Cbold}^T\plock{\du{\assemop}}\block{u}_b=\begin{pmatrix} 0 \\ I_p \end{pmatrix}^T \begin{pmatrix}\plock{\du{\assemop}}_d\block{u}_d \\ \plock{\du{\assemop}}_p\block{u}_p \end{pmatrix} = \plock{\du{\assemop}}_p\block{u}_p
\end{equation}

Like in the hybrid approach this constraint is verified using a
unique displacement field on the primal part of the interface:
$\block{u}_p=\plock{{\assemop}}_p^T \ubold_p$. Interface system
then reads like in the hybrid approach, with the additional
assumption that the constraints are so that the local problem
possesses enough Dirichlet conditions to make it invertible.
\begin{equation}
\Sbold_{p\du{d}}
\begin{pmatrix} \ubold_p \\ \du{\lambdabold}_d \end{pmatrix} =
\begin{pmatrix} \bbold_p \\ -\du{\bbold}_d \end{pmatrix}
\end{equation}
Introducing following notations for blocks composing $\block{S}_{pd}$:
\begin{equation}
\block{S}_{pd}=\begin{pmatrix}
\block{s}_{pp} & \block{s}_{pd} \\ -\block{s}_{dp} & \block{s}_{dd}
\end{pmatrix}\end{equation} then \begin{equation}  \Sbold_{p\du{d}}= \begin{pmatrix}
\sbold_{pp} & \sbold_{p\du{d}} \\ -\sbold_{\du{d}p} & \sbold_{\du{d}\du{d}}
\end{pmatrix}= \begin{pmatrix} \plock{\assemop}_p & 0 \\ 0 & \plock{\du{\assemop}}_d
\end{pmatrix} \begin{pmatrix}
\block{s}_{pp} & \block{s}_{pd} \\ -\block{s}_{dp} & \block{s}_{dd}
\end{pmatrix} \begin{pmatrix} \plock{\assemop}_p & 0 \\ 0 & \plock{\du{\assemop}}_d
\end{pmatrix}^T
\end{equation}
Unknown $\ubold_p$ is condensed on the remaining of the interface:
\begin{eqnarray}
\ubold_p&=&\sbold_{pp}^{-1}\left(\bbold_p-\sbold_{p\du{d}}\lambdabold_d \right)\\
\du{\sbold}_d\du{\lambdabold}_d=\left(\sbold_{\du{d}\du{d}}+\sbold_{\du{d}p}\sbold_{pp}^{-1}\sbold_{p\du{d}}\right)\du{\lambdabold}_d&=&-\du{\bbold}_d+\sbold_{\du{d}p}\sbold_{pp}^{-1}\bbold_p
\end{eqnarray}
The latest equation is solved using an iterative solver, operator
$\du{\sbold}_d$ has the same properties as the restriction of the
dual operator to the $d$-part of the interface (semi-definition,
symmetry and positivity). Operator $\du{\sbold}_d$ is the assembly
of local contributions  \begin{equation}
\du{\sbold}_d=\plock{\du{\assemop}}_d\block{s}_d\plock{\du{\assemop}}_d^T=\plock{\du{\assemop}}_d\left(\block{s}_{dd}+\block{s}_{dp}\plock{\assemop}_p^T\sbold_{pp}^{-1}\plock{\assemop}_p\block{s}_{pd}\right)\plock{\du{\assemop}}_d^T
\end{equation}
Using operator $\du{\sbold}_d$ requires the computation of the
inverse of matrix
$\sbold_{pp}=\plock{\assemop}_p\block{s}_{pp}\plock{\assemop}_p^T$,
which is an assembled matrix. Then this formulation includes a
global coarse problem set on primal variables.

The recommended preconditioner for such an approach is directly
inspired by the dual approach: it consists in solving local
Dirichlet problems with scaled imposed displacement on the
$d$-part of the interface and null displacement on the $p$ part of
the interface and extracting the average reaction of the $d$-part
of the interface. Then the preconditioner reads:
\begin{equation}
\du{\tilde{\sbold}}^{-1}_d = \begin{pmatrix}0_p & \plock{\du{\tilde{\assemop}}_d}\end{pmatrix}\block{S}_p\begin{pmatrix}0_p & \plock{\du{\tilde{\assemop}}_d}\end{pmatrix}^T
\end{equation}

\begin{figure}[ht]\centering
\includegraphics[width=0.8\textwidth]{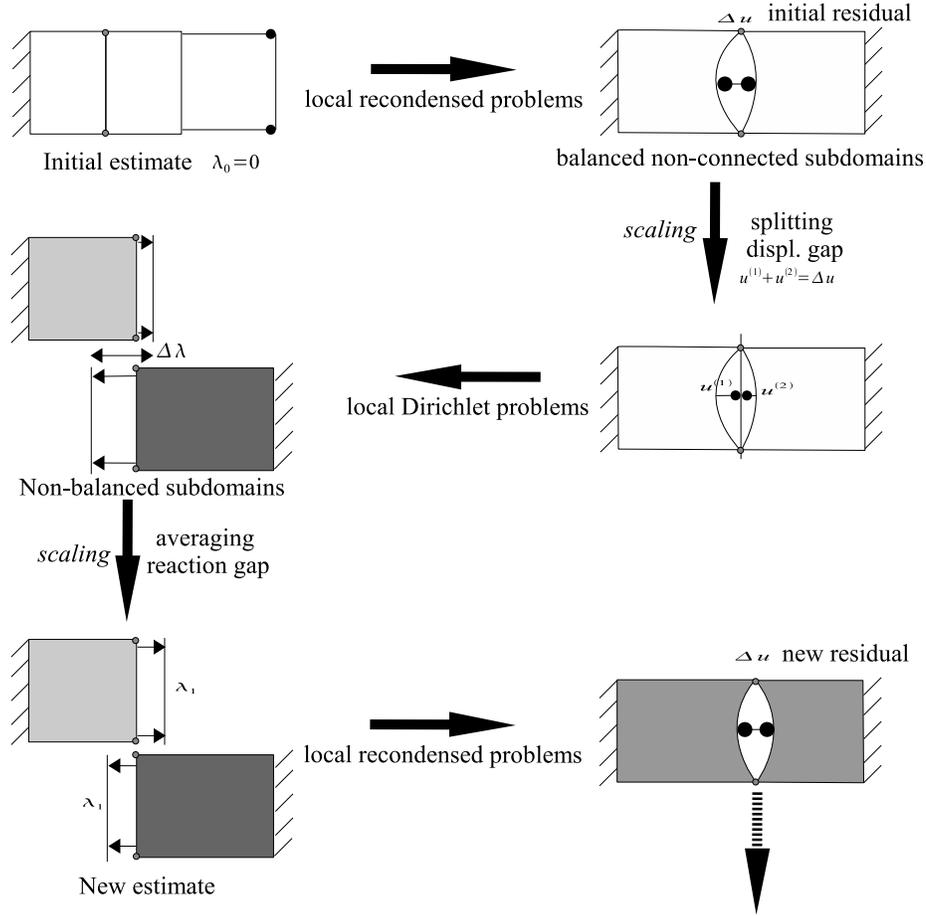}\caption{Representation of the first iteration of preconditioned FETIDP}\label{fig:sch_fetidp}
\end{figure}

Figure \ref{fig:sch_fetidp} provides schematic representation of the first iteration of preconditioned FETID method.

\subsubsection{More complex constraints}
We now consider constraints which are not limited to one degree of
freedom, for instance one can consider that we want to ensure that
the average jump of displacement on one edge is equal to zero,
which involves all the degrees of freedom of the edge.

The classical solution \cite{KLAWONN:2005:SCR} is to realize a
change of basis of degrees of freedom (denoted by matrix
$\block{T}$) so that each constraint is represented by one
"modified" degree of freedom. The change of basis is the same
local operation realized on every subdomain, then we can define a
global change $\Tbold$ so that
$\plock{\assemop}\block{T}=\Tbold\plock{\assemop}$
\begin{eqnarray}\label{eq:change_basis}
\Cbold^T \plock{\assemop} \block{u}_b = \Cbold^T \plock{\assemop}\block{T}\block{\hat{u}}_b =\Cbold^T\Tbold \plock{\assemop}\block{\hat{u}}_b &=& \hat{\Cbold}^T \plock{\assemop}\block{\hat{u}}_b \\ \text{with\ } \hat{\Cbold}&=&\begin{pmatrix}0\\I_{\hat{p}}\end{pmatrix}
\end{eqnarray}
After the change of basis is realized the same algorithm can be
employed. Because constraints most often respect a certain
locality of data (for instance independent constraints on each
edge), change of basis is not a too expensive operation, and does
not make too poor the sparsity of local matrices.

\subsection{Adding "constraints" to the preconditioner}
This subsection deals with the dualization of the recondensation
strategies \cite{CROS:2002:PSC,DOHRMANN:2003:PSB}. For instance
the balanced domain decomposition with constraints (BDDC) is a
primal version of FETIDP: during the preconditioning step, the
continuity of displacement is ensured at specific degrees of
freedom (which can be the result of a local change of basis), so
that the local Neumann operator remains fully invertible. So
solving classical primal approach problem
\begin{equation}
\Sbold_p \ubold_b = \plock{\assemop}\block{S}_p\plock{\assemop}^T \ubold_b = \bbold_p
\end{equation}
the preconditioner reads
\begin{equation}
\tilde{\Sbold}_p^{-1} = \begin{pmatrix}I_p & 0\\0& \plock{\displaystyle{\tilde{\assemop}}}_d\end{pmatrix} \begin{pmatrix} \sbold_{pp}^{-1} & -\sbold_{pp}^{-1}\plock{\assemop}_p\block{s}_{pd} \\ -\block{s}_{dp}\plock{\assemop}_p^T\sbold_{pp}^{-1}& \block{s}_d \end{pmatrix} \begin{pmatrix}I_p & 0 \\ 0 &  \plock{\displaystyle{\tilde{\assemop}}}_d\end{pmatrix}^T
\end{equation}

Figure \ref{fig:sch_fetidp} provides a schematic representation of
the first iteration of preconditioned BDDC method.

\begin{figure}[ht]\centering
\includegraphics[width=0.8\textwidth]{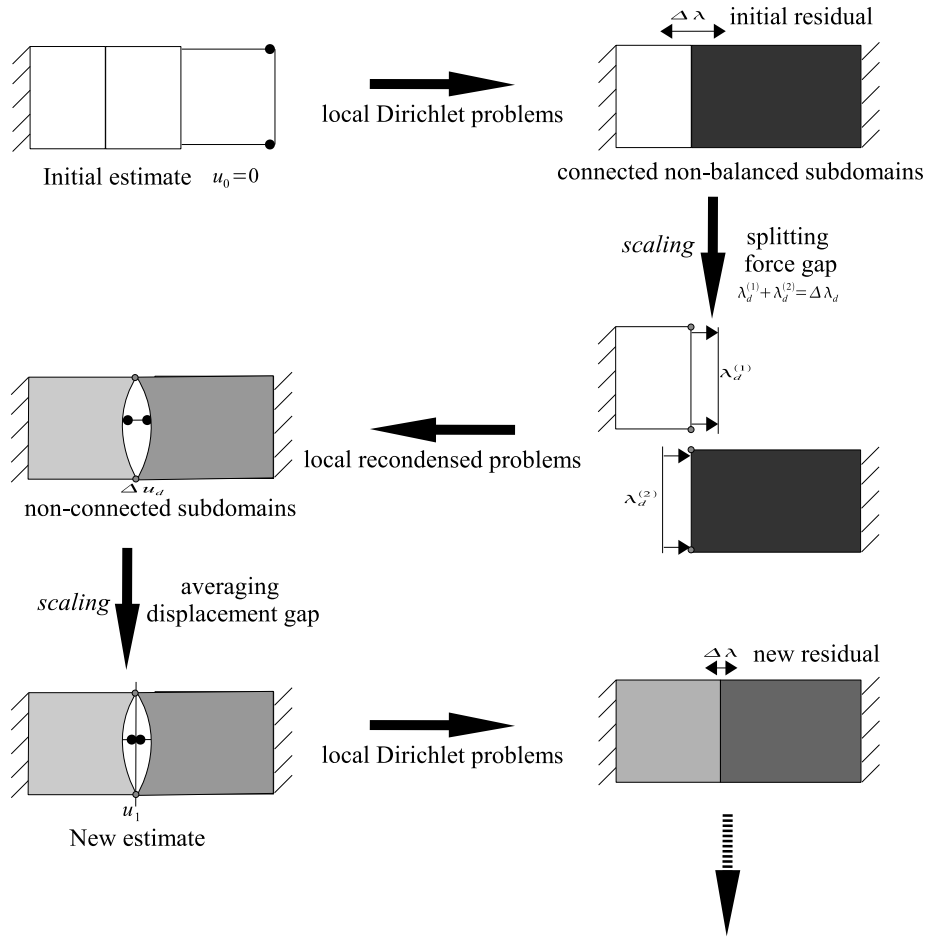}\caption{Representation of the first iteration of preconditioned BDDC}\label{fig:sch_bddc}
\end{figure}

\section{Classical issues}\label{sec:issu}
In this section we try to provide answers to questions that
commonly arise when using domain decomposition methods.
\subsection{Rigid body motion detection}
Handling floating substructures is definitely a very accurate
issue. This difficulty is one of the reason of the success of
methods which regularize the stiffness of subdomains like FETIDP
or mixed approaches, leading to fully invertible matrices. Anyhow
basic primal and dual methods remain competitive (mostly because
zero-energy modes provide a very natural coarse problem), hence
providing an efficient algorithm for the computation of rigid body
motion is essential. Many strategies can be used
\cite{FARHAT:1996:CNS,FARHAT:2000:FETINL}, and this review does
not claim to be exhaustive.

First we have to discuss the exact composition of the possible kernel of substructures. What can be found is:
\begin{itemize}
\item rigid body motions of floating substructures,
\item internal mechanisms of substructures,
\item weird things due to nonlinearities (buckling, exotic behaviors, ...)
\item numerical zero-energy modes.
\end{itemize}
An internal mechanism can for instance occur when a substructure
is made out of two parts connected by one linear edge (pivot) or
one singular point (kneecap). Methods exist to avoid such
substructures either inside the decomposition algorithm or as
external programs regularizing a given decomposition
\cite{FRAGAKIS:2002:UFF}, and then should be employed.

The last two kinds of kernel are non-standard and can only be
detected using fully algebraic methods (like Gauss pivoting). The
problem with algebraic methods is their high sensitivity to the
condition number of the stiffness matrix of substructures. The
condition number can be influenced by the aspect ratio and the
material composition of the substructures, even after
adimensionalization the quality of the methods is hard to
warranty.

Finally we only develop here two strategies to handle zero-energy
modes. The first one belongs to the purely algebraic methods, it
is very simple to implement and can lead to satisfying results for
not-too-complex problems, it can handle more than solid body
motions but it is strongly dependent on \textit{a priori} selected
degrees of freedom. The second method is purely geometric, it is
very robust but only suited to detect solid body motions.

In order to simplify notations, we consider the research of the
zero-energy modes of matrix $K$ (which should be a local stiffness
matrix).

\subsubsection{Simple algebraic approach}
This method is based on fundamental relationship between the
kernel of a matrix and the kernel of Schur complement
\eqref{eq:kernels}. The principle is to preselect a small set of
degrees of freedom which we will denote by subscript $N$ (the
other degrees of freedom are denoted with subscript $O$). Then
compute explicitly primal Schur complement associated to these
degrees of freedom: $S=K_{NN}-K_{NO}K_{OO}^{-1}K_{ON}$). If
$N$-degrees of freedom are selected so that $K_{00}$ is invertible
(if only solid body motions have to be detected then it is
sufficient to take the degrees of freedom associated to three
non-aligned nodes) then $S$ is well defined and its kernel is
linked to the kernel of matrix $K$.

Since $N$ is a "small" set (12 degrees of freedom is often
sufficient), computing the kernel of matrix $S$ using "exact"
algorithm like singular value decomposition is rather cheap and
then kernel of matrix $K$ can be computed using equation
\eqref{eq:kernels}.

\subsubsection{Geometric approach}
The basic idea of this method is that kinematically admissible
solid body motions can be deduced from boundary conditions imposed
on one subdomain. Let $R_c$ be a basis of candidate rigid body
motions (would be solid body motions if no Dirichlet conditions
were applied on the subdomain, $R_c$ is a 6-column matrix in 3D
and 3-column matrix in 2D), and let $E$ be the matrix of Dirichlet
boundary conditions: each column of $E$ represents one
(combination of) blocked degree of freedom. Kinematically
admissible rigid body motions $R$ are linear combinations of
candidate rigid body motions (hence $R=R_c Q$) which do not make
Dirichlet boundary conditions work (\textit{ie} $E^TR=0$). In
order to find such linear combination, we  compute singular value
decomposition of matrix $E^T R_c = U D V^T$ and set $Q = V_0$
where $V_0$ is the submatrix of $V$ associated to negligible
singular values. Because $E^TR_c$ is a matrix only made out of
geometric considerations, it is well conditioned and the criterion
to detect "zero" singular values is well defined (there is a large
gap between zero and non-zero singular values).

\subsubsection{Generalized inverse}
The two methods presented above led to $r$-ranked basis $R$ of the
kernel of matrix $K$. To compute generalized inverse $K^+$, the
most classical way is to select $r$ degrees of freedom, so that if
they were added Dirichlet conditions, rigidity matrix would be
well defined. For instance after pivoting and renumber one could
get $R=\begin{pmatrix} R^* \\ I_r\end{pmatrix}$ and then the $r$
last degrees of freedom would suit. Then we have:
\begin{equation}
K=\begin{pmatrix} K_{\bar{r}\bar{r}} & K_{\bar{r}r} \\ K_{r\bar{r}}& K_{rr} \end{pmatrix} \qquad \text{and} \qquad K^+=\begin{pmatrix} K_{\bar{r}\bar{r}}^{-1} & 0 \\ 0 & 0 \end{pmatrix}
\end{equation}
This is just one instance of generalized inverse, other can be
built using penalization or other modification to matrix $K$.
Though not theoretically prohibited, choosing "blocked" degrees of
freedom on the interface is often a bad idea from a practical
point of view.

\subsection{Choice of optional constraints}
As seen in section \ref{sec:constraints}, constraints can be
imposed either using augmentation (using one or two projectors,
sections \ref{sub:augmkry} and \ref{sub:consaugmkry}) or using
recondensation algorithms. In the case of recondensation
algorithms, constraints have to be sufficient in order to suppress
rigid body motions and then regularize the local stiffness matrix.
In the case of augmentation algorithm, constraints have to be
independent from rigid body motions which are already handled by
the formulation.

Because constraints are expensive to handle, they have to be
chosen with care; anyhow, except in a few cases, there are no
general results on how to choose them. Mechanical comprehension of
studied phenomena and anticipation of convergence difficulties may
lead to efficient strategies. In the case of solving several
linear systems (even with different left hand sides) interesting
strategies exist
\cite{SAAD:1987:LMS,ERHEL:2000:ACG,GOSSELET:2002:SRK,REY:2003:SLN,RISLER:1998:RRV}.

The next two subsections deal with very common strategies, while
the last subsection describes another framework for constraints
inspired by the LaTIn method \cite{LAD99b,NOUY:2003:THES}.

\subsubsection{Forth order elasticity}
As plate and shell models are often used in structural mechanics,
forth order problems have been carefully studied
\cite{FARHAT:1996:COR1,FARHAT:1996:COR2,ROUX:1997:DIRECT2,LETALLEC:1997:SHELL,FARHAT:2000:FETINL}.
Such problems are characterized by the appearance of singularities
at the corner of substructures (so-called "corner modes") which
are destroying the scalability of usual methods. The classical
solution consists in enforcing the continuity of the (most often
only normal) displacement field at the corners in order to
regularize the problem. From a practical point of view, corners
are most often defined as multiple points (nodes shared by more
than two substructures), that set can be enriched by extremities
of edges.

The implementation, just like the singularity, is strongly linked to the chosen formulation.
\begin{description}
\item[Dual approach: ]
since the projected residual corresponds to the displacement jump
between substructures \eqref{eq:dualresid}, one just has to use
augmentation algorithm with one constraint for each pair of
neighbor at each corner. Matrix $\du{\Cbold}$ is then made out of
columns with one coefficient $1$ on one corner degree of freedom
and $0$ elsewhere. Because the dual description of a $m$-multiple
point leads to $(m-1)$ relationships, such a coarse space is
rather large.
\item[Primal approach: ]
In order to regularize the displacement field, the constraints
have to be imposed on the preconditioned residual (assuming
Neumann-Neumann preconditioner is employed). The aim is then to
have the local contributions of preconditioned residual equal to
zero on corner points. So if $\block{C}$  denotes the local
interface matrix made out of columns with one coefficient $1$ on
one corner degree of freedom and $0$ elsewhere, and
$\block{\tilde{C}}$ the same matrix scaled according to the
scaling used inside the preconditioner, the constraints read $C =
\plock{\tilde{\assemop}} \block{S}_d \block{\tilde{C}}$. Then a
$m$-multiple degree of freedom leads to $m$ constraints.

\item[Recondensed approaches: ]
FETIDP or BDDC were first designed in this context, from the
consideration that (extended) corners constraints were sufficient
to suppress rigid body motions, then the first level constraints
could be avoided. So the methods directly apply since they consist
in constraining the displacement field. Here whatever the
multiplicity of a corner may be, it always leads to just $1$
constraint.
\end{description}

\subsubsection{Second order elasticity}
Because classical methods are already scalable in the frame of
second order elasticity, optional constraints are not often used
in such a context. Furthermore, it is hard to predict what
constraint should be imposed. In some cases, efficient strategies
have been proposed, such as in \cite{GOSSELET:2002:SRK} for
nonlinear problems using Newton-Raphson solver where
approximations of eigen vectors are used.

The question of optional constraints arose when willing to extend
recondensed algorithms (FETIDP and BDDC) to such problems, mostly
because the previous definition of corners lead to significant
problems in 3D (too many constraints, poor convergence...). The
first solution was proposed in \cite{Lesoinne:1999:FES}, the idea
was to select 3 non-aligned nodes on each face (interface between
2 subdomains) which maximized the surface of the triangle they
defined. The current solution, the scalability of which is
mathematically and numerically proved, is to enforce average
convergence on edges \cite{KLAWONN:2005:SCR}, which is realized by
a change of basis described in \eqref{eq:change_basis}. In order
to take into account heterogeneity on the interface, the average
may be scaled by a coefficient representing the stiffness of the
subdomains. For more difficult problems, first order moments of
edges can also be added.

\subsubsection{Link with homogenization theory}
This paragraph intends to present a mechanical vision of optional
constraint which, though hard to implement in the framework of the
presented method, may lead to better understanding of what
optional constraints and associated coarse problems can provide to
the methods. This analysis is inspired by the multiscale version
of the LaTIn method \cite{LAD99b}.

The underlying question when choosing optional constraints (except
from specific numerical questions like in the forth order
elasticity) is "what global information should be transmitted to
the whole structure ?" or more precisely "what should far
substructures know from one substructure". A meaningful answer is
provided by Saint-Venant principle and homogenization theory: at a
first order development, a substructure can be represented by its
rigid body motions and its constant strain states (simple traction
and shearing states). Such an idea adds six (3 in 2D) more
constraints per subdomains; as these constraints are somehow
complex to build, they can be approximated by interface modes
(but of course the number of constraints then grows quickly).

\subsection{Linear Multiple Points Constrains}
Multiple points constraints are relationships defined between some
degrees of freedom, they are often used in order to connect
nonconforming meshes, to represent boundary conditions (for
instance periodicity), to model contact or apply control laws. In
the case of linear(ized) constraints we can write, on the whole
structure scale:
\begin{eqnarray}
Ku&=&f\\
C u &=& a
\end{eqnarray}
What seems most suited to the domain decomposition context is to
dualize the constraint and introduce Lagrange multiplier $\mu$ in
order to enforce the condition. System then reads:
\begin{equation}
\begin{pmatrix}
K & C^T \\ C & 0
\end{pmatrix}\begin{pmatrix}
u \\ \mu
\end{pmatrix}=\begin{pmatrix}
f \\ a
\end{pmatrix}
\end{equation}
After decomposition we have
\begin{equation}
\begin{pmatrix}
\block{K} & -\block{I} & {\plock{C}}^T \\ \plock{\du{\assemop}}\block{\traceop} & 0 & 0 \\ 0 & \plock{\assemop}\block{\traceop} & 0\\ {\plock{C}} & 0 & 0
\end{pmatrix}\begin{pmatrix}
\block{u} \\ \block{\lambda} \\ \mu
\end{pmatrix}=\begin{pmatrix}
\block{f} \\ \du{0} \\ 0 \\ a
\end{pmatrix}
\end{equation}
with $\plock{C}$ so that $\plock{C}\block{u}=Cu=a$ which implies (since $\block{u}_b=\plock{\assemop}^Tu_b$):
\begin{eqnarray}
Cu=\begin{pmatrix}
{\plock{C}}_i & C_b
   \end{pmatrix} \begin{pmatrix}
\block{u}_i \\ u_b
   \end{pmatrix} &=& \plock{C}\block{u} = \begin{pmatrix}
{\plock{C}}_i & {\plock{C}}_b
   \end{pmatrix} \begin{pmatrix}
\block{u}_i \\ \block{u}_b
   \end{pmatrix} \\
   \text{then}\qquad{\plock{C}}_b\plock{\assemop}^T&=&C_b
\end{eqnarray}
Or in other words, if matrix $C$ deals with interface degrees of
freedom, the associated constraints have to be correctly
distributed between sharing subdomains. The constraint can be
interpreted as specific (non-boolean) assembly operator which
explains the chosen notation. Using MPCs with domain decomposition
methods was studied in \cite{RIXEN:2002:EPF} in the dual method
context.

In order to provide general methodology to apply MPCs, we then
incorporate inside hybrid domain decomposition method (equations
\eqref{eq:localhybrid:1} to \eqref{eq:localhybrid:3}):
\begin{eqnarray}
\begin{pmatrix}
\block{K}_{\bp\bp} & \block{K}_{\bp p} \plock{\assemop}_p^T & {\plock{C}}_{\bp}^T \\
\plock{\assemop}_p\block{K}_{p \bp} & \plock{\assemop}_p\block{K}_{p p}\plock{\assemop}_p^T & \Cbold_p^T \\
{\plock{C}}_{\bp} & \Cbold_p & 0
\end{pmatrix}
\begin{pmatrix}
\block{u}_{\bp} \\  \ubold_p \\ \mu
\end{pmatrix} &=&
\begin{pmatrix}
 \block{f}_{\bp}+{\block{\traceop}_d}^T\plock{\du{\assemop}}_d^T \lambda_{d} \\ \plock{\assemop}_p \block{f}_p \\ a
\end{pmatrix}\\ \plock{\du{\assemop}}_d{\block{\traceop}_d}\block{u}_{\bp}&=&\du{0}
\end{eqnarray}
The elimination of $\block{u}_{\bp}$ leads to, with classical hybrid notations:
\begin{equation}\label{eq:hybridmpc}
\begin{pmatrix}
\Sbold_{p\du{d}} &   \begin{array}{c} \mathcal{C}_p^T \\ \mathcal{C}_d^T  \end{array}& \begin{array}{c} \Gbold_p \\ \du{\Gbold}_d \end{array}\\
\begin{array}{cc} -\mathcal{C}_p & \mathcal{C}_d  \end{array} &\mathcal{C}_\mu & \Gbold_\alpha \\\begin{array}{cc} -\Gbold_p^T & \du{\Gbold}_d^T \end{array} &  \Gbold_\alpha^T & 0 \end{pmatrix}
\begin{pmatrix} \ubold_p \\ \du{\lambdabold}_d \\ \mu \\ \block{\alpha} \end{pmatrix} =
\begin{pmatrix} \bbold_p \\ -\du{\bbold}_d  \\ h \\ -\block{e} \end{pmatrix}
\end{equation}
with
\begin{eqnarray}
\mathcal{C}_p & = & \Cbold_p-{\plock{C}}_{\bp}{\block{K}_{\bp \bp}}^{+}\block{K}_{\bp p}\plock{\assemop}_p^T \\
\mathcal{C}_d & = &  {\plock{C}}_{\bp}{\block{K}_{\bp \bp}}^{+}{\block{\traceop}}^T_d\plock{\du{\assemop}}_d^T \\
\mathcal{C}_\mu & = & {\plock{C}}_{\bp} {\block{K}_{\bp \bp}}^{+}{\plock{C}}_{\bp}^T \\
\Gbold_\alpha & = & {\plock{C}}_{\bp}R_{\bp} \\
h&=&a-{\plock{C}}_{\bp}{\block{K}_{\bp \bp}}^{+}\block{f}_{\bp}
\end{eqnarray}

Various strategies can be used in order to solve system
\eqref{eq:hybridmpc}, which combine elimination of constraints
(rigid body motions and/or MPCs) by projection methods (FETI-like
approaches) and/or by recondensation methods (FETIDP like
approaches). All these methods correspond to solving rather
complex coarse problems but they are suited to traditional
preconditioners. We propose, after \cite{RIXEN:2002:EPF}, to use
classical projection method to handle rigid body motions then use
iterative solver to find simultaneously $(\ubold_p,
\du{\lambdabold}_d,\mu)$ and provide efficient preconditioner to
this problem.

System reads with trivial notations (for simplicity reasons we
suppose that the rigid body motion constraints have been
symmetrized, which is always possible and which is naturally the
case if $G_p=0$ like in many applications):
\begin{equation}
\begin{pmatrix}
\Sbold_{p\du{d}\mu} & \mathcal{G} \\ \mathcal{G}^T & 0
\end{pmatrix}  \begin{pmatrix} \xbold \\ \block{\alpha} \end{pmatrix} =
\begin{pmatrix} \bbold \\ -\block{e} \end{pmatrix}
\end{equation}
with
\begin{equation}
\Sbold_{p\du{d}\mu} = \begin{pmatrix} \plock{\assemop}_p & 0 \\ 0 & \plock{\du{\assemop}}_d\block{\traceop}_d \\ 0 & {\plock{C}}_{\bp} \end{pmatrix}
\block{S}_{pd}
\begin{pmatrix} \plock{\assemop}_p & 0 \\ 0 & \plock{\du{\assemop}}_d\block{\traceop}_d \\ 0 & {\plock{C}}_{\bp} \end{pmatrix}^T + \begin{pmatrix} 0 & 0 & \Cbold_p \\ 0& 0 &0\\-\Cbold^T_p &0 &0\end{pmatrix}
\end{equation}
As can  be seen, MPCs have very different actions wether they are set on primal interface degrees of freedom or not: matrix $C_p$ modifies the structure of the hybrid system while dual and internal constraints lead to classical hybrid approach with modified dual trace and assembly operator $\plock{\mathcal{A}}=\begin{pmatrix} \plock{\du{\assemop}}_d\block{\traceop}_d \\ {\plock{C}}_{\bp}\end{pmatrix}$. The definition of efficient preconditioner inspired by classical methods is then much simplified if no constraints are set on primal degrees of freedom ($C_p=0$), which is what we suppose now:
\begin{equation}
\Sbold_{p\du{d}\mu} = \begin{pmatrix} \plock{\assemop}_p & 0 \\ 0 & \plock{\mathcal{A}} \end{pmatrix}
\block{\mathcal{S}}_{pd}
\begin{pmatrix} \plock{\assemop}_p & 0 \\ 0 &\plock{\mathcal{A}} \end{pmatrix}^T
\end{equation}
The $\block{\mathcal{S}}_{pd}$ notation is due to the association
of the trace operator with the assembly operator, which in fact is
equivalent to defining "extended interface dual degrees of
freedom" made out of dual degrees of freedom and degrees of
freedom involved in MPCs. Since, in this hypothesis, system reads
like classical hybrid approach, that is (modified) assembly of
local contributions, the proposed preconditioner is a scaled
assembly of local inverses.
\begin{equation}
\tilde{\Sbold}_{p\du{d}\mu}^{-1} = \begin{pmatrix} \plock{\tilde{\assemop}}_p & 0 \\ 0 & \plock{\tilde{\mathcal{A}}} \end{pmatrix}
{\block{\mathcal{S}}_{dp}}
\begin{pmatrix} \plock{\tilde{\assemop}}_p & 0 \\ 0 &\plock{\tilde{\mathcal{A}}} \end{pmatrix}^T
\end{equation}
The primal scaled assembly operator can be directly imported from the primal approach. Concerning the dual approach, according to previous definitions, we have:
\begin{equation}
\plock{\tilde{\mathcal{\assemop}}} = \left(\plock{\mathcal{\assemop}} {\block{M}_{\bp}}^{-1}\plock{\mathcal{\assemop}}^T\right)^{+}\block{\mathcal{\assemop}}{\block{M}_{\bp}}^{-1}
\end{equation}
Where $\block{M}_{\bp}$ is a diagonal matrix chosen like in the classical methods. The matrix to inverse reads:
\begin{equation}
\left(\plock{\mathcal{\assemop}} {\block{M}_{\bp}}^{-1}\plock{\mathcal{\assemop}}^T\right) = \begin{pmatrix}
\plock{\du{\assemop}}_d\block{\traceop}_d{\block{M}_{\bp}}^{-1}{\block{\traceop}_d}^T\plock{\du{\assemop}}_d^T & \plock{\du{\assemop}}_d\block{\traceop}_d{\block{M}_{\bp}}^{-1}{\plock{C}}_{\bp}^T \\
{\plock{C}}_{\bp} {\block{M}_{\bp}}^{-1} {\block{\traceop}_d}^T\plock{\du{\assemop}}_d^T &
{\plock{C}}_{\bp} {\block{M}_{\bp}}^{-1} {\plock{C}}_{\bp}^T
\end{pmatrix}
\end{equation}
The idea is then to make this system easy to inverse, having the off-diagonal blocks equal to zero. We have ${\plock{C}}_{\bp}=\begin{pmatrix}{\plock{C}}_i & {\plock{C}}_d \end{pmatrix}$ with ${\plock{C}}_d{\plock{\du{\assemop}}}_d^T=C_d$; if we choose
\begin{equation}
{\plock{C}}_d=C_d\left({\plock{\du{\assemop}}}_d{\block{\traceop}_d}\block{M}_{\bp}{\block{\traceop}_d}^T{\plock{\du{\assemop}}}_d^T\right)^{-1}{\plock{\du{\assemop}}}_d{\block{\traceop}_d}\block{M}_{\bp}{\block{\traceop}_d}^T
\end{equation}
then $\plock{\du{\assemop}}_d\block{\traceop}_d{\block{M}_{\bp}}^{-1}{\plock{C}}_{\bp}^T =\du{0}$ and the scaled assembly operator is non expensive to compute.

In others words, one simply has to split constraints on interface
degrees of freedom between subdomains according to the scaling
used inside the preconditioner.

\subsection{Choice of decomposition}
Decomposing a given structure in order to use the algorithms
presented in this paper is a complex problem. Algorithms and
softwares were proposed
\cite{FARHAT:1993:TOP,WALSHAW:1995:PAO,KARYPIS:1998:MAM}, which
mostly refer to graph theory. Such approaches enable to take into
account load balance between processors (supposing each processor
is assigned to one subdomain, this is equivalent to making each
local problem as easy to solve as others) and to minimize the
dimension of the interface (so that the global condensed problem
is as small as possible). They also enable to avoid internal rigid
body motions (so called mechanisms).

Anyhow another important point is to have local problems as well
conditioned as possible, so having subdomains with good aspect
ratio (ratio between largest and smallest characteristical
dimensions of the subdomain) is considered as an important point.
Anyhow one has to realize that good aspect ratio is often linked
to large local bandwidth and then to somehow expensive local
problems to solve.

\begin{figure}[ht]\centering\psfrag{U}{$\Upsilon$}
\subfigure[Interface avoiding heterogeneity]{\includegraphics[width=0.2\textwidth]{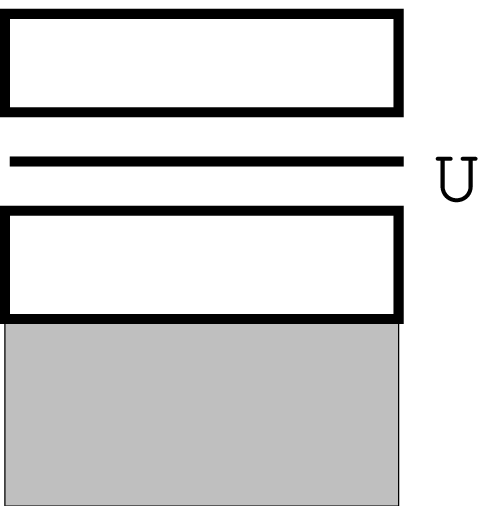}}\qquad
\subfigure[Interface matching heterogeneity]{\includegraphics[width=0.2\textwidth]{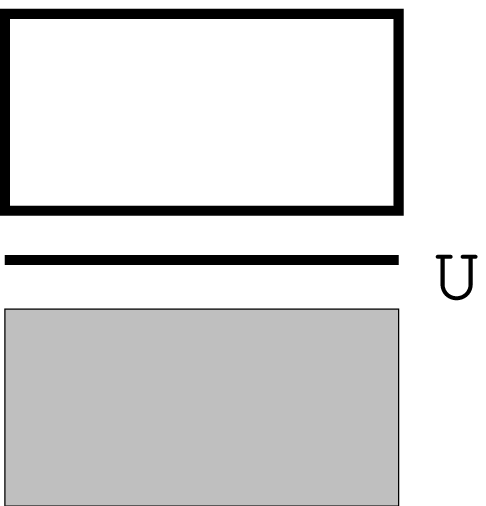}}\qquad
\subfigure[Interface crossing heterogeneity]{\includegraphics[width=0.2\textwidth]{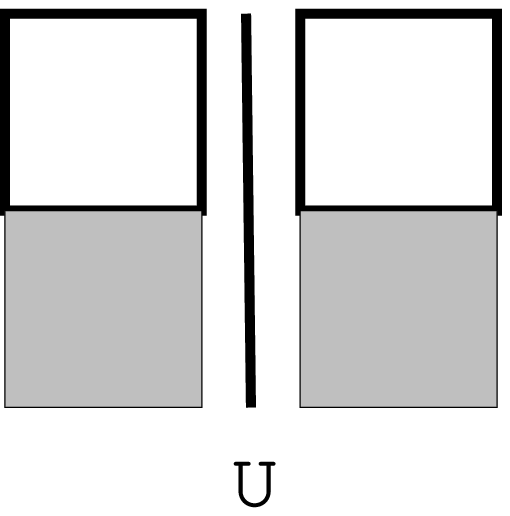}}\caption{Different kinds of heterogeneity in domain decomposition context}\label{fig:hete}
\end{figure}

An even more difficult point to take into account high
heterogeneities (see figure \ref{fig:hete}): using stiffness
scaling enables to correctly handle heterogeneity when interface
between subdomains matches interface between materials, anyhow
when interface between subdomains "crosses" interface between
materials then numerical difficulties may occur. For now
scaled-average optional constraints seem to be the best solution
to handle these difficulties but it leads to large coarse problems.

Finally finding the best decomposition is still a rather open
problem and mechanical sense is often a necessary complement to
efficient automatic decomposing algorithm.

\subsection{Extensions}
\subsubsection{Nonsymmetric problems}
Nonsymmetry occurs in many physical modeling: plasticity, nonlocal
models for fracture \cite{GERMAIN:2005:MCNL}, frictional contact
\cite{BARBOTEU:2001:DDS}. The use of the domain decomposition
methods presented in this paper just requires more care in the
implementation because some simplifications are not available (for
instance coarse problem matrices are nonsymmetric), and of course
the use of well suited iterative solver like GMRes, OrthoDir or
BiCG because Schur complements are no longer symmetric.

Globally methods show good numerical performance results. Anyhow a
real problem is the absence of theoretical results to ensure good
convergence properties (this is mainly due to the fact that proofs
for classical methods rely on the construction of an interface
inner-product related to Schur complements which is no longer
possible).

\subsubsection{Nonlinear problems}
For now we have considered the solution to linear systems. To adapt the
method to nonlinear problems, a classical solution is to use
Newton-Raphson linearization scheme: linearized stiffness matrices
are computed independently on each subdomain then the linearized
system is solved using domain decomposition
\cite{REY:1998:RKS,GOSSELET:2002:SRK,RISLER:1998:RRV,FARHAT:2000:FETINL,GOSSELET:2002:DDM}.
For such approaches, domain decomposition methods can be seen as
efficient black-boxed linear solvers.

One critical point in such a method is that, depending on the
formulation (for instance fully Lagragian formulation of a large
deformation elasticity problem) rigid body motions may vary from
one system to the other. In the proposed context, what has been
proved is that translations always belong to zero energy modes,
what has been observed is that rotations only appeared as
zero-energy-modes in the first system (which corresponds to linear
elasticity problem). This might be penalizing because the size of
the information transmitted inside the coarse problem is decreased
after the first system; moreover, rotations are often converted to
"negative-energy modes" which, if they are in small number, can be
handled by fully-reorthogonalized conjugate gradient (though
convergence will be slower). One classical solution is to reinject
disappeared rotations as optional constraints (via augmentation
algorithms).

In the case where nonlinearity is localized in few substructures,
an interesting strategy can be to carry out subiterations of the
nonlinear solver independently in those substructures
\cite{CAI:2002:NPI,CRESTA:2005:CMP}.

\subsection{Implementation issues}\label{sec:implementation}
Implementing domain decomposition methods from existing code is
not a too complex task. We give a few details on our
software architecture though practical solutions are far from
being unique. Our code is a plug-in to ZeBuLoN object-oriented
finite element computational software
\cite{ZEBUUSER:2001,ZEBUDEVEL:2001}, it takes advantage of
Frederic Feyel's previous work
\cite{FEYEL:1998:THESE,FEYEL:2005:HDR}. Our implementation aims at
being as generic as possible, so for now hybrid domain
decomposition method has been developed (and mixed approaches are
under construction), and separation between formulation and solver
(so that any iterative solver can be used to solve the interface
problem). All classical projectors and preconditioners are
implemented.

The most basic pieces of the code are:
\begin{itemize}
\item topological description of the interface, \textit{ie} ability to realize trace operations;
\item exchange library (PVM, MPI), \textit{ie} ability to realize assembly operations;
\item classical FE code, \textit{ie} ability on any given subdomain to get any local field given sufficient boundary conditions.
\end{itemize}

\subsubsection{Organization of the topological information}
In order to implement hybrid approach, we propose to define a
specific class for "interface degrees of freedom" which wraps
classical degrees of freedom and provide information on:
\begin{itemize}
\item the number of subdomains which share this degree of freedom
\item the kind of treatment (primal, dual...) which is applied
\end{itemize}

Then degrees of freedom are set together neighborwise, defining
"interface" class made out of a list of pointer to "interface
degrees of freedom", and the global number of the neighbor (that
number enables to identify the subdomain and to realize
exchanges).

A "subdomain" is defined as a collection of "interfaces" and a
classical domain in the sense of usual FE software (mainly mesh),
it possesses its own global identification number.

Note that with such a description of the decomposition, the local
interfaces are redundant, multiple degrees of freedom appear in
several interfaces. It is then necessary to take certain care to
define some operations (transposed trace for dual degrees of
freedom). A specific class can then be used to ease the management
of multiple degrees of freedom.

\begin{figure}\centering
\includegraphics[width=6cm]{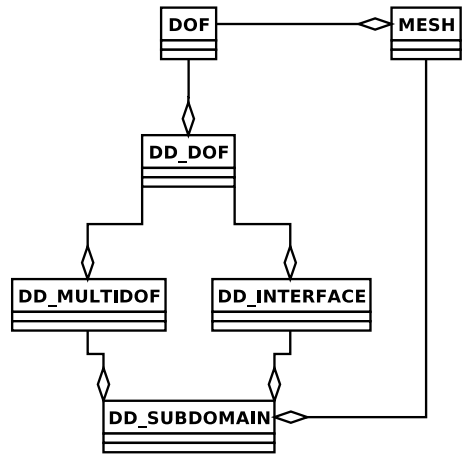}
\caption{Topological interface information}\label{fig:dd_topo}
\end{figure}

\begin{figure}\centering
\includegraphics[width=9cm]{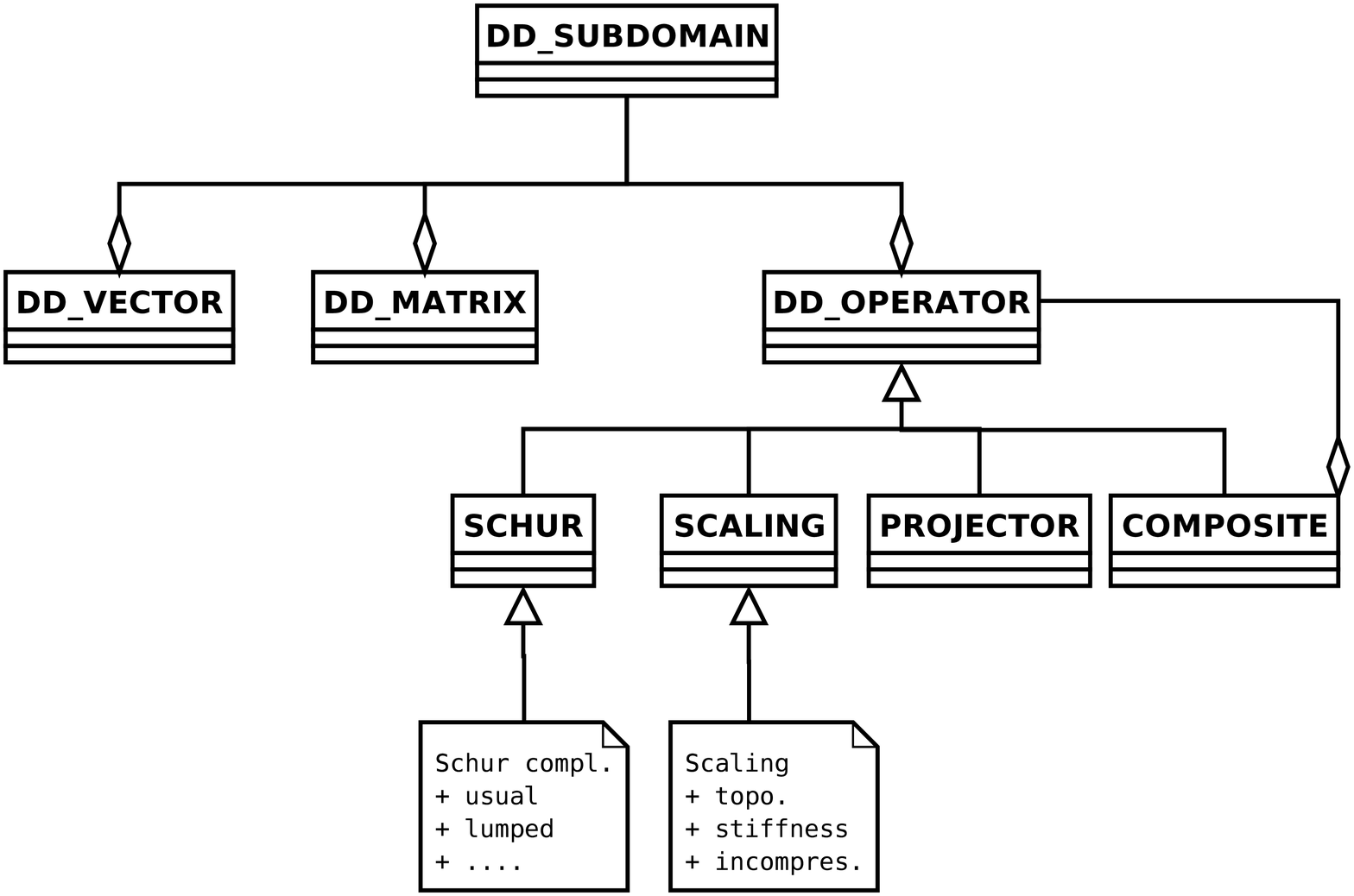}
\caption{Algebraic interface objects}\label{fig:dd_alge}
\end{figure}

\subsubsection{Defining algebraic interface objects (fig. \ref{fig:dd_alge})}
In order to easily connect the domain decomposition formulation to
an iterative solver, we propose to define "interface vectors"
(displacement, intereffort...), "interface matrices" (trace of
rigid body motions..., can be seen as a collection of interface
vectors) and "interface operators" (square interface matrices),
with all classical operations (basically sum, difference, product
and transposed product).

The particularity of these objects is to be defined on the
interface and then data is shared between subdomains, so all the
previous operations sometimes require to assemble data (an
interesting idea is to have a boolean member indicating the
assembled state of data). Because of the choice of description of
the interface, the assembly operation requires certain care for
primal multiple degrees of freedom (these degrees of freedom shall
have the same value at their different occurrences). Note that an
object like "interface matrix" can highly be optimized (mainly in
terms of memory storage).

Interface operators mainly know how to multiple vectors and
matrices, they are used to define Schur complements, scaling
operators, projectors. Composite design pattern can be used to
simplify the succession of operations.

In order to let the user choose the various configurations of the
domain decomposition method, we use inheritance and object factory
design pattern.

\begin{figure}\centering
\includegraphics[width=7cm]{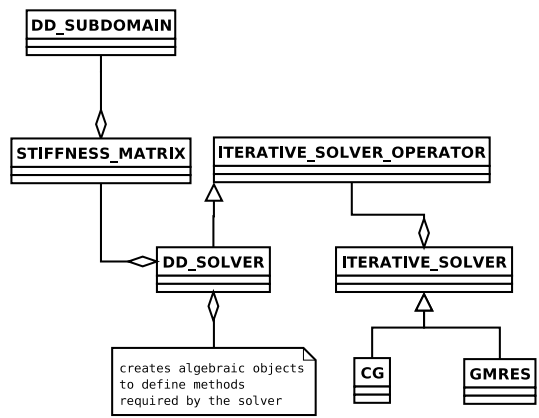}
\caption{Articulation between formulation and solver}\label{fig:dd_solver}
\end{figure}

\subsubsection{Articulation between formulation and solver (fig. \ref{fig:dd_solver})}
What we propose is to have client/server relationship between
solver and formulation: basically an iterative solver needs to
know how to initialize, how to multiply, how to precondition, how
to do inner products, how to evaluate convergence. All these
operations are implemented inside the "interface formulation"
object which is linked to one subdomain (topology and stiffness)
and creates "interface algebraic objects" in order to define
required operations.

\section{Assessments}\label{sec:asse}
The assessments we propose here are based on the code described in
the previous section. Basically, we have implemented the hybrid
approach which lets us assess the classical primal and dual
approaches with most classical preconditioners and projectors.

We first present a sequence of academical tests in order to
recover classical numerical performance results (scalability and
relative efficiency of the different approaches): two dimensional
plane stress problem, plate problem, three dimensional problem
with heterogeneity and unstructured decomposition. In these
problems, $H$ denotes the characteristic size of the subdomains
and $h$ the characteristic size of the elements. We also present
results on non-academic problem (bitraction test specimen).

For all those tests, in order to compare all the approaches
(including the hybrid approach), GMRes solver is used and
convergence is monitored by the norm of the residual as given by
the solver (with $\varepsilon$ set to $10^{-6}$). In other cases
(when hybrid approach is not assessed), the convergence is
monitored through the evaluation of global primal residual
$|Ku-f|/|f|<\varepsilon$ with $\varepsilon$ set to $10^{-6}$.
Depending on the method, a different coarse problem may be
introduced, we denote by CS:a+b the size of the coarse problems (a
for the admissibility coarse problem and b for the optimal
preconditioning coarse problem) or
$\textit{number\_of\_iterations}_{\text{total\_number\_of\_constraints}}$.
Note that the hybrid approach deals with two independent coarse
problems, so their solutions is much cheaper than the solution to
a unique large coarse problem.

We test the primal approach with or without optimality coarse
problem, the dual approach with lumped or Dirichlet preconditioner
and identity or superlumped or Dirichlet projector (denoted by
P(I), P(W) and P(D) respectively). As for such tests no physical
consideration can guide the choice of hybrid treatment to the
interface, in order to show the potential of the hybrid approach,
we present results where all degrees of freedom of one direction
($U_1$, $U_2$ or $U_3$) are treated in the same way. For instance
"D-P" stands for a dual treatment for degrees of freedom
associated to direction $U_1$, and a primal treatment for degrees
of freedom associated to direction $U_2$.

\subsection{Two dimensional plane stress problem}
\begin{figure}
    \centering
    \includegraphics[width=5cm ]{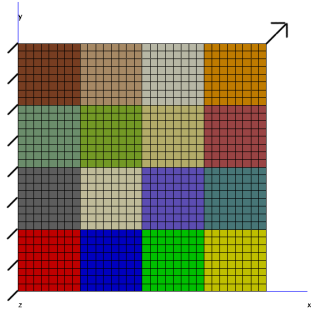}
    \caption{16 subdomains decomposed square}
    \label{fig:carre}
\end{figure}

We consider a simple second order two-dimensional problem, the
structure is an homogeneous square decomposed in square
substructures meshed with linear square finite element (Q1
Lagrange). The behavior is linear elastic (Young modulus
$E=200000$ MPa and Poisson coefficient $\nu=0.3$), the loading
consists in clamping on the left side and punctual effort on the
top right corner (figure \ref{fig:carre}).

\begin{table}[ht]\centering
\begin{tabular}[b]{|l|l|c|c|c|c|}
\hline
\multicolumn{2}{|c|}{\backslashbox{\qquad Method\qquad }{\qquad H/h\qquad}} & 8 & 16 & 32 & 64 \\
\hline \multirow{2}{2cm}{Primal}
& SC:0+0 & 44 & 45 & 45 & 45 \\
& SC:0+36 & 11 & 12 & 14 & 15 \\\hline
\multirow{3}{2cm}{Dual \\
SC:36+0}
& Lumped - P(I) & 14 & 25 & 32 & 42  \\
& Dirichlet - P(I) & 13 & 15 & 17 & 20 \\
& Dirichlet - P(D) & 12 & 14 & 15 & 17 \\\hline
\multirow{2}{2cm}{Hybrid D-P}
& SC:12+0 - P(D)& 29 & 30 & 33 & 35 \\
& SC:12+12 - P(D)& 12 & 14 & 16 & 18 \\\hline
\multirow{2}{2cm}{Hybrid P-D}
& SC:12+0 - P(I) & 29 & 32 & 36 & 38 \\
& SC:12+12 - P(I) & 14 & 17 & 20 & 22 \\
& SC:12+0 - P(D) & 26 & 29 & 31 & 33 \\
& SC:12+12 - P(D) & 12 & 14 & 16 & 18 \\\hline
\end{tabular}\caption{Scalability results in 2D / 16 subdomains}\label{tab:sd16}
\end{table}

\begin{table}[ht]\centering
\begin{tabular}[b]{|l|l|c|c|c|c|c|c|c|}
\hline
\multicolumn{2}{|c|}{\backslashbox{Method}{nb. subdomains}} & 4 & 9 & 16 & 25 & 36 & 49 & 64 \\
\hline \multirow{2}{2cm}{Primal\\(Neumann$^2$)}
& No opt. coarse & $13_{0}$ & $29_{0}$ & $45_{0}$ & $63_{0}$ & $83_{0}$ & $102_{0}$ & $126_{0}$  \\
& With opt. coarse & $8_{6}$ & $10_{18}$ & $12_{36}$ & $13_{60}$ &
$14_{90}$ & $14_{126}$ & $15_{168}$
\\\hline
\multirow{3}{2cm}{Dual}
& Lumped - P(I) & $18_{6}$ & $24_{18}$ & $26_{36}$ & $27_{60}$ & $29_{90}$ & $29_{126}$ & $31_{168}$ \\
& Dirichlet - P(I) & $9_{6}$ & $13_{18}$& $15_{36}$ & $16_{60}$ & $17_{90}$ & $18_{126}$ & $19_{168}$ \\
& Dirichlet - P(D) & $9_{6}$ & $12_{18}$ & $14_{36}$ & $15_{60}$ &
$16_{90}$ & $17_{126}$ & $18_{168}$ \\\hline
 \multirow{2}{2cm}{Hybrid D-P\\ P(D)}
& No opt. coarse & $9_{2}$ & $21_{6}$ &$ 30_{12}$ & $40_{20}$  & $50_{30}$ & $60_{42}$ &  $67_{56}$ \\
& With opt. coarse   & $7_{4}$ & $12_{12}$ & $14_{24}$ & $16_{40}$
& $17_{60}$ & $18_{84}$ & $18_{112}$
\\\hline \multirow{2}{2cm}{Hybrid P-D\\ P(D)}
& No opt. coarse & $9_{2}$ & $20_{6}$ & $29_{12}$ & $38_{20}$ & $48_{30}$ & $57_{42}$ &  $57_{56}$ \\
& With opt. coarse  & $7_{4}$ & $12_{12}$ & $14_{24}$ & $16_{40}$
& $17_{60}$ & $18_{84}$ & $17_{112}$\\\hline
\end{tabular}\caption{Performance results in 2D for given $\frac{H}{h}=16$}\label{tab:hH16}
\end{table}

Table \ref{tab:sd16} shows the number of iterations of available
strategies for different number of elements per subdomain (for a
16-subdomain decomposition), and table \ref{tab:hH16} for
different number of subdomains (for given ratio $H/h=16$).
Globally all approaches (primal, dual and hybrid) equipped with
their best preconditioner and projector behave similarly and are
scalable. Note that even in its optimal configuration the hybrid
approach requires a smaller coarse space (for instance, table
\ref{tab:sd16}, two $12\times 12$ coarse problems against one
$36\times 36$ coarse problem for primal or dual approaches) for
equivalent efficiency. As expected if the optimality coarse
problem is suppressed, performance results decay and scalability
is lost. Finally for such simple problems, the simplified versions
of the dual approach gives excellent results.

\subsection{Bending plate}
We consider a forth order plate problem, the structure is an
homogeneous square decomposed in squared substructures meshed with
square Mindlin plate element. The behavior is linear elastic
(Young modulus $E=200000$ MPa and Poisson coefficient $\nu=0.3$),
the loading consists in clamping on the left side and punctual
normal effort on the top right corner.

\begin{table}[ht]\centering
\begin{tabular}[b]{|l|l|c|c|c|c|c|c|}
\hline \multicolumn{2}{|c|}{\backslashbox{Method}{nb. subdomains}}
& 4 & 9 & 16 & 25 & 36 & 49  \\\hline
\multirow{2}{3cm}{Primal \\
(Neumann$^2$)}
& No corner & $15_{12}$ & $24_{36}$ & $32_{72}$ & $40_{120}$ & $47_{180}$ & $55_{252}$\\
 & With corners & $13_{16}$ & $17_{64}$ & $20_{108}$ & $23_{184}$ & $24_{280}$ & $26_{396}$  \\\hline
 \multirow{4}{3cm}{Dual \\ (Dirichlet)}
& P(I) - No corner & $17_{12}$ & $31_{36}$ & $43_{72}$ & $59_{120}$ & $75_{180}$ & $91_{252}$     \\
& P(I) - With corners & $16_{15}$& $24_{48}$ & $27_{99}$ & $29_{168}$ & $31_{255}$ & $33_{360}$     \\
& P(D) - No corner & $15_{12}$ & $25_{36}$ & $34_{72}$ & $43_{120}$ & $51_{180}$ & $59_{252}$    \\
& P(D) - With corners & $14_{15}$ & $21_{48}$ & $28_{99}$ &
$31_{168}$ & $31_{255}$ & $36_{360}$ \\\hline

\end{tabular}\caption{Bending plate: performance results for given $\frac{H}{h}=8$}\label{tab:plate}
\end{table}

Table \ref{tab:plate} presents the number of iterations for the
dual and primal approaches, with or without optional corner
constraints (the subscript indicates the total size of coarse
problems, \textit{ie} rigid body motions and corner modes). As
predicted, corner constraints are essential in order to make the
algorithms scalable. Anyhow the dimension of the coarse space
associated to corners quickly explodes which makes the methods
less interesting from a CPU time point of view, which justifies
the FETIDP philosophy which leads to much smaller coarse spaces.

\subsection{Heterogeneous 3D problem}
\begin{figure}[ht]\centering
\subfigure[27 subdomains structured decomposition of a
heterogeneous
cube]{\includegraphics[width=0.3\textwidth]{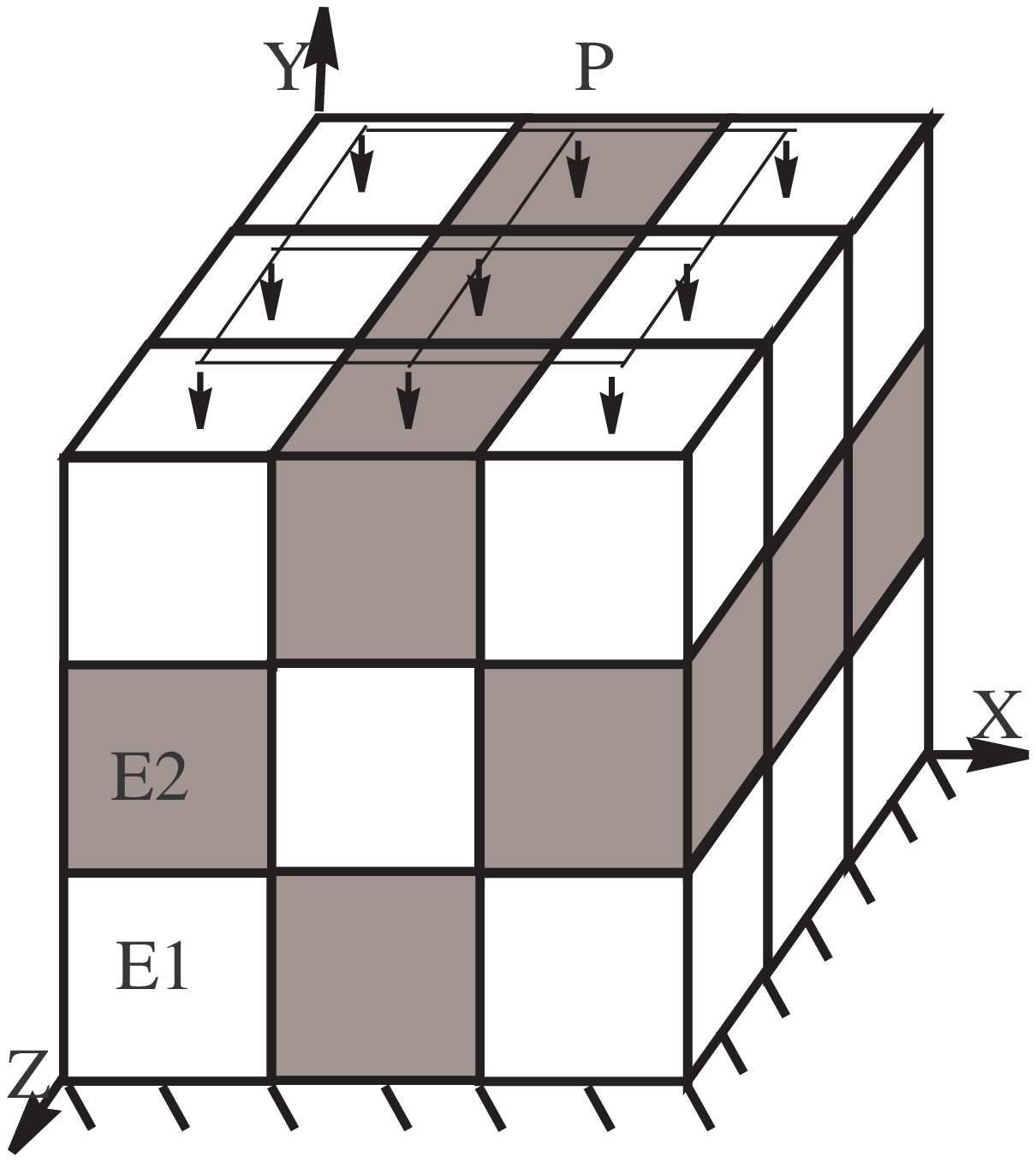}}\label{fig:cubhet}
\qquad \subfigure[27 subdomains unstructured decomposition of a
cube]{\includegraphics[width=0.4\textwidth]{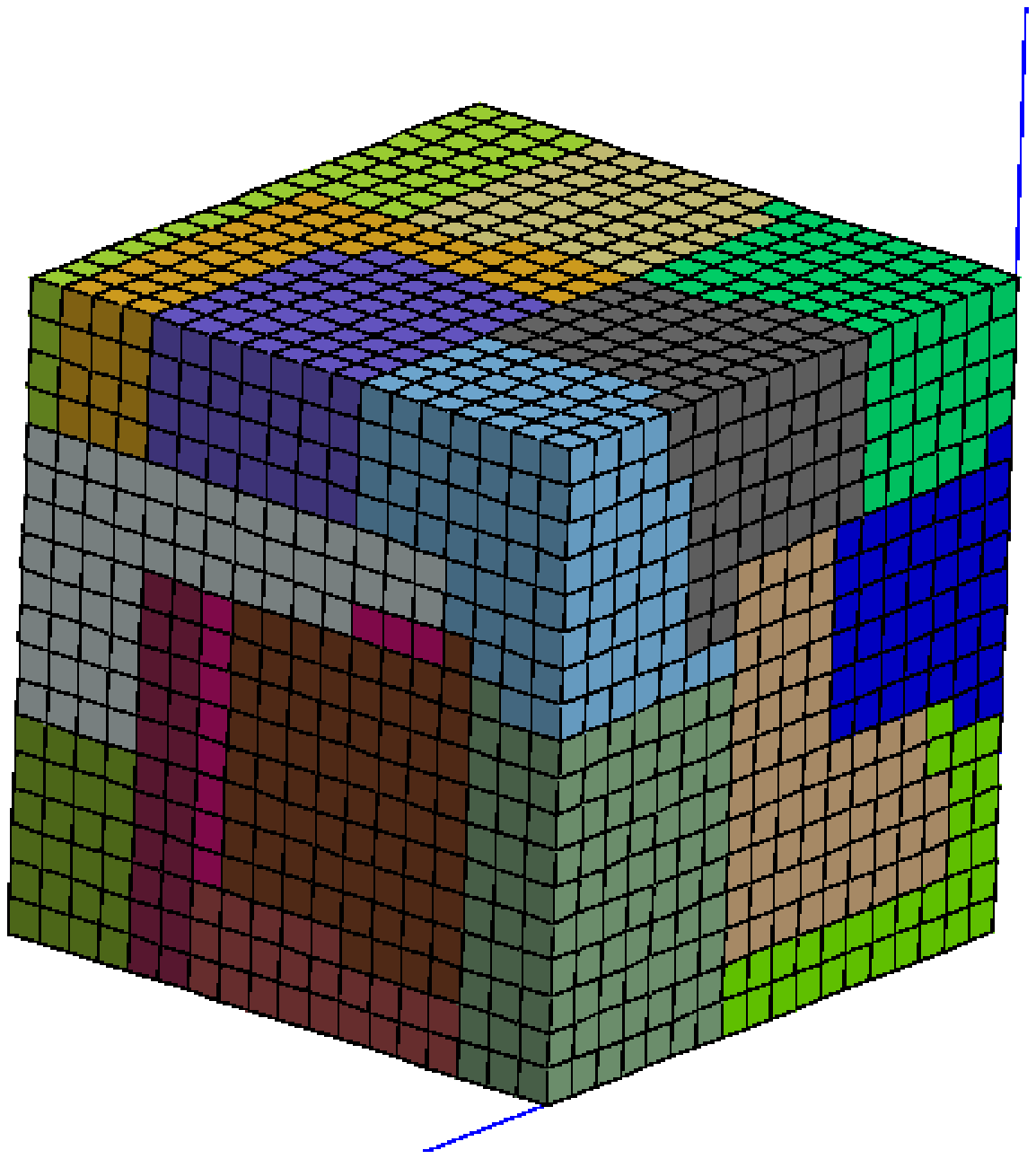}}\label{fig:cubunst}
\caption{3D assessments}\label{fig:cubes}
\end{figure}

We consider a 3D problem, the structure is an heterogeneous cube
decomposed in $3\times3\times3$ cubic substructures meshed with
$3\times3\times3$ Q2-Lagrange cubic elements (27 nodes per
element). The heterogeneity pattern is described in figure
\ref{fig:cubes}a, behaviors are linear elastic (Young modulus
$E_1=200000$ MPa, $E_2=2$ MPa and Poisson coefficient $\nu=0.3$),
the loading consists in clamping on the bottom side and constant
pressure on top side.

\begin{table}[ht]\centering
\begin{tabular}[b]{|l|l|c|}
\hline
\multicolumn{2}{|c|}{Method} & Number of iterations  \\
\hline
\multicolumn{2}{|l|}{Primal}  & 19\\\hline
\multirow{3}{4cm}{Dual  P(D)}
& No splitting & 28    \\
& Classical splitting & 28    \\
& Condensed splitting & 18    \\\hline
\multirow{3}{3cm}{Dual P(W)}
& No splitting & 21    \\
& Classical splitting & 21    \\
& Condensed splitting & 20    \\\hline
\multirow{3}{4cm}{Dual P(I)}
& No splitting & 74    \\
& Classical splitting & 74    \\
& Condensed splitting & 73    \\\hline
\end{tabular}\caption{Heterogeneous cube}\label{tab:cubhet}
\end{table}

Table \ref{tab:cubhet} presents the number of iterations for the
conjugate gradient to converge. Assessed methods are classical
primal approach and dual approach with different projectors for
all splittings (or equivalent initializations) presented in
section \ref{subsub:fetiini}, of course stiffness scaling is
employed. What appears clearly is the good behavior of the
approaches face to heterogeneity except the identity projector of
the dual approach (which is definitely not suited to heterogeneous
problems), and the efficiency of the condensed initialization. For
such a problem the superlumped projector leads to very good
results, anyhow for more complex cases Dirichlet projector is
necessary and shall be improved at no extra computational cost by
the condensed initialization.

\subsection{Homogeneous non-structured 3D problem}
We consider a 3D problem, the structure is an homogeneous cube
meshed with Q1-Lagrange cubic elements (8 nodes per element). The
behavior is linear elastic (Young modulus $E=200000$ MPa and
Poisson coefficient $\nu=0.3$), the loading consists in clamping
on the bottom side and constant pressure on top side. We consider
two kinds of decomposition: either structured decompositions
($3\times3\times3$ or $4\times4\times4$ cubic substructures) or so
called "unstructured" decompositions realized by Metis software
(\url{http://www-users.cs.umn.edu/~karypis/metis/}), see figure
\ref{fig:cubes}b.

\begin{table}[ht]\centering
\begin{tabular}[b]{|l|c|c|c|c|}
\hline
\backslashbox{Method}{\multirow{2}{3.5 cm}{\hfill Decomposition}} & \multicolumn{2}{|c|}{Structured} &  \multicolumn{2}{|c|}{Unstructured} \\ & 27 sd. & 64 sd. & 27 sd. & 64 sd. \\\hline
Primal Neumann-Neumann & 11$_{108}$  & 14$_{288}$ & 42 & 67  \\\hline %
Dual Dirichlet P(I) & 12$_{108}$ & 16$_{288}$ & 44 & 69  \\
Dual Dirichlet P(D) & 12$_{108}$ & 16$_{288}$ & 43 & 70  \\\hline %
Hybrid D-D-P P(I)& 14$_{72}$  & 17$_{192}$ & - & -  \\
Hybrid P-P-D P(I)& 15$_{72}$  & 19$_{192}$ & - & -  \\
Hybrid D-P-D P(I)& 13$_{72}$  & 17$_{192}$ & - & -  \\\hline
\end{tabular}\caption{Homogeneous cube / influence of the decomposition}\label{tab:cubhom}
\end{table}

Table \ref{tab:cubhom} enables to highlight the fundamental role
played by the decomposition: scalabily result only holds for
structured decomposition; moreover there might by a huge
performance gap between two decompositions with the same number of
subdomains (factor 3 for 27 subdomains, factor 4 for 64
subdomains).

\subsection{Bitraction test specimen}
In order to assess "real life" problems, we consider the
bitraction test specimen presented in figure \ref{fig:croix} (this
structure, courtesy of ONERA -- Pascale Kanout\'e --, was optimized
with ZeBuLoN software in order to have stress field as homogeneous
as possible in its center). It was decomposed with Metis software
into 4 or 16 subdomains.

\begin{table}[ht]\centering
\begin{tabular}[b]{|l|l|c|c|}
\hline \multicolumn{2}{|c|}{Method} & 4 sd. & 16 sd.   \\
\hline \multicolumn{2}{|l|}{Primal Neumann-Neumann} & 23$_{0+12}$
& 30$_{0+69}$
\\\hline
\multirow{3}{4cm}{Dual}
& Lumped P(I) & 32$_{12}$  & 41$_{69}$  \\
& Dirichlet P(I) & 25$_{12}$  & 32$_{69}$  \\
& Dirichlet P(Q) & 24$_{12}$  & 32$_{69}$  \\\hline
\multirow{2}{3cm}{Hybrid P(Q)}
& P-P-D & 31$_{8}$& 44$_{46}$    \\
& D-D-P & 25$_{8}$& 37$_{46}$   \\\hline
\end{tabular}\caption{Bitraction test specimen}\label{tab:croix}
\end{table}

As shown in table \ref{tab:croix} all methods give excellent
performance results on non-academical problem. Note the good
behavior of the hybrid approach which gives equivalent performance
with much smaller coarse problem, even if no physical
consideration could guide the choice of the treatment of interface
degrees of freedom.

\section{Conclusions}\label{sec:conc}
In this paper, we have reviewed most used non-overlapping domain
decomposition methods. These methods are perfectly suited to
modern computational hardware, they are based on very close
concepts which we tried to outline. We introduced the hybrid
framework to include as many methods as possible: the principle is
to assign to each interface degree of freedom its own treatment,
for now primal and dual treatments have been implemented, and
mixed and recondensed should follow. Hybriding methods also enable
to define physic-friendly approaches for multifield problems.

Because of the conceptual proximity of all methods, assessments
showed very close numerical performance results. Once equipped
with convenient preconditioner and coarse problem, all the methods 
proved their ability to handle second and forth order
elasticity in presence of strong heterogeneities. Though from a
computational point of view some combination may be more
interesting: dual approach with lumped preconditioner or
simplified projector (if these are sufficient to ensure fine rate
of convergence), hybrid approach (which generates smaller coarse
space). We also outlined the importance of the decomposition even
for very simple problems. The methods have also proved their
efficiency on industrial cases, some of them were implemented in
computational softwares.

In this paper we limited to the solution to linearized systems,
which anyhow enables to solve nonlinear problems. Another strategy
is to commute the nonlinear solver and the domain decomposition
method so that nonlinear problems can be solved independently on
each subdomain. Another evolution of domain decomposition
philosophy is the decomposition of the time interval
\cite{LIONS:2001:RED} for evolution problems.

\section*{ACKNOWLEDGEMENTS}
We wish to thank Frederic Feyel and Daniel Rixen for coauthoring
some of our works. More specifically (but not exclusively)
Frederic for his contribution when implementing our methods inside
ZeBuLoN software. Moreover, working with them is always a great
pleasure. First author wish to thank ECCOMAS and ARCME for their
kind proposal to write this paper.
\bibliographystyle{plain}
\bibliography{pi}

\begin{thebibliography}{100}

\bibitem{ACHDOU:1995:SPF}
Yves Achdou and Yuri~A. Kuznetsov.
\newblock Substructuring preconditioners for finite element methods on
  nonmatching grids.
\newblock {\em East-West J. Numer. Math.}, 3(1):1--28, 1995.

\bibitem{ACHDOU:1996:MIS}
Yves Achdou, Yvon Maday, and Olof~B. Widlund.
\newblock M\'ethode it\'erative de sous-structuration pour les \'el\'ements
  avec joints.
\newblock {\em C.R. Acad. Sci. Paris}, I(322):185--190, 1996.

\bibitem{BARBOTEU:2001:DDS}
Mickael Barboteu, Pierre Alart, and Marina Vidrascu.
\newblock A domain decomposition strategy for nonclassical frictional
  multicontact problems.
\newblock {\em Comp. Meth. App. Mech. Eng.}, 190:4785--4803, 2001.

\bibitem{BARRET:1994:TEM}
Richard Barrett, Michael Berry, Tony~F. Chan, James Demmel, June Donato, Jack
  Dongarra, Victor Eijkhout, Roldan Pozo, Charles Romine, and Henk~Van der
  Vorst.
\newblock {\em Templates for the Solution of Linear Systems: Building Blocks
  for Iterative Methods}.
\newblock SIAM, 1994.

\bibitem{BELYTSCHKO:1996:MMO}
T.~Belytschko, Y.~Krongauz, D.~Organ, M.~Fleming, , and P.~Krysl.
\newblock Meshless methods: An overview and recent developments.
\newblock {\em Computer Methods in Applied Mechanics and Engineering},
  139:3--47, 1996.

\bibitem{BERNARDI:1990:MORT}
C.~Bernardi, Y.~Maday, and T.~Patera.
\newblock A new non conforming approach to domain decomposition: the {Mortar
  Element Method}.
\newblock In H.~Brezis and J.L. Lions, editors, {\em {N}onlinear Partial
  Differential Equations and their Applications}. Pitman, London, 1989.

\bibitem{BHARDWAJ:1999:PQ_LRV}
M.~Bhardwaj, D.~Day, C.~Farhat, M.~Lesoinne, K.~Pierson, and D.~Rixen.
\newblock Application of the {FETI} method to {ASCI} problems: Scalability
  results on a thousand-processor and discussion of highly heterogeneous
  problems.
\newblock {\em Int. J. Num. Meth . Eng.}, 47(1-3):513--536, 2000.

\bibitem{BOLANDER:2005:ILM}
J.~E. Bolander and N.~Sukumar.
\newblock Irregular lattice model for quasistatic crack propagation.
\newblock {\em Phys. Rev. B}, 71, 2005.

\bibitem{BREITKOPT:2002:MPA}
Piotr Breitkopt and Antonio Huerta, editors.
\newblock {\em Meshfree and particle based approaches in computational
  Mechanics}, volume~11 of {\em REEF special release}.
\newblock Hermes, 2002.

\bibitem{BRENNER:2005:LBD}
Suzanne Brenner.
\newblock Lower bounds in domain decomposition.
\newblock In {\em Proceedings of the 16th international conference on domain
  decomposition methods}, 2005.

\bibitem{BREZZI:1993:3FIELD}
F.~Brezzi and L.D. Marini.
\newblock A three-field domain decomposition method.
\newblock In {\em Proceedings of the sixth international conference on domain
  decomposition methods}, pages 27--34, 1993.

\bibitem{WALSHAW:1995:PAO}
M.~Cross C.~Walshaw and M.~G. Everett.
\newblock A parallelisable algorithm for optimising unstructured mesh
  partitions.
\newblock Technical report, School of Mathematics, Statistics \& Scientific
  Computing, University of Greenwich, London, 1995.

\bibitem{CAI:2002:NPI}
X-C Cai and David Keyes.
\newblock Nonlinearly preconditioned inexact newton algorithms.
\newblock {\em SIAM J. Sci. Comp.}, 24:183--200, 2002.

\bibitem{CHAPMAN:1996:DAK}
A.~Chapman and Y.~Saad.
\newblock Deflated and augmented {K}rylov subspace techniques.
\newblock {\em Numer. Linear Algebra Appl.}, 1997.

\bibitem{CIARLET:1979:FEM}
P.G. Ciarlet.
\newblock {\em The finite element method for elliptic problems}.
\newblock North Holland, 1979.

\bibitem{SALENCON:1988:MMC}
J.~Salen\c con.
\newblock {\em M\'ecanique des milieux continus}.
\newblock Ecole polytechnique. Ellipses, 1988.

\bibitem{CRAIG:1968:CSD}
R.~Craig and M.~Bampton.
\newblock Coupling of substructures for dynamic analysis.
\newblock {\em AIAA Journal}, 6:1313--1319, 1968.

\bibitem{CRESTA:2005:CMP}
Philippe Cresta, Olivier Allix, Christian Rey, and St\'ephane Guinard.
\newblock Comparison of multiscale and parallel nonlinear strategies based on
  domain decomposition for post buckling analysis.
\newblock {\em Comp. Meth. App. Mech. Eng.}, 196(8):1436--1446, 2005.

\bibitem{CROS:2002:PSC}
Jean-Michel Cros.
\newblock A preconditioner for the schur complement domain decomposition
  method.
\newblock In Herrera, Keyes, and Widlund, editors, {\em Proceedings of the
  $14^{th}$ international conference on domain decomposition methods}, pages
  373--380, 2002.

\bibitem{DADDETTA:2004:PCH}
G.A. D'Addetta, E.~Ramm, S.~Diebels, and W.~Ehlers.
\newblock A particle center based homogenization strategy for granular
  assemblies.
\newblock {\em Int J for Computer-Aided Engineering}, 21(2-4):360--383, 2004.

\bibitem{DELABOURDONNAYE:1998:FETIH}
A.~de~La~Bourdonnaye, C.~Farhat, A.~Macedo, F.~Magoules, and F.-X. Roux.
\newblock {\em Advances in Computational Mechanics with High Performance
  Computing}, chapter A method of finite element tearing and interconnecting
  for the Helmholtz problem, pages 41--54.
\newblock Civil-Comp Press, Edinburgh, United Kingdom, 1998.

\bibitem{DELAPLACE:2005:FDF}
Arnaud Delaplace.
\newblock Fine description of fracture by using discrete particle model.
\newblock In {\em Proceedings of ICF 11 - 11th International Conference on
  Fracture}, 2005.

\bibitem{DOHRMANN:2003:PSB}
Clarck Dohrmann.
\newblock A preconditioner for substructuring based on constrained energy
  minimization.
\newblock {\em SIAM J. Sci. Comp.}, 25(1):246--258, 2003.

\bibitem{dolean-2005}
Victorita Dolean, Frederic Nataf, and Gerd Rapin.
\newblock New constructions of domain decomposition methods for systems of
  pdes.
\newblock {\em C. R. Acad. Sci. Paris}, 340(1):693--696, 2005.

\bibitem{DOSTAL:1988:CGM}
Z.~Dostal.
\newblock Conjugate gradient method with preconditioning by projector.
\newblock {\em Int. J. Comput. Math.}, 23:315--323, 1988.

\bibitem{DOSTAL:2005:OSF}
Zdenek Dostal, David Horak, and Dan Stefanica.
\newblock An overview of scalable feti-dp algorithms for variational
  inequalities.
\newblock In {\em Proceedings of the 16th conference on domain decomposition
  methods}, 2005.

\bibitem{DUREISSEIX:2001:CONTACT}
D.~Dureisseix and C.~Farhat.
\newblock A numerically scalable domain decomposition method for the solution
  of frictionless contact problems.
\newblock {\em Internat. J. Num. Meth. Engin.}, 50(12):2643--2666, 2001.

\bibitem{DUVAUT:1990:MMC}
Georges Duvaut.
\newblock {\em M\'ecanique des milieux continus}.
\newblock Masson, 1990.

\bibitem{ERHEL:1995:RGM}
J.~Erhel, K.~Burrage, and B.~Pohl.
\newblock Restarted {GMR}es preconditioned by deflation.
\newblock {\em J. Comput. Appl. Math.}, 69:303--318, 1996.

\bibitem{ERHEL:2000:ACG}
Justine Erhel and F.~Guyomarc'h.
\newblock An augmented conjugate gradient method for solving consecutive
  symmetric positive definite linear systems.
\newblock {\em SIAM J. Matrix Anal. Appl.}, 21(4):1279--1299, 2000.

\bibitem{FARHAT:1992:SADDLE}
C.~Farhat.
\newblock A saddle-point principle domain decomposition method for the solution
  of solid mechanics problems.
\newblock In D.~Keyes, T.F. Chan, G.A. Meurant, J.S. Scroggs, and R.G. Voigt,
  editors, {\em Domain Decomposition Methods for Partial Differential
  Equations}, pages 271--292, 1992.

\bibitem{FARHAT:1996:COR2}
C.~Farhat, P.-S. Chen, and F.-X. Roux.
\newblock The two-level {FETI} method - part {II}: Extension to shell problems.
  parallel implementation and performance results.
\newblock {\em J. Comp Meth. Appl. Mech. Eng.}, 155:153--180, 1998.

\bibitem{FARHAT:1996:CNS}
C.~Farhat and M.~G\'eradin.
\newblock On the computation of the null space and generalized invers of large
  matrix, and the zero energy modes of a structure.
\newblock Technical Report CU-CAS-96-15, Center for aerospace structures, may
  1996.

\bibitem{FARHAT:2001:FETI_DP}
C.~Farhat, M.~Lesoinne, P.~LeTallec, K.~Pierson, and D.~Rixen.
\newblock {FETI-DP}: a dual-primal unified {FETI} method - part i: a faster
  alternative to the two-level {FETI} method.
\newblock {\em Int. J. Num. Meth . Eng.}, 50(7):1523--1544, 2001.

\bibitem{FARHAT:2000:FETI_DP}
C.~Farhat, M.~Lesoinne, and K.~Pierson.
\newblock A scalable dual-primal domain decomposition method.
\newblock {\em Numer. Linear Algebra Appl.}, 7(7-8):687--714, 2000.

\bibitem{FARHAT:2000:TLD}
C.~Farhat, A.~Macedo, and M.~Lesoinne.
\newblock A two-level domain decomposition method for the iterative solution of
  high frequency exterior {H}elmholtz problems.
\newblock {\em Numerische Mathematik}, 85:283--308, 2000.

\bibitem{FARHAT:1999:FETIH}
C.~Farhat, A.~Macedo, and R.~Tezaur.
\newblock {FETI-H}: a scalable domaine decomposition method for high frequency
  exterior {H}elmholtz problems.
\newblock In {\em Domain decomposition methods in science and engineering}.
  Domain decomposition press, 1999.

\bibitem{FARHAT:1996:COR1}
C.~Farhat and J.~Mandel.
\newblock The two-level {FETI} method for static and dynamic plate problems -
  part {I}: An optimal iterative solver for biharmonic systems.
\newblock {\em J. Comp Meth. Appl. Mech. Eng.}, 155:129--152, 1998.

\bibitem{FARHAT:1994:OPT}
C.~Farhat, J.~Mandel, and F.X. Roux.
\newblock Optimal convergence properties of the {FETI} domain decomposition
  method.
\newblock {\em Comp. Meth. Appl. Mech. Eng.}, 115:365--385, 1994.

\bibitem{FARHAT:2000:FETINL}
C.~Farhat, K.~Pierson, and M.~Lesoinne.
\newblock The second generation {FETI} methods and their application to the
  parallel solution of large-scale linear and geometrically non-linear
  structural analysis problems.
\newblock {\em J. Comp Meth. Appl. Mech. Eng.}, 184(2-4):333--374, 2000.

\bibitem{FARHAT:1991:FETI}
C.~Farhat and F.-X. Roux.
\newblock A method of finite tearing and interconnecting and its parallel
  solution algorithm.
\newblock {\em Int. J. Num. Meth . Eng.}, 32:1205--1227, 1991.

\bibitem{FARHAT:1994:RFETI}
C.~Farhat and F.-X. Roux.
\newblock The dual {S}chur complement method with well-posed local {N}eumann
  problems.
\newblock {\em Contemporary Mathematics}, 157:193--201, 1994.

\bibitem{FARHAT:1994:ADV}
C.~Farhat and F.~X. Roux.
\newblock Implicit parallel processing in structural mechanics.
\newblock {\em Computational Mechanics Advances}, 2(1):1--124, 1994.
\newblock North-Holland.

\bibitem{FARHAT:1993:TOP}
C.~Farhat and H.~D. Simon.
\newblock Top/domdec - a software tool for mesh partitioning and parallel
  processing.
\newblock Technical report, NASA Ames, 1993.

\bibitem{FEYEL:1998:THESE}
Fr\'ed\'eric Feyel.
\newblock {\em Application du calcul parall\`ele aux mod\`eles \`a grand nombre
  de variables internes}.
\newblock Th\`ese de doctorat, Ecole Nationale Sup\'erieure des Mines de Paris,
  1998.

\bibitem{FEYEL:2005:HDR}
Fr\'ed\'eric Feyel.
\newblock {\em Quelques multi-probl\`emes en m\'ecanique des mat\'eriaux et
  structures}.
\newblock Habilitation \`a diriger des recherches, Universit\'e Pierre et Marie
  Curie, 2005.

\bibitem{FORTIN:1982:MLA}
Michel Fortin and Roland Glowinski.
\newblock {\em M\'ethodes de lagrangien augment\'e - applications \`a la
  r\'esolution num\'erique de probl\`emes aux limites}.
\newblock Dunod, 1982.

\bibitem{FRAGAKIS:2002:UFF}
Yannis Fragakis and Manolis Papadrakakis.
\newblock A unified framework for formulating domain decomposition methods in
  structural mechanics.
\newblock Technical report, Institute of Structural Analysis and Seismic
  Research, National Technical University of Athens, 2002.

\bibitem{FRAGAKIS:2004:MHP1}
Yannis Fragakis and Manolis Papadrakakis.
\newblock The mosaic of high-performance domain decomposition methods for
  structural mechanics -- part i: Formulation, interrelation and numerical
  efficiency of primal and dual methods.
\newblock {\em Comp. Meth. Appl. Mec. Eng.}, 192(35-36):3799--3830, 2003.

\bibitem{FRAGAKIS:2004:MHP2}
Yannis Fragakis and Manolis Papadrakakis.
\newblock The mosaic of high-performance domain decomposition methods for
  structural mechanics -- part ii: Formulation enhancements, multiple
  right-hand sides and implicit dynamics.
\newblock {\em Comp. Meth. Appl. Mec. Eng.}, 193(42-44):4611--4662, 2004.

\bibitem{FRAYSSE:1998:IOS}
V.~Frayss\'e, L.~Giraud, and H.~Kharraz-Aroussi.
\newblock On the influence of the orthogonalization scheme on the parallel
  performance of {GMR}es.
\newblock Technical report, {CERFACS}, 1998.

\bibitem{GERMAIN:2005:MCNL}
Norbert Germain, Jacques Besson, and Fr\'ed\'eric Feyel.
\newblock M\'ethodes de calcul non local : Application aux structures
  composites.
\newblock In {\em Actes du 7\`eme colloque national en calcul des structures},
  Giens, 2005.

\bibitem{GERMAIN:1986:M}
P.~Germain.
\newblock {\em M\'ecanique}.
\newblock Ecole polytechnique. Ellipses, 1986.

\bibitem{GLO90}
R.~Glowinski and P.~Le Tallec.
\newblock Augmented lagrangian interpretation of the nonoverlapping {S}chwartz
  alternating method.
\newblock In {\em Third International Symposium on Domain Decomposition Methods
  for Partial Differential Equations}, pages 224--231, 1990.

\bibitem{GOLDFELD:2002:BNN}
P.~Goldfeld.
\newblock Balancing {N}eumann-{N}eumann for (in)compressible linear elasticity
  and (generalized) {S}tokes -- parallel implementation.
\newblock In {\em Proceedings of the $14^{th}$ international conference on
  domain decomposition method}, pages 209--216, 2002.

\bibitem{GOSSELET:2003:AHD}
P.~Gosselet, V.~Chiaruttini, C.~Rey, and F.~Feyel.
\newblock Une approche hybride de d\'ecomposition de domaine pour les
  probl\`emes multiphysiques : application \`a la poro\'elasticit\'e.
\newblock In {\em Actes du sixi\`eme colloque national en calcul de
  structures}, volume~2, pages 297--304, 2003.

\bibitem{GOSSELET:2004:DDM}
P.~Gosselet, Vincent Chiaruttini, C.~Rey, and F.~Feyel.
\newblock A monolithic strategy based on an hybrid domain decomposition method
  for multiphysic problems. application to poroelasticity.
\newblock {\em Revue europ\'eenne des \'elements finis}, 13:523--534, 2004.

\bibitem{GOSSELET:2002:SRK}
P.~Gosselet and C.~Rey.
\newblock On a selective reuse of krylov subspaces in newton-krylov approaches
  for nonlinear elasticity.
\newblock In {\em Proceedings of the $14^{th}$ conference on domain
  decomposition methods}, pages 419--426, 2002.

\bibitem{GOSSELET:2002:DDM}
P.~Gosselet, C.~Rey, P.~Dasset, and F.~L\'en\'e.
\newblock A domain decomposition method for quasi incompressible formulations
  with discontinuous pressure field.
\newblock {\em Revue europ\'eenne des \'elements finis}, 11:363--377, 2002.

\bibitem{GOSSELET:2003:IEI}
P.~Gosselet, C.~Rey, and D.~Rixen.
\newblock On the initial estimate of interface forces in {FETI} methods.
\newblock {\em Comp. meth. appl. mech. engrg.}, 192:2749--2764, 2003.

\bibitem{GOSSELET:2003:THESE}
Pierre Gosselet.
\newblock {\em M\'ethodes de d\'ecomposition de domaine et m\'ethodes
  d'acc\'el\'eration pour les probl\`emes multichamp en m\'ecanique
  non-lin\'eaire}.
\newblock PhD thesis, Universit\'e P. et M. Curie, 2003.

\bibitem{KARYPIS:1998:MAM}
George Karypis and Vipin Kumar.
\newblock Multilevel algorithms for multi-constraint graph partitioning.
\newblock Technical report, University of Minnesota, Department of Computer
  Science, 1998.

\bibitem{KLAWONN:2005:SCR}
A.~Klawonn, O.~Rheinbach, and O.B. Widlund.
\newblock Some computational results for dual-primal feti methods for three
  dimensional elliptic problems.
\newblock {\em Lect. Notes Comput. Sci. Eng.}, 40:361 -- 368, 2005.

\bibitem{KLAWONN:1999:DP}
A.~Klawonn and O.B. Widlund.
\newblock Dual and dual-primal {FETI} methods for elliptic problems with
  discontinuous coefficients.
\newblock In {\em Proceedings of the 12th International Conference on Domain
  Decomposition Methods, Chiba, Japan}, October 1999.

\bibitem{KLAWONN:2001:FNN}
A.~Klawonn and O.B. Widlund.
\newblock {FETI} and {N}eumann-{N}eumann iterative substructuring methods:
  connections and new results.
\newblock {\em Comm. pure and appl. math.}, LIV:0057--0090, 2001.

\bibitem{LAD99b}
P.~Ladev\`eze.
\newblock {\em Nonlinear Computational Structural Mechanics - New Approaches
  and Non-Incremental Methods of Calculation}.
\newblock Springer Verlag, 1999.

\bibitem{LAD01}
P.~Ladev\`eze, O.~Loiseau, and D.~Dureisseix.
\newblock A micro-macro and parallel computational strategy for highly
  heterogeneous structures.
\newblock {\em Int. J. Num. Meth. Engnrg.}, 52:121--138, 2001.

\bibitem{LADEVEZE.2006.1}
P.~Ladev\`eze, D.~N\'eron, and P.~Gosselet.
\newblock On a mixed and multiscale domain decomposition method.
\newblock {\em Computer Methods in Applied Mechanics and Engineering},
  196(8):1526--1540, 2006.

\bibitem{Lesoinne:1999:FES}
Michel Lesoinne and Kendall Pierson.
\newblock {FETI-DP}: An efficient, scalable, and unified {D}ual-{P}rimal {FETI}
  method.
\newblock In {\em Domain Decomposition Methods in Sciences and Engineering},
  pages 421--428, 1999.

\bibitem{LI:2002:DPF}
Jing Li.
\newblock A dual-primal feti method for solving stokes/navier-stokes equations.
\newblock In {\em Proceedings of the $14^{th}$ international conference on
  domain decomposition method}, pages 225--233, 2002.

\bibitem{LINGEN:2000:EGS}
F.J. Lingen.
\newblock Efficient {G}ram-{S}cmidt orthonormalisation on parallel computers.
\newblock {\em Com. Numer. Meth. Engng.}, 16:57--66, 2000.

\bibitem{LIONS:2001:RED}
JL. Lions, Y.~Maday, and G.~Turinici.
\newblock R\'esolution d'edp par un sch\'ema en temps parar\'eel.
\newblock {\em C. R. Acad. Sci. Paris}, 333(1):1--6, 2001.

\bibitem{LIU:2005:IMM}
Gui-Rong Liu and Yuan-Tong Gu.
\newblock {\em An Introduction to Meshfree Methods and Their Programming}.
\newblock Springer, 2005.

\bibitem{MANDEL:1993:BAL}
J.~Mandel.
\newblock Balancing domain decomposition.
\newblock {\em Comm. Appl. Num. Meth. Engrg.}, 9:233--241, 1993.

\bibitem{MANDEL:1996:COEF}
J.~Mandel and M.~Brezina.
\newblock Balancing domain decomposition for problems with large jumps in
  coefficients.
\newblock {\em Math. Comp.}, 65(216):1387--1401, 1996.

\bibitem{MANDEL:1996:OPT}
J.~Mandel and R.~Tezaur.
\newblock Convergence of a substructuring method with {L}agrange multipliers.
\newblock {\em Numerische Mathematik}, 73:473--487, 1996.

\bibitem{MANDEL:2000:DP}
J.~Mandel and R.~Tezaur.
\newblock On the convergence of a dual-primal substructuring method.
\newblock UCD/CCM Report 150, Center for Computational Mathematics, University
  of Colorado at Denver, April 2000.
\newblock to appear in Numer. Math.

\bibitem{ZEBUDEVEL:2001}
Northwest Numerics.
\newblock {\em Z-set developper manual}, 2001.

\bibitem{ZEBUUSER:2001}
Northwest Numerics.
\newblock {\em Z-set user manual}, 2001.

\bibitem{NOUY:2003:THES}
A.~Nouy.
\newblock {\em Une strat\'egie de calcul multi\'echelle avec
  homog\'en\'eisation en temps et en espace pour le calcul de structures
  fortement h\'et\'erog\`enes}.
\newblock PhD thesis, ENS de Cachan, 2003.

\bibitem{PARK:1997:ALG}
K.C. Park, M.R. Justino, and C.A. Felippa.
\newblock An algebraically partitioned {FETI} method for parallel structural
  analysis: algorithm description.
\newblock {\em Int. J. Num. Meth. Eng.}, 40(15):2717--2737, 1997.

\bibitem{PARK:1997:PERF}
K.C. Park, M.R. Justino, and C.A. Felippa.
\newblock An algebraically partitioned {FETI} method for parallel structural
  analysis: performance evaluation.
\newblock {\em Int. J. Num. Meth . Eng.}, 40(15):2739--2758, 1997.

\bibitem{PAZ:2005:ISP}
Rodrigo Paz and Mario Storti.
\newblock An interface strip preconditioner for domain decomposition methods:
  application to hydrology.
\newblock {\em Int. J. Numer. Meth. Engng.}, 62:1873--1894, 2005.

\bibitem{REY:2003:SLN}
C.~Rey and P.~Gosselet.
\newblock Solution to large nonlinear systems: acceleration strategies based on
  domain decomposition and reuse of krylov subspaces.
\newblock In {\em Proceedings of the $6^th$ {ESAFORM} conference on material
  forming}, 2003.

\bibitem{REY:1998:RKS}
C.~Rey and F.~L\'en\'e.
\newblock Reuse of krylov spaces in the solution of large-scale non linear
  elasticity problems.
\newblock In {\em Domain Decomposition Methods in Sciences and Engineering},
  pages 465--471, 1998.

\bibitem{Rey:1998:RRP}
C.~Rey and F.~Risler.
\newblock A {R}ayleigh-{R}itz preconditioner for the iterative solution to
  large scale nonlinear problems.
\newblock {\em Numerical Algorithms}, 17:279--311, 1998.

\bibitem{Risler:2000:IAA}
Franck Risler and Christian Rey.
\newblock Iterative accelerating algorithms with {K}rylov subspaces for the
  solution to large-scale nonlinear problems.
\newblock {\em Numerical Algorithms}, 23:1--30, 2000.

\bibitem{RISLER:1998:RRV}
Frank Risler and Christian Rey.
\newblock On the reuse of {R}itz vectors for the solution to nonlinear
  elasticity problems by domain decomposition methods.
\newblock In {\em DD10 Proceedings, Contemporary Mathematics}, volume 218,
  pages 334--340, 1998.

\bibitem{RIXEN:1997:THESE}
D.~Rixen.
\newblock {\em Substructuring and dual methods in structural analysis}.
\newblock PhD thesis, University of Li\`ege, Belgique, 1997.

\bibitem{RIXEN:2004:DCB}
D.~Rixen.
\newblock A dual craig-bampton method for dynamic substructuring.
\newblock {\em J. Comput. Appl. Math.}, 168:383--391, 2004.

\bibitem{RIXEN:1998:SUPERL}
D.~Rixen and C.~Farhat.
\newblock A simple and efficient extension of a class of substructure based
  preconditioners to heterogeneous structural mechanics problems.
\newblock {\em Int. J. Num. Meth. Eng.}, 44(4):489--516, 1999.

\bibitem{RIXEN:1999:AFETI}
D.~Rixen, C.~Farhat, R.~Tezaur, and J.~Mandel.
\newblock Theoretical comparison of the {FETI} and algebraically partitioned
  {FETI} methods, and performance comparisons with a direct sparse solver.
\newblock {\em Int. J. Num. Meth. Eng.}, 46(4):501--534, 1999.

\bibitem{RIXEN:2002:EPF}
Daniel Rixen.
\newblock Extended preconditioners for the feti method applied to constrained
  problems.
\newblock {\em Int. Journal for Numerical methods in engineering}, 54(1):1--26,
  2002.

\bibitem{ROUX:1997:DIRECT2}
F.-X. Roux.
\newblock Parallel implementation of direct solution strategies for the coarse
  grig solvers in 2-level {FETI} method.
\newblock Technical report, {ONERA}, Paris, France, 1997.

\bibitem{SAAD:1997:AAK}
Y.~Saad.
\newblock Analysis of augmented {K}rylov subspace methods.
\newblock {\em SIAM J. Matrix Anal. Appl.}, 18(2):435--449, April 1997.

\bibitem{SAAD:2000:IMS}
Y.~Saad.
\newblock {\em Iterative methods for sparse linear systems}.
\newblock PWS Publishing Company, 3rd edition, 2000.

\bibitem{SAAD:1986:GMR}
Y.~Saad and M.~H. Schultz.
\newblock {GMR}es: a generalized minimal residual algorithm for solving
  nonsymmetric linear systems.
\newblock {\em SIAM J. Sci. Comput.}, 7:856--869, 1986.

\bibitem{SAAD:1987:LMS}
Youssef Saad.
\newblock On the {L}anczos method for solving symmetric linear systems with
  several right hand sides.
\newblock {\em Math. Comp.}, 48:651--662, 1987.

\bibitem{SERIES:2003:MDD2}
L.~Series, F.~Feyel, and F.-X. Roux.
\newblock Une m\'ethode de d\'ecomposition de domaine avec deux multiplicateurs
  de {L}agrange.
\newblock In {\em Actes du $16^{eme}$ congr\`es fran\c cais de m\'ecanique},
  2003.

\bibitem{SERIES:2003:MDD1}
L.~Series, F.~Feyel, and F.-X. Roux.
\newblock Une m\'ethode de d\'ecomposition de domaine avec deux multiplicateurs
  de {L}agrange, application au calcul des structures, cas du contact.
\newblock In {\em Actes du sixi\`eme colloque national en calcul des
  structures}, volume III, pages 373--380, 2003.

\bibitem{STEFANICA:1999:FME}
D.~Stefanica and A.~Klawonn.
\newblock The {FETI} method for mortar finite elements.
\newblock In {\em Proceedings of 11th International Conference on Domain
  Decomposition Methods}, pages 121--129, 1999.

\bibitem{LETALLEC:1994:DDM}
P.~Le {T}allec.
\newblock Domain-decomposition methods in computational mechanics.
\newblock {\em Computational Mechanics Advances}, 1(2):121--220, 1994.
\newblock North-Holland.

\bibitem{LETALLEC:1997:SHELL}
P.~Le {T}allec, J.~Mandel, and M.~Vidrascu.
\newblock A {N}eumann-{N}eumann domain decomposition algorithm for solving
  plate and shell problems.
\newblock {\em SIAM J. Num. Ana.}, 35(2):836--867, April 1998.

\bibitem{LETALLEC:1991:DDM}
P.~Le {T}allec, Y.-H.~De {R}oeck, and M.~Vidrascu.
\newblock Domain-decomposition methods for large linearly elliptic three
  dimensional problems.
\newblock {\em J. Comp. Appl. Math.}, 34:93--117, 1991.
\newblock Elsevier Science Publishers, Amsterdam.

\bibitem{LETALLEC:1993:DDM}
P.~Le {T}allec and M.~Vidrascu.
\newblock M\'ethodes de d\'ecomposition de domaines en calcul de structures.
\newblock In {\em Actes du premier colloque national en calcul des structures},
  volume~I, pages 33--49, 1993.

\bibitem{LETALLEC:1996:BERG}
P.~Le {T}allec and M.~Vidrascu.
\newblock Generalized {N}eumann-{N}eumann preconditioners for iterative
  substructuring.
\newblock In {\em Proceedings of the ninth conference on Domain Decomposition},
  June Bergen 1996.
\newblock to appear.

\bibitem{VANDERSLUIS:1986:RCCG}
A.~van~der Sluis and H.~can~der Vorst.
\newblock The rate of convergence of conjugate gradients.
\newblock {\em Numer. Math.}, 48:543--560, 1986.

\bibitem{WESSELING:2004:IMM}
Pieter Wesseling.
\newblock {\em An Introduction to Multigrid Methods}.
\newblock R.T. Edwards, Inc, 2004.

\bibitem{ZIENKIEWICZ:1989:FEM}
O.C. Zienkiewicz and R.L. Taylor.
\newblock {\em The finite element method}.
\newblock Mc Graw-Hill Book COmpagny, 1989.

\end{thebibliography}
\newpage
\appendix
\section{Krylov iterative solvers}\label{sec:solvers}
\def\thesubsection{\Alph{section}.\arabic{subsection}}
Krylov iterative solvers for the resolution of linear systems have
been widely studied. The aim of this section is just to briefly
present important results and algorithms, reader interested in
wider documentation can refer to \cite{SAAD:2000:IMS}, and to
\cite{BARRET:1994:TEM} for shorter explanation.

Krylov methods belong to the projection class of iterative
algorithms, which consist in approximating solution $S^{-1}b$ of
system $Sx=b$ by vector $p(S)b$ where $p$ is a smartly built
polynomial.

In this section we consider the iterative solution to system
$Sx=b$. $S$ is a $n\times n$ matrix and $b$ a vector in
$\range(S)$. The $i^{th}$ iteration leads to approximation $x_i$
of the solution, associated residual is $r_i=b-Sx_i=S(x-x_i)$.
Initialization is $x_0$ (most often $x_0=0$). Canonical
(orthonormal) basis of $\mathbb{R}^n$ reads $(e_1,\ldots,e_n)$.

\subsection{Principle of Krylov solvers}
Krylov solvers are based on the iterative construction of
so-called "Krylov subspace" $\mathcal{K}_m(S,r_0)$ defined by:
\begin{equation}\label{eq:def_krylov:1}
    \mathcal{K}_m(S,r_0)=\vect\left(r_0,\ldots,S^{m-1}r_0\right)
\end{equation}

The solution to linear system consists in searching $x_m$ under the following constraints:
\begin{equation}\label{eq:def_krylov:2}
    \left\{ \begin{array}{l} x_m \in x_0 + \mathcal{K}_m(S,r_0) \\
    r_m \perp_? \mathcal{K}_m(S,r_0) \end{array} \right.
\end{equation}
where the choice of the orthogonality relationship enables to
define various approaches.

\subsection{Most used solvers}
We herein present two of the principal Krylov solvers. First GMRes
\cite{SAAD:1986:GMR} which is suited to any type of matrix, then
conjugate gradient which is adapted to symmetric definite positive
matrices.

Of course, the iterative solution to a linear system assumes that
a convergence criterion is employed, and that a limit of precision
is set so that the system is supposed to have converged once the
criterion is below this precision. We note $\varepsilon$ this
limit value of the criterion.

\subsection{GMRes}
Algorithm GMRes (alg. \ref{alg:gmres:1}) consists in an oblique
projection based on the construction of Krylov subspace
$\mathcal{K}_m(S,v_0)$ with $v_0=r_0/\|r_0\|_2$. The research
principle is:
\begin{equation}\label{eq:def_krylov:4}
    \left\{ \begin{array}{l} x_m \in x_0 + \mathcal{K}_m(S,r_0) \\
    r_m \perp S \mathcal{K}_m(S,r_0) \end{array} \right.
\end{equation}
which is equivalent to finding $x_m\in x_0 + \mathcal{K}_m(S,r_0)$
minimizing $\|r_m\|_2$.

\begin{algorithm}\caption{GMRes}\label{alg:gmres:1}
\begin{algorithmic}[1]
\STATE Compute $r_0=b-Sx_0$, $v_0=r_0/\|r_0\|_2$
\FOR{$j=0,\ldots,m-1$}%
  \STATE Compute $w_j=Sv_j$%
  \FOR{$i=0,\ldots,j$}%
    \STATE $h_{ij}=(v_i,w_j)$
    \STATE $w_j=w_j-h_{ij}v_i$%
  \ENDFOR%
  \STATE $h_{(j+1)j}=\|w_j\|_2$%
  \IF{$\|r_j\|_2 \leqslant \varepsilon$} \STATE stop \ELSE \STATE $v_{j+1}=w_j/h_{(j+1) j}$ \ENDIF%
\ENDFOR%
\STATE Compute $y_m$ minimizing $\left\| \|r_0\|_2 e_1 - \bar{H}_m
y\right\|_2$ and set $x_m=x_0+V_m y_m$
\end{algorithmic}
\end{algorithm}

A particulary striking property of GMRes is not to compute the
approximation at each iteration, a smart implementation of GMRes
enables to directly access the norm of the residual $\|r_{j}\|_2$.
Only the final approximation $x_m$ is computed (by the inversion
of a $m\times m$ upper triangular matrix). From the computation
complexity point of view, each iteration consists in a full
orthonormalization of vector $w_j$ with respect to
$\mathcal{K}_j$.

Algorithm GMRes(m) (or restarted GMRes) consists in stopping
computation before convergence at \textit{a priori} fixed step $m$
and restarting computation using previous $x_m$ as initialization.
This strategy aims at minimizing orthogonalization computations by
limiting the size of Krylov subspaces \cite{ERHEL:1995:RGM}. This
method may lead to stagnation for non positive definite matrices.

\subsection{Conjugate gradient}
Let $S$ be a symmetric positive definite matrix, a conjugate
gradient algorithm consists in an orthogonal projection. The
research principle is:
\begin{equation}\label{eq:def_krylov:5}
    \left\{ \begin{array}{l} x_m \in x_0 + \mathcal{K}_m(S,r_0) \\
    r_m \perp \mathcal{K}_m(S,r_0) \end{array} \right.
\end{equation}
which is equivalent to finding $x_m\in x_0 + \mathcal{K}_m(S,r_0)$
minimizing $\|x_m-x\|_S$.

Because of the properties of $S$, conjugation (orthogonality) properties appear, leading to algorithm
\ref{alg:cg:1}.
\begin{algorithm}\caption{Conjugate gradient}\label{alg:cg:1}
\begin{algorithmic}[1]
\STATE Compute $r_0=b-Sx_0$, set $w_0=r_0$%
\FOR{$j=0,\ldots,m$} %
  \STATE $\alpha_j=(r_j,r_j)/(Sw_j,w_j)$%
  \STATE $x_{j+1}=x_j+\alpha_j w_j$%
  \STATE $r_{j+1}=r_j-\alpha_j Sw_j$%
  \STATE $\beta_j=(r_{j+1},r_{j+1})/(r_j,r_j)$
  \STATE $w_{j+1}=r_{j+1}+\beta_j w_j$%
\ENDFOR
\end{algorithmic}
\end{algorithm}
The algorithm is based on the construction of various basis of
$\mathcal{K}_m(S,r_0)$ : $(r_m)$ (residual basis) is orthogonal,
$(w_m)$ (research direction basis) is $S$-orthogonal. Step $6-7$
of algorithm \ref{alg:cg:1} is the $S$-orthogonalization of
$w_{j+1}$ with respect to $w_{j}$ which theoretically implies the
orthogonality of $w_{j+1}$ with respect to all previous research
directions. However this orthogonality property is numerically
lost as the number of iterations increases, it is then more suited
to use a full orthogonalization of research directions leading to
algorithm \ref{alg:cg:2}.

\begin{algorithm}\caption{Reorthogonalized conjugate gradient}\label{alg:cg:2}
\begin{algorithmic}[1]
\STATE Compute $r_0=b-Sx_0$, set $w_0=r_0$%
\FOR{$j=0,\ldots,m$} %
  \STATE $\alpha_j=(r_j,r_j)/(Sw_j,w_j)$%
  \STATE $x_{j+1}=x_j+\alpha_j w_j$%
  \STATE $r_{j+1}=r_j-\alpha_j Sw_j$%
  \STATE For $0\leqslant i\leqslant j$, $\beta_j^i=-(r_{j+1},Sw_{i})/(w_i,Sw_i)$
  \STATE $w_{j+1}=r_{j+1}+\sum_{i=1}^j {\beta_j^i w_i}$%
\ENDFOR
\end{algorithmic}
\end{algorithm}

Full reorthogonalization is often compulsory for complex
simulations. Various implementations are available (among others
Gram-Schmidt, modified Gram-Schmidt, iterative Gram-Schmidt
\cite{LINGEN:2000:EGS,FRAYSSE:1998:IOS}) depending on the chosen
ratio between precision and computational cost. Our experience
leads us to prefer modified Gram-Schmidt algorithm (the one used
in algorithm \ref{alg:gmres:1}) to classical Gram-Schmidt (the one
used in algorithm \ref{alg:cg:2}). Note that once fully
reorthogonalized, conjugate gradient is almost as expensive as
GMRes. Anyhow conjugate gradient provides the approximation at
each iteration, which can be very useful (see for instance section
\ref{subsub:error_dual}, where such an information enables the
computation of relevant convergence criterion).

\subsection{Study of the convergence, preconditioning}
Because of their error-minimization property, conjugate gradient
and GMRes have convergence theorems with known minimal convergence
rate, for instance for conjugate gradient:
\begin{equation}\label{eq:cg_cv:1}
    \|x-x_m\|_S \leqslant 2 \left[
    \frac{\sqrt{\kappa}-1}{\sqrt{\kappa}+1}\right]^m  \|x-x_0\|_S
\end{equation}
where $\kappa$ is the condition number of matrix $S$. Condition number is the ratio between the biggest and the smallest eigenvalues.
\begin{equation}\label{eq:cond}
    \kappa = \left| \frac{\lambda_n}{\lambda_1} \right| \text{\ with\ } |\lambda_1 |
    \leqslant |\lambda_2| \leqslant \ldots \leqslant |\lambda_n|
    \text{\ eigenvalues\ of\ }S
\end{equation}

Moreover performance results of Krylov iterative solvers are
strongly linked to the spectrum of matrix $S$. More precisely only
the active spectrum (set of eigenvalues which the right hand side
is not orthogonal to the associated eigenvectors) influences the
convergence; condition number $\kappa$ can be replaced with active
condition number $\kappa_{act}$ inside relation \eqref{eq:cg_cv:1}
leading to better convergence range. More precise study would lead
to the introduction of Ritz spectrum and effective condition
number \cite{VANDERSLUIS:1986:RCCG}.

These simple considerations are sufficient to explain the interest
of preconditioning the system: the idea is to solve equivalent
system $\widetilde{S}^{-1}Sx=\widetilde{S}^{-1}b$ where
$\widetilde{S}^{-1}$ is a well-chosen matrix providing the system
with better spectral properties (if $\widetilde{S}^{-1}\approx
S^{-1}$ then condition number is optimal, which justifies the
notation).

For conjugate gradient, the use of preconditioner may seem
problematic since the symmetry is \textit{a priori} lost. However
if preconditioner $\widetilde{S}^{-1}$ is symmetric definite
positive, applying conjugate gradient to nonsymmetric system is
equivalent to a symmetric resolution
($(L^{-T}SL^{-1})(Lx)=L^{-T}b$ with Cholesky factorization
$\widetilde{S}=L^TL$) and the method is still relevant.

So preconditioning the above two algorithms is simply realized
replacing $S$ by $\widetilde{S}^{-1}S$ and $b$ by
$\widetilde{S}^{-1}b$ in lines 1 and 3 of algorithm
\ref{alg:gmres:1}, and in lines 1 and 5 of algorithm
\ref{alg:cg:2} (anyhow the research directions are still
$S$-orthogonal). Of course the main problem remains the definition
of an efficient preconditioner.

\subsection{Constrained Krylov methods, projector implementation}\label{sub:conskry}
We may deal (for instance in the dual approach) with constrained systems such as:
\begin{equation}\label{eq:kryconst:1}
\begin{pmatrix}
S & G \\ G^T & 0
\end{pmatrix} \begin{pmatrix}
x \\ \alpha
\end{pmatrix} = \begin{pmatrix}
b \\ e \end{pmatrix}
\end{equation}
Because constraint $G^Tx_0=e$ is compulsory, it is often referred
to as "admissibility constraint". A classical solution is to find
an initialization $x_0$ which satisfies constraint and then ensure
that the remainder of the solution is researched inside a
supplemental space: $G^T(x_i-x_0)=0$. A projected algorithm
naturally arises:
\begin{eqnarray}
x & = & x_0 + Px^* \\
G^T x_0 & = & e \\
G^T P & = &0
\end{eqnarray}
which leads to:
\begin{eqnarray}
x_0 & = & QG\left(G^TQG\right)^{-1}e \\
P & = & I-QG\left(G^TQG\right)^{-1}G^T
\end{eqnarray}
where $Q$ is a matrix so that matrix $G^TQG$ is invertible. Iterative system then reads:
\begin{eqnarray}
P^T S P x^* & = & P^T\left(b-Sx_0\right)
\end{eqnarray}
$\alpha$ can be post-computed $\alpha=\left(G^TQG\right)^{-1}G^TQ\left(b-Sx\right)$.

\subsection{Augmented-Krylov methods, projector implementation}\label{sub:augmkry}
Augmented-Krylov methods \cite{CHAPMAN:1996:DAK,SAAD:1997:AAK} are
employed to add optional constraints to the resolution of a
system. The principle is to set subspace $\mathcal{C}$ in
$\mathbb{R}^n$ of dimension $n_c$ represented by $n\times n_c$
rectangular matrix $C$ ($\range(C)=\mathcal{C}$, for more
simplicity we suppose that  $C$ is full-ranked-column), then to
define augmented-Krylov subspace
$\widetilde{\mathcal{K}}_m(S,r_0,C)=\mathcal{K}_m(S,r_0)+\range(C)$,
and to use the following research principle:
\begin{equation}\label{eq:def_krylov_aug:1}
    \left\{ \begin{array}{l} x_m \in x_0 + \widetilde{\mathcal{K}}_m(S,r_0,C) \\
    r_m \perp_? \widetilde{\mathcal{K}}_m(S,r_0,C) \end{array} \right.
\end{equation}

Augmented-Krylov methods can be implemented either by
reorthogonalization schemes or by projection methods which are the
one we propose to present here. The research space is separated
into two subspaces: $\range(C)$ and a supplemental subspace. The
part of the solution in $\range(C)$ is detected during
initialization, while the remainder is iteratively looked for,
projector $P$ ensures the research is realized inside the correct
subspace.
\begin{eqnarray}\label{eq:krylov_aug:2}
    x&=&x_0 + P x^* \\
    C^T r_0 & = & C^T (b-S x_0) = 0 \\
    C^T S P & = & 0
\end{eqnarray}
Which leads to:
\begin{eqnarray}\label{eq:krylov_aug:3}
    x_0 & = &C \left(C^TSC\right)^{-1}C^T b \\
    P & = & I - C\left(C^TSC\right)^{-1}C^TS
\end{eqnarray}
system then reads:
\begin{eqnarray}\label{eq:krylov_aug:4}
    SPx^* & = & b-Sx_0 \\
\text{\ or\ }    P^TSP x^* & =& P^T \left(b-Sx_0\right)=P^T b
\end{eqnarray}

Though it can be proved that projected system is better
conditioned than original problem \cite{DOSTAL:1988:CGM}, the
efficiency of the method essentially depends on the choice of
matrix $C$, which is most often an opened problem.  Within the
framework of domain decomposition methods, this choice can be
guided by several considerations. In the framework of
multiresolution, reuse of previous numerical information can lead
to very interesting performance results
\cite{REY:1998:RKS,GOSSELET:2002:SRK,GOSSELET:2002:DDM,REY:2003:SLN,Rey:1998:RRP,Risler:2000:IAA}

\subsection{Constrained augmented Krylov methods}\label{sub:consaugmkry}
We here consider solving constrained system \eqref{eq:kryconst:1}
with $C$-augmented algorithm. The admissibility constraint is
often referred to as first level constraint and augmentation as
second level constraint. Two strategies are possible, the first
consists in mixing levels together while the second respects the
hierarchy between constraints.
\begin{description}
\item[One projector strategy:] set $J=\begin{pmatrix} G & S^T H \end{pmatrix}$ and $\tilde{e}=\begin{pmatrix} -e \\ H^T b \end{pmatrix}$, then system reads:
\begin{equation}
\begin{pmatrix}
S & J \\ J^T & 0
\end{pmatrix} \begin{pmatrix}
x \\ \tilde{\alpha}
\end{pmatrix} = \begin{pmatrix}
b \\ \tilde{e} \end{pmatrix}
\end{equation}
the following initialization/projection are employed ($Q$ is a parameter to tune)
\begin{eqnarray}
x_0&=&J\left(J^TQJ\right)^{-1}\tilde{e}\\ P&=&I-QJ\left(J^TQJ\right)^{-1}J^T
\end{eqnarray}Because $Q$ is not easy to interpret and choose, this method is hardly ever used.
\item[Two-projector strategy:] the two conditions are imbricated. First ensure admissibility constraint
\begin{eqnarray}
x&=&x_0+P x^*\\
x_0&=&G\left(G^TQG\right)^{-1}e\\P&=&I-QG\left(G^TQG\right)^{-1}G^T
\end{eqnarray}
then set
\begin{eqnarray}
x^*&=&x^*_0+P^*x^{**}\\
x^*_0&=&C\left(C^TP^TSPC\right)^{-1}C^TP^T\left(b-Sx_0\right)\\
P^*&=&I-PC\left(C^TP^TSPC\right)^{-1}C^TP^TS
\end{eqnarray}
so that optimality constraint is verified. As can be seen such an
approach is equivalent to classical augmentation with making
second level constraints consistent with the first (setting
$C^*=PC$).
\end{description}







\end{document}